
%
%
\input amstex
\documentstyle{amsppt}
\magnification = \magstep 1
\loadbold
\def\grad{\text{\rm grad}\,}
\def\sgrad{\text{\rm sgrad}\,}
\def\hess{\text{\rm Hess}\,}
\def\tr{\text{\rm tr}\,}
\def\ad{\text{\rm ad}\,}
\def\Exp{\text{\rm Exp}}
\def\pc{\preceq}
\def\pcc{\prec}
\def\ol{\overline}
\def\R{\boldkey R}
\def\C{\boldkey C}
\def\Z{\boldkey Z}
\def\Q{\boldkey Q}
\def\p{\frak p}
\def\n{\frak n}
\def\D{\Delta}
\def\wt{\widetilde}
\def\A{\Cal A}
\def\J{\Cal J}
\def\I{\Cal I}

\topmatter
\title On a class of K\"ahler manifolds whose
geodesic flows are integrable
\endtitle
\rightheadtext{}
\author Kazuyoshi Kiyohara\endauthor
\leftheadtext{}
\address Department of Mathematics, Faculty of Science,
Hokkaido University, Sapporo 060, Japan
\email kiyohara\@euler.math.hokudai.ac.jp\endemail
\endaddress
\abstract We study $n$-dimensional K\"ahler manifolds whose geodesic
flows
possess $n$ first integrals in involution that are
fibrewise hermitian forms and
simultaneously normalizable. Under some mild assumption, one can
associate
with such a manifold an $n$-dimensional commutative Lie
algebra of infinitesimal
automorphisms.  This, combined with the given $n$ first integrals,
makes the geodesic flow integrable.
If the manifold is
compact, then it becomes a toric variety.
\endabstract
\toc
\head {} Introduction
\endhead
\head {} Notations and preliminary remarks
\endhead
\head 1. Local calculus on $M^1$
\endhead
\head 2. Summing up the local data
\endhead
\head 3. Structure of $M-M^1$
\endhead
\head 4. Torus action and the invariant hypersurfaces
\endhead
\head 5. Properties as a toric variety
\endhead
\head 6. Bundle structure associated with a subset
of $\Cal A$
\endhead
\head 7. The case where $\#\Cal A=1$
\endhead
\head 8. Existence theorem
\endhead
\endtoc

\endtopmatter
\document

\specialhead Introduction
\endspecialhead
It is known that the geodesic flow of
the complex projective space $\C P^n$
equipped with the standard K\"ahler metric is (completely)
integrable in the
sense of symplectic geometry (or in Liouville's sense)
(cf. Thimm [13], see also
[9], [10], [5]).
The first integrals are given as follows: Let
$c_0,\dots,c_n$ be constants such
that $1=c_0>c_1>\dots>c_n=0$, and let $(z_0,\dots,z_n)$
be the homogeneous
coordinate system. Then, by putting $\partial_i=\partial/\partial z_i$,
$$\gather
\sum_{0\le j\le n\atop j\ne i}\frac{(z_i\bar\partial_j-z_j\bar\partial_i)
(\bar z_i\partial_j-\bar z_j\partial_i)}{c_j-c_i}\quad (1\le i\le n-1),\\
\sum_{i,j}(z_i\bar\partial_j-z_j\bar\partial_i)
(\bar z_i\partial_j-\bar z_j\partial_i),
\quad \sqrt{-1}(z_i\partial_i-\bar z_i
\bar\partial_i)\quad (1\le i\le n)
\endgather$$
are well-defined symmetric tensor fields and vector fields
on $\C P^n$. Regarded
as functions on the cotangent bundle $T^*\C P^n$,
they become a system of first
integrals in involution of the geodesic flow (note
that the first integrals given here
are slightly different from those in [13]).

In this paper, we shall define a class
of K\"ahler manifolds whose geodesic flows
are integrable by a set of first integrals
possessing similar properties to those
for $\C P^n$, and study the properties of
such manifolds. We shall call them
K\"ahler-Liouville manifolds (of type (A)).
The precise definition is as follows. Let $M$
be an $n$-dimensional complete
K\"ahler manifold, $I$ its complex structure, and $E$
its energy function
(the hamiltonian of the geodesic flow).
Let $\Cal F$ be an $n$-dimensional vector space
of sections of $S^2TM$
(the symmetric tensor product over $\R$ of two
copies of the tangent bundle $TM$).
Then we say that $(M,\Cal F)$
is K\"ahler-Liouville manifold if it satisfies the following
conditions:
\roster
\item"(1)" $E\in\Cal F$;
\item"(2)" $\{F,H\}=0$ for any $F,H\in\Cal F$;
\item"(3)" every $F\in\Cal F$ is ``hermitian'';
\item"(4)" $\Cal F_p=\{F_p\ |\ F\in\Cal F\}$ is simultaneously
normalizable
for every $p\in M$;
\item"(5)" $\Cal F_p$ is $n$-dimensional at some $p\in
M$.
\endroster
Here, $F_p\in S^2T_pM$ is the value of $F$
at $p\in M$. Since sections of $S^2TM$
are naturally regarded as functions on the cotangent
bundle $T^*M$, the Poisson
bracket in (2) makes sense. Also, (3) means
that the function $F$, restricted to
each fibre $T^*_pM$, is a hermitian form. We
say that two K\"ahler-Liouville
manifolds $(M,\Cal F)$ and $(M',\Cal F')$ are mutually
isomorphic if there is an
isomorphism $\phi:M\to M'$ of K\"ahler manifolds that maps
$\Cal F$ to $\Cal F'$.

As is immediately seen, only $n$ first integrals
are given in the definition.
Nevertheless, it will turn out that other $n$
first integrals appear
automatically if some non-degeneracy condition is assumed. A
K\"ahler-Liouville
manifold satisfying this assumption is called {\it of
type (A)} (for the precise
definition, see Section 1). The main purpose of
this paper is to investigate
local and global properties of K\"ahler-Liouville manifolds of
type (A).
Compact, 2-dimensional K\"ahler-Liouville manifolds were studied by Igarashi
[4], in which he adopted another type of
non-degeneracy condition and obtained
similar results to ours in this case. The
results indicate that the condition
adopted in [4] is equivalent to our condition
of type (A).

We now explain the various results in this
paper. In the following,
K\"ahler-Liouville manifolds are assumed to be of type
(A), unless otherwise
stated.

\proclaim{1 (cf. Proposition 1.9)} A finite, partially ordered
set $\A$ is naturally
associated with $(M,\Cal F)$. Also, a positive integer
$|\alpha|$ is assigned
to each $\alpha\in\A$ so that $\sum_{\alpha\in\A} |\alpha|=n$. For
any
$\alpha\in\A$, the subset $\{\beta\in\A\ |\ \beta\pcc\alpha\}$ is totally
ordered.
\endproclaim

\proclaim{2 (cf. Proposition 1.16, Theorems 3.1, 3.2)} An
$n$-dimensional
commutative Lie algebra
$\frak k$ of infinitesimal automorphisms of $(M,\Cal F)$
is naturally defined,
possessing the property that $\frak k$ and $\Cal
F$ are elementwise commutative
with respect to the Poisson bracket. With $\frak
k$ and $\Cal F$ the geodesic
flow of $M$ becomes integrable.
\endproclaim

Up to now, further results are obtained only
for compact K\"ahler-Liouville
manifolds. In this case, we obtain the results
below. Put $\frak g=\frak k+
I\frak k$, and let $K$ and $G$ be
the transformation group of $M$ generated by
$\frak k$ and $\frak g$ respectively.

\proclaim{3 (Theorems 4.2, 4.18)} $K$ and $G$ are
isomorphic to $U(1)^n$ and
$(\C^{\times})^n$ respectively. With the action of $G$, $M$
becomes a toric
variety.
\endproclaim

The structure of $M$ as toric variety is
completely investigated in Section 4,
and as a consequence we obtain the notion
of ``toric variety of KL-A type''
(cf. Section 5). The toric variety associated with
a K\"ahler-Liouville manifold
of type (A) is necessarily of KL-A type.
Conversely, we have the following result.

\proclaim{4 (cf. Theorem 8.3)} Suppose that a toric
variety of KL-A type is given.
Then there exists a K\"ahler-Liouville manifold of type
(A) such that the
associated toric variety is isomorphic to the given
one.
\endproclaim

Thanks to the general theory for toric varieties
(cf. [1], [2], [3], [11]), we can know
what kind of complex manifold $M$ is. For
the detail, see Section 5.
In particular, we have the following bundle structures.

\proclaim{5 (Proposition 5.4)} There is a holomorphic principal
fibre bundle
$$\prod_{\alpha\in\A}(\C^{|\alpha|+1}-\{0\})\to M$$
whose structure group is isomorphic to $(\C^{\times})^{\#\A}$.
\endproclaim

\proclaim{6 (Proposition 5.5, Theorems 6.3, 6.4, 6.5)} Let
$\A'$ be a subset of
$\A$ possessing the property that if $\alpha\in\A'$ and
$\beta\pcc\alpha$, then
$\beta\in\A'$. Put $\A''=\A-\A'$. Then, associated with $\A'$, there
naturally
exist K\"ahler-Liouville manifolds $(M',\Cal F')$, $(M'',\Cal F'')$, and
a holomorphic
fibre bundle $\pi: M\to M''$ whose typical fibre
is $M'$. They possess the following
properties: (1) $(M'',\Cal F'')$ is of type (A)
and the associated partially
ordered set is $\A''$; (2) if $(M',\Cal F')$
is of type (A), then the associated
partially ordered set is $\A'$; (3) there is
a homomorphism $G\to G''$ so that
$\pi: M\to M''$ is equivariant, where $G$ and
$G''$ denote the algebraic tori acting
on $M$ and $M''$ respectively; (4) even if
$(M',\Cal F')$ is not of type (A), $M'$
possess the structure of toric variety inherited from
$M$ so that the maximal
compact subgroup of the algebraic torus acting on
$M'$ preserves the metric and
each element of $\Cal F'$.
\endproclaim

The property (4) stated above indicates that the
geodesic flow of $(M',\Cal F')$
is integrable even if it is not of
type (A). In this case we shall say that the
K\"ahler-Liouville manifold $(M',\Cal F')$ is of type (B)
(cf. Section 6).
Such a manifold will be necessary for the
study of K\"ahler-Liouville manifold
of type (A) only when $\dim M'=1$.
Using the result above successively, we consequently obtain
a family of
K\"ahler-Liouville manifolds $(M_{\alpha},
\Cal F_{\alpha})$ $(\alpha\in\A)$ such
that the partially ordered set associated with $(M_{\alpha},\Cal
F_{\alpha})$
consists of one element $\{\alpha\}$.
In this case, it turns out that $M_{\alpha}$
is isomorphic to $\C P^{|\alpha|}$ as
toric variety. It also turns out that it
is of type (A) if and only if $|\alpha|\ge 2$.

The result 4 mentioned above is actually given
in much finer form in Theorem 8.3.
There, besides the toric variety $M$, we prescribe
the K\"ahler-Liouville manifolds
$(M_{\alpha},\Cal F_{\alpha})$ $(\alpha\in\A)$ and some constants, and prove
the uniqueness of $(M,\Cal F)$ as well as
the existence. The reason for
formulating the ``existence theorem'' in this form is
that K\"ahler-Liouville
manifolds such that the associated partially ordered sets
consist of one element
are easily understandable by using our work [6]
on (real) Liouville
manifolds. The result is roughly stated as follows.

\proclaim{7 (cf. Theorem 7.2)} The isomorphism classes of
K\"ahler-Liouville
manifolds such that $\#\A=1$ are completely classified by
means of several
constants and a function on a circle.
\endproclaim

We now briefly explain the organization of the
paper. In Section 1 we first
formulate some non-degeneracy condition (depending on points) on
a
K\"ahler-Liouville manifold $(M,\Cal F)$. Denoting by $M^1$ the
set of points
where the condition is satisfied, we say that
$(M,\Cal F)$ is of type (A) if
$M^1\ne\emptyset$. We perform local calculus on $M^1$, and
introduce almost all
basic quantities related to $(M,\Cal F)$.
In Sections 2 and 3 we sum up
the local data given in Section 1, and describe
properties of the basic quantities and the Lie
algebra $\frak k$ in global form.
The result 2 stated above is proved in
Section 3.

Through Sections 4--8 we assume that $M$ is
compact. In Section 4 we prove the
result 3 stated above and the results that
determine the structure of $M$ as
toric variety. In particular, we show that $M^1$
is the unique open $G$-orbit,
and $M-M^1$ is the union of $n+\#\A$ closed
hypersurfaces that are $G$-invariant
and totally geodesic.

In Section 5 we describe the various properties
of the toric variety $M$.
There we specify the fan of $M$, and
give the definition of toric variety of KL-A
type. This section contains 3 subsections;
``The fan of $M$'', ``Fibre bundles associated with
$M$'', and ``Line bundles''.
In Section 6 we prove the result 5
stated above. There we also prove its converse
(Theorem 6.11), which plays a crucial role in
the proof of Theorem 8.3.

Section 7 is devoted to the proof of
Theorem 7.2 (see the result 7 above). We
establish the one-to-one correspondence between the isomorphism classes
of K\"ahler-Liouville manifolds of type (A) such that
$\#\A=1$ and the
isomorphism classes of Liouville manifolds of rank one
and type (B) that satisfy
a certain condition. The definition
and the classification of Liouville manifolds of rank
one are given in [6].
By using them, the theorem is proved.
In Section 8 we prove Theorem 8.3 mentioned
above. It is no longer hard by
virtue of Theorems 6.11 and 7.2.

\specialhead Notations and preliminary remarks
\endspecialhead
Let $M$ be an $n$-dimensional K\"ahler manifold, and
let $g$ and $I$ be its
K\"ahler metric and complex structure respectively. Then the
K\"ahler form
$\omega$ and the energy function $E$ are given
as follows:
$$\omega(X,Y)=g(IX,Y) \quad (X,Y\in T_pM,p\in M)\qquad
E(\lambda)=\frac12 |\lambda|^2 \quad (\lambda\in T^*M),$$
where $|\cdot|$ denotes the norm function on $T^*M$
associated with the metric
$g$. The energy function $E$ is the hamiltonian
of the geodesic flow of $M$.
For a function $h$ on $M$ we define
vector fields $\grad h$ and $\sgrad h$ by
the following formulae:
$$i_{\grad h}g=dh,\qquad i_{\sgrad h}\omega=-dh,$$
where $i$ denotes the inner derivation. We have
$\sgrad h=I\grad h$.

Let $p\in M$, and let $S^2T_pM$ be the
symmetric tensor product of two copies
of the tangent space $T_pM$.
Let $F$ be an element of $S^2T_pM$, and
suppose that, regarded as a
quadratic form on $T^*_pM$, $F$ is a hermitian
form. Then there is an orthonormal
basis $V_j,IV_j$ $(1\le j\le n)$ of $T_pM$ and
constants $a_1,\dots,a_n$ such that
$$F=\sum_ja_j(V_j^2+(IV_j)^2).$$
We define the endomorphism $F^e$ of $T_pM$ by
$$F^e(V_i)=a_iV_i,\quad F^e(IV_i)=a_iIV_i.$$
Clearly it is independent of the choice of
$\{V_i\}$. Regarding
$T_pM$ as a complex vector space (by identifying
$I$ with $\sqrt{-1}$), $F^e$
becomes $\C$-linear. We define $\tr F$ and $\det
F$ by the trace and the
determinant of $F^e$ over $\C$ respectively ($\tr F=\sum_ja_j$,
$\det F=\prod_ja_j$).

Let $(M,\Cal F)$ be a K\"ahler-Liouville manifold. Then
the condition (4) in the
definition is equivalent to that $\{F_p^e\ |\ F\in\Cal
F\}$ is commutative with
respect to the composition of endomorphisms for every
$p\in M$.

We shall use the term ``smooth'' in the
same meaning as ``of class $C^{\infty}$''.

\specialhead 1. Local calculus on $M^1$
\endspecialhead
Let $(M, \Cal F)$ be a K\"ahler-Liouville manifold
of dimension $n$. Put
$$M^0=\{p\in M\ |\ \dim \Cal F_p=n\}\qquad M^s=M-M^0.$$
Then $M^0$ is an open subset of $M$,
which is not empty because of the condition
(5) in the definition of K\"ahler-Liouville manifold. Let
$F_1,\dots,F_n$ be a
basis of $\Cal F$, and let $p$ be
a point of $M^0$. Then there
are an orthonormal frame $V_i, IV_i$ $(i=1,\dots, n)$
and $n^2$ functions $f_{ij}$
around $p$ such that
$$F_i=\sum_{j=1}^n f_{ij} \ (V_j^2+ (IV_j)^2).$$
Putting $(a_{ij})=(f_{ij})^{-1}$, we have
$$\sum_{j=1}^n a_{ij} F_j= V_i^2+ (IV_i)^2.$$
Let $D_i$ be the subbundle of $TM$ defined
around $p$ spanned by $V_i$ and
$IV_i$.

\proclaim{Proposition 1.1} There are positive functions
$a_1,\dots,a_n$ around
the point $p$ such that, putting
$$b_{ij}=\frac{a_{ij}}{a_i},\qquad W_i=\frac{V_i}{\sqrt{a_i}},$$
$$W_jb_{ik}=(IW_j)b_{ik}=0\quad (j\ne i,\ \text{any } k),\tag 1.1$$
$$\{W_i^2+(IW_i)^2,W_j^2+(IW_j)^2\}=0\quad (\text{any }i,j).\tag 1.2$$
The function $a_i$ can be chosen to be
one of $|a_{i1}|,\dots,|a_{in}|$ that is
non-zero around $p$. Moreover, if $\{a'_i\}$ possess the
same properties as
above, then
$$V_j\frac{a'_i}{a_i}= (IV_j)\frac{a'_i}{a_i}=0\quad (j\ne i).$$
\endproclaim
\demo{Proof} We have
$$\align
& \{V_i^2+(IV_i)^2,V_j^2+(IV_j)^2\}=
\sum_{k,l=1}^n\left(\{a_{ik},F_l\}a_{jl}F_k
+\{F_k,a_{jl}\}a_{ik}F_l\right)\\
& =\sum_k\{a_{ik},V_j^2+(IV_j)^2\}F_k+\sum_l\{V_i^2+(IV_i)^2,a_{jl}\}F_l
\tag 1.3 \endalign$$
Note that each term in the formula above
is a homogeneous polynomial of degree 3
in the variables $V_k, IV_k$ $(1\le k\le n)$.
Since the left-hand side belongs to
the ideal $(V_iV_j,V_iIV_j,V_jIV_i,IV_iIV_j)$ of the polynomial algebra, and
since $F_k$ are linear combinations of $V_l^2+(IV_l)^2$ $(1\le
l\le n)$, it follows
that
$$\sum_k\{a_{ik},V_j\}F_k=c_{ij}(V_i^2+(IV_i)^2),\quad
\sum_k\{a_{ik},IV_j\}F_k=d_{ij}(V_i^2+(IV_i)^2)\tag 1.4$$
for some functions $c_{ij}$ and $d_{ij}$, provided $i\ne
j$. Hence we have
$$\{a_{ik},V_j\}=c_{ij}a_{ik},\quad \{a_{ik},IV_j\}=d_{ij}a_{ik}.\tag 1.5$$

Let $a_i$ be one of the functions $|a_{i1}|,\dots,|a_{in}|$
that does not vanish
around the point $p$. Then by (1.5) we
have
$$\left\{\frac{a_{ik}}{a_i}, V_j\right\}=
\left\{\frac{a_{ik}}{a_i}, IV_j\right\}=0
\qquad (i\ne j).$$
This implies that
$$\sum_k\{a_{ik},V_j^2+(IV_j)^2\}F_k=\frac1{a_i}\{a_i,V_j^2+(IV_j)^2\}
(V_i^2+(IV_i)^2).$$
Hence, by (1.3) we obtain
$$\{\frac1{a_i}(V_i^2+(IV_i)^2),\frac1{a_j}(V_j^2+(IV_j)^2)\}=0.$$
The remaining part is clear.
\qed
\enddemo

\proclaim{Proposition 1.2} $[V_i,IV_i]\equiv\sgrad (\log a_i)\quad \mod D_i$.
\endproclaim
\demo{Proof} We use the K\"ahler form $\omega$. We
have
$$0=d\omega(V_i,IV_i,V_j)=-\omega([V_i,IV_i],V_j)+\omega([V_i,V_j],IV_i)
-\omega([IV_i,V_j],V_i).$$
Since (1.2) implies
$$\aligned
[W_i,W_j]=\alpha_{ij}IW_i-\alpha_{ji}IW_j,\quad
& [IW_i,IW_j]=\beta_{ij}W_i-\beta_{ji}W_j,\\
[W_i,IW_j]=-\beta_{ij}IW_i+\alpha_{ji}W_j,\quad
& [IW_i,W_j]=-\alpha_{ij}W_i+\beta_{ji}IW_j
\endaligned \tag 1.6$$
for some functions $\alpha_{ij}$ and $\beta_{ij}$ $(i\ne j)$,
it follows that
$$\omega([V_i,V_j],IV_i)=-\omega([IV_i,V_j],V_i)=-V_j\log \sqrt{a_i}.$$
Hence we have
$$\omega([V_i,IV_i],V_j)=-V_j\log a_i, \tag 1.7$$
provided $j\ne i$. Replacing $V_j$ with $IV_j$, we
also have
$$\omega([V_i,IV_i],IV_j)=-IV_j\log a_i.\tag 1.8$$
{}From (1.7) and (1.8) it thus follows that
$$[V_i,IV_i]\equiv \sum_{j\ne i}
\left(-(IV_j\log a_i)V_j+(V_j\log a_i)IV_j\right)
\quad \mod D_i.$$
\qed
\enddemo

We now consider the following condition for points
on $M^0$:
$$\text{For any }i\text{ there is some }j \text{
such that }da_j|_{D_i}\ne 0
\text{ at }p.\tag 1.9$$
Note that this condition is independent of the
choice of $\{a_i\}$. Put
$$M^1=\{p\in M^0\ |\ (1.9)\text{ holds at }p\}.$$
We shall say that a K\"ahler-Liouville manifold $(M,\Cal
F)$ is {\it of type (A)} if
$$M^1\ne \emptyset.$$
{}From now on (until the end of Section
4) K\"ahler-Liouville manifolds are assumed
to be of type (A), unless otherwise stated.

Let $p\in M^1$.

\proclaim{Proposition 1.3} Let $j,k\ne i$. If $d\log a_j|_{D_i}\ne
0$ and
$d\log a_k|_{D_i}\ne 0$ at $q\in M^1$ near $p$,
then they are linearly dependent
at $q$.
\endproclaim
\demo{Proof} Note first that $\sum_j a_{jk}$ is constant
for every $k$, because
$E\in\Cal F$. Since $a_{jk}=a_jb_{jk}$, this implies that $\{a_i\}$
and $\{a_{ij}\}$
are written as rational functions of $\{b_{ij}\}$. Hence
for each $i$ there is
some $l$ such that $db_{il}\ne 0$ around $p$.
On the other hand, the kernel of
$d\log a_j$ on $D_i$ is spanned by
$$-(IV_i\log a_j)V_i+(V_i\log a_j)IV_i=[V_j,IV_j]_{D_i},$$
where the right-hand side denotes the $D_i$-component of
$[V_j,IV_j]$.
Since $V_jb_{il}=IV_jb_{il}=0$, we also have
$$0=[V_j,IV_j]b_{il}=[V_j,IV_j]_{D_i}b_{il}.$$
Hence the kernel of $d\log a_j$ on $D_i$
coincides with that of $db_{il}$ on $D_i$.
Since the latter does not depend on $j$,
the proposition follows.
\qed
\enddemo

By virtue of the proposition above we can
take the orthonormal frame $V_i,IV_i$
$(i=1,\dots,n)$ around $p\in M^1$ so that
$$d\log a_j(IV_i)=0\qquad \text{for any }j\ne i.$$
Note that $V_i$ are uniquely determined up to
sign (and the numbering).
We shall assume that $V_i$ are taken in
this way. Let $D^+$ (resp. $D^-$) be
the subbundle of $TM$ spanned by $V_1,\dots,V_n$ (resp.
$IV_1,\dots,IV_n$). $D^+$ and $D^-$ are well-defined over $M^1$;
$$TM^1=D^+\oplus D^-.$$

\proclaim{Proposition 1.4} \roster
\runinitem"(1)" $da_i,da_{ij},db_{ij}$ are zero on $D^-$.
\item"(2)" For any $i$ there is some $j$
such that $W_ib_{ij}\ne 0$.
\item"(3)" $[W_i,W_j]=[IW_i,IW_j]=[W_i,IW_j]=0$ $(i\ne j)$. In particular,
$D^+$ and $D^-$ are integrable.
\endroster
\endproclaim
\demo{Proof} (1) and (2) are clear from the
proof of Proposition 1.3. Suppose that
$W_ib_{ik}\ne 0$. Then by (1.6),
$$0=[IW_i,IW_j]b_{ik}=\beta_{ij}W_ib_{ik},\quad
0=[IW_i,W_j]b_{ik}=-\alpha_{ij}W_ib_{ik}.$$
Hence $\alpha_{ij}=\beta_{ij}=0$, and (3) follows.
\qed
\enddemo

\proclaim{Proposition 1.5} $[W_i,IW_i]\in D^-$.
\endproclaim
\demo{Proof} By virtue of Propositions 1.2 and 1.4,
it suffices to prove that the
$D_i$-component of $[W_i,IW_i]$ belongs to $D^-$. Choose $j(\ne
i)$ such that
$d\log a_j|_{D_i}\ne 0$. Then the $D_i$-component of $[W_j,IW_j]$
is not zero.
Describing
$$[W_j,IW_j]= \alpha W_j +\beta IW_j + \sum_{k\ne j}\gamma_k
IW_k,$$
we have
$$0=[W_i,[W_j,IW_j]]_{D_i}=(W_i\gamma_i)IW_i +\gamma_i[W_i,IW_i]_{D_i}.$$
Since $\gamma_i\ne 0$, the proposition follows.
\qed
\enddemo

\proclaim{Proposition 1.6} Maximal integral manifolds of $D^+$ are
(locally)
totally geodesic.
\endproclaim
\demo{Proof} Since $IW_k<W_i,W_j>=0$ and $[W_i,IW_k]\in D^-$,
it follows that
$$<\nabla_{W_i}W_j, IW_k>=0$$
for any $i,j,k$, where $\nabla$ denotes the Levi-Civita
covariant derivative. Hence
the proposition follows.
\qed
\enddemo

By virtue of Proposition 1.5 the $D_i$-component of
$[W_i,IW_i]$ is of the form
$c_iIW_i$, $c_i$ being the function around $p$.

\proclaim{Proposition 1.7} \roster
\runinitem"(1)" For $i,j$ such that $i\ne j$ and
$W_ia_j\ne 0$,
$$c_i=-W_i\log a_i +W_i\log a_j- \frac{W_i^2\log a_j}{W_i\log a_j}.$$
\item"(2)" $(IW_k)c_i=0\qquad (\text{any } k)$.
\item"(3)" $W_kc_i= -(W_k\log a_i)(W_i\log a_k)\qquad (k\ne i)$.
\item"(4)" $W_jW_i\log a_i= (W_i\log a_j)(W_j\log a_i)\qquad (i\ne j).$
\item"(5)" $W_iW_j\log a_k=0\qquad (i\ne j\ne k\ne i).$
\item"(6)" $(W_i\log a_k)(W_j\log a_k)= (W_i\log a_j)(W_j\log a_k)
+(W_j\log a_i)(W_i\log a_k)$ $(i\ne j\ne k\ne i)$.
\endroster
\endproclaim
\demo{Proof} By Propositions 1.2 and 1.4 we have
$$[W_i,IW_i]=c_iIW_i+\sum_{j\ne i}\frac{a_j}{a_i}(W_j\log a_i)IW_j.$$
Then, computing $[W_k,[W_i,IW_i]]$ for $k\ne i$, we obtain
$$\align
0 & =(W_kc_i+(W_k\log a_i)(W_i\log a_k))IW_i\\
& + (c_k\frac{a_k}{a_i} W_k\log a_i+W_k(\frac{a_k}{a_i}W_k\log a_i))IW_k\\
& + \sum_{j\ne i,k}\left(W_k(\frac{a_j}{a_i}W_j\log a_i)+
\frac{a_j}{a_i}(W_k\log a_i)(W_j\log a_k)\right)IW_j.
\endalign$$
This formula implies (1), (3), and
$$\align
& W_kW_j\log a_i =(W_k\log a_i)(W_j\log a_i) \tag 1.10
\\
& -(W_k\log a_j)(W_j\log a_i) -(W_j\log a_k)(W_k \log a_i)
\endalign$$
for mutually distinct $i,j,k$. Also, (2) follows from
(1) and Proposition 1.4.
To prove (4), (5), and (6), we recall
that $\sum_i a_ib_{il}$ are constants.
Differentiating these by $W_j$ and $W_k$ successively $(j\ne
k)$
we have
$$\align
& \sum_i(W_j\log a_i)a_ib_{il} + a_jW_jb_{jl}=0,\\
& \sum_i(W_kW_j\log a_i)a_ib_{il}+(W_j\log a_k)a_kW_kb_{kl}\\
& +\sum_i(W_j\log a_i)(W_k\log a_i)a_ib_{il}+ (W_k\log a_j)a_jW_jb_{jl}=0.
\endalign$$
Thus,
$$\align
& W_kW_j\log a_i =(W_j\log a_k)(W_k\log a_i)\tag 1.11\\
& +(W_k\log a_j)(W_j\log a_k) -(W_k\log a_i)(W_j\log a_i).
\endalign$$
{}From (1.10) and (1.11) the formulae (5) and
(6) follows. Also, since $i$ is
arbitrary in the formula (1.11), (4) follows by
putting $i=k$ in (1.11).
\qed
\enddemo

We now define the binary relation $\pc$ on
the set of indices $\{1,\dots,n\}$:
$$i\pc j\Longleftrightarrow i\ne j\text{ and }\ W_i\log a_j\ne
0
\text{ at }p\in M^1,\text{ or } i=j.$$
When it is necessary to clarify the point-dependence,
we shall say that $i\pc j$
at $p$. Also, we write $i\sim j$ if
$i\pc j$ and $j\pc i$.

\proclaim{Lemma 1.8}
\roster
\runinitem"(1)" If $i\pc j$ and $j\pc k$, then
$i\pc k$.
\item"(2)" The relation $\sim$ is an equivalence relation.
\endroster
\endproclaim
\demo{Proof} Assume that $i\not\pc k$ and $j\ne i,k$.
Then, by
Proposition 1.7 (6) we have
$$(W_i\log a_j)(W_j\log a_k)=0.$$
Hence $i\not\pc j$ or $j\not\pc k$, and (1)
follows. (2) is an immediate
consequence of (1).
\qed
\enddemo

Let $\A$ be the set of the equivalence
classes. It is clear from Lemma 1.8 (1)
that the relation $\pc$ induces a binary relation
(denoted by the same symbol)
on the set $\A$. For $\alpha\in \A$, let
$|\alpha|$ denote the number of indices
contained in the equivalence class $\alpha$.

\proclaim{Proposition 1.9} Let $\alpha,\beta,\gamma\in \A$.
\roster
\item"(1)" If $\alpha\pc\beta$ and $\beta\pc\gamma$, then $\alpha\pc\gamma$.
Namely, the relation $\pc$ is a partial order
on $\A$.
\item"(2)" If $\alpha\pc\gamma$ and $\beta\pc\gamma$, then $\alpha\pc\beta$
or $\beta\pc\alpha$. Namely, for any fixed $\gamma$, the
set of $\delta\in\A$
such that $\delta\pc\gamma$ is a totally ordered subset
of $\A$.
\item"(3)" If $\alpha$ is a maximal element, then
$|\alpha|\ge 2$.
\endroster
\endproclaim
\demo{Proof} (1) is clear from Lemma 1.8 (1).

(2) Let $i\in\alpha$, $j\in\beta$, and $k\in\gamma$, and assume
that
$\alpha\not\pc\beta$ and $\beta\not\pc\alpha$. Then $W_i\log a_j=W_j\log a_i
=0$. Hence by Proposition 1.7 (6), we have
either $W_i\log a_k=0$ or
$W_j\log a_i=0$.

(3) Suppose $\alpha$ is maximal and $i\in\alpha$. Then
by the condition (1.9)
we see that there is some $j(\ne i)$
such that $W_i\log a_j\ne 0$. Since
$\alpha$ is maximal, it follows that $j\in\alpha$.
\qed
\enddemo

Let $W_i^*$ $(1\le i\le n)$ be 1-forms such
that $W_i^*(IW_j)=0$
and $W_i^*(W_j)=\delta_{ij}$. Then Propositions 1.4 and 1.5 imply
that $dW_i^*=0$.
Hence there is a system of functions $(x_1,\dots,x_n)$
such that $dx_i=W_i^*$.
Clearly, $a_i$ are functions of $(x_1,\dots,x_n)$, and $W_j\log
a_i$ is nothing but
the derivative of $\log a_i$ with respect to
the variable $x_j$. For simplicity,
we shall write $\partial_j$ instead of $\partial/\partial x_j$.

\proclaim{Lemma 1.10} For any $i,j$ $(i\ne j)$, there
are functions $h_{ij}(x_i)$
such that $h_{ij}-h_{ji}\ne 0$ at $p$, and
$$\partial_i\log a_j=-\partial_i\log |h_{ij}-h_{ji}|.$$
Moreover, if $\partial_i\log a_j\ne 0$ and $\partial_i\log a_k\ne
0$ at $p$
$(j,k\ne i)$, then there are constants $c(\ne 0)$
and $d$ such that
$$h_{ik}=ch_{ij}+d.$$
\endproclaim
\demo{Proof} By Proposition 1.7 (5) the function $\partial_i\log
a_j$ depends
only on the variables $x_i$ and $x_j$. Since
$$\partial_j\partial_i\log a_j=\partial_i\partial_j\log a_i=
(\partial_i\log a_j)(\partial_j\log a_i),$$
it follows that there is a positive function
$H=H(x_i,x_j)$ such that
$$\partial_i\log a_j=\partial_i\log H,
\quad \partial_j\log a_i=\partial_j\log H,
\quad \partial_i\partial_j\log H=(\partial_i\log H)(\partial_j\log H).$$
The last formula implies that
$$\partial_i\partial_j\left(\frac1H\right)=0.$$
Hence there are functions $h_{ij}(x_i)$ and $h_{ji}(x_j)$ such
that $H^{-1}=
|h_{ij}-h_{ji}|$.

Now, let us assume that $\partial_i\log a_j\ne 0$
and $\partial_i\log a_k\ne 0$
at $p$. Then the derivatives $h'_{ij}$ and $h'_{ik}$
do not vanish at $p$. By
virtue of Proposition 1.7 (1), the function
$$\partial_i\log a_j-\frac{\partial_i^2\log a_j}{\partial_i\log a_j}=
-\frac{h''_{ij}}{h'_{ij}}$$
does not depend on $j$. Hence it follows
that
$$\frac{h''_{ij}}{h'_{ij}}=\frac{h''_{ik}}{h'_{ik}},$$
which proves the latter half of the lemma.
\qed
\enddemo

\proclaim{Proposition 1.11} For any $p\in M^1$ there is
a neighborhood $U$ of $p$
such that the relation $\pc$ is stable on
$U$. Namely,
$$i\pc j \text{ at } p\Longleftrightarrow i\pc j\text{
at } q$$
for any $q\in U$.
\endproclaim
\demo{Proof} Let $U(\subset M^1)$ be a connected neighborhood
of $p$ such that
every $\partial_i\log a_j$ $(i\ne j)$ that does not
vanish at $p$ also does
not vanish everywhere on $U$, and that all
the functions $h_{ij}$ are defined there.

Now, let us assume that $i\not\pc j$ at
$p$ and $i\pc j$ at some $q\in U$.
Let $k(\ne i)$ be a number such that
$i\pc k$ at $p$. Then $i\pc k$ at $q$, and by
Lemma 1.10,
$$h'_{ij}=ch'_{ik}\tag 1.12$$
on a neighborhood of $q$ for some constant
$c\ne 0$. Let $U'$ be the
set of all $q'\in U$ such that the
formula (1.12) is valid at $q'$. Clearly, $U'$
is closed in $U$. Moreover, if $q'\in U'$,
then $h'_{ij}\ne 0$ at $q'$. This implies
$i\pc j$ at $q'$ by Lemma 1.10, and
thus (1.12) holds on a neighborhood of $q'$.
Hence $U'$ is open. Since $q\in U'$ and
$U$ is connected, it follows that $U'=U$.
But since $p\not\in U'$, it is a contradiction.
\qed
\enddemo

We shall write $\alpha\pcc\beta$ $(\alpha,\beta\in \A)$ if $\alpha\pc\beta$
and $\alpha\ne \beta$.

\proclaim{Proposition 1.12} Let $p\in M^1$, and let $U$
be a small neighborhood
of $p$. Then there are functions $h_i$ $(1\le
i\le n)$ on $U$ and constants
$e_{\alpha\beta}$ $(\alpha,\beta\in\A$, $\alpha\pcc\beta)$ satisfying the
following four conditions:
\roster
\item"(1)" $dh_i|_{D^-}=dh_i|_{D_j}=0\quad (i\ne j),\qquad W_ih_i\ne 0\quad
\text{everywhere on }U$;
\item"(2)" $e_{\alpha\beta}=e_{\alpha\gamma}\qquad (\alpha\pcc\beta\pc
\gamma)$;
\item"(3)" $h_i\ne h_j\quad (i\sim j,\ i\ne j),\quad h_i+e_{\alpha\beta}\ne
0
\quad (i\in\alpha,\alpha\pcc\beta)\quad\text{everywhere on }U$;
\item"(4)" the functions $a_i$ can be taken in
the form
$$a_i=\left| \prod_{j\in\alpha\atop j\ne i}(h_j-h_i)\right|^{-1}
\left| \prod_{\gamma\pcc\alpha}\prod_{k\in\gamma}(h_k+e_{\gamma\alpha})
\right|^{-1} \quad (i\in\alpha).$$
\endroster
If $\{\wt h_i\}$ and $\{\wt e_{\alpha\beta}\}$ also satisfy
the conditions above,
then there are constants $c_{\alpha}\ne 0$ and $d_{\alpha}$
such that
$$\wt h_i=c_{\alpha}h_i-d_{\alpha}, \quad \wt e_{\alpha\beta}=c_{\alpha}
e_{\alpha\beta}+d_{\alpha}\qquad (i\in\alpha).$$
\endproclaim
\demo{Proof} We determine the functions $h_i$ and the
constants
$e_{\alpha\beta}$ inductively with respect to the partial order
$\pcc$. Let
$h_{ij}(x_i)$ $(i\ne j)$ and $U$ be as in
Lemma 1.10 and Proposition 1.11
respectively.

Let $\alpha\in\A$ be a minimal element. Then the
functions $a_i$ $(i\in\alpha)$
depend only on the variables $x_j$ $(j\in\alpha)$. Hence,
if $|\alpha|=1$, then
we can put $a_i=1$ (cf. Proposition 1.1). In
this case the function $h_i$ does not
appear yet. Also, if $|\alpha|=2$, then putting
$h_i=h_{ij}$ and $h_j=h_{ji}$ $(i,j\in\alpha)$, we can define $a_i$
and $a_j$ as
$$a_i=a_j=\frac1{|h_i-h_j|}.$$
Now, suppose that $|\alpha|\ge 3$.
By Lemma 1.10 we know that the ratios
$h'_{ij}/h'_{ik}$ are non-zero constants
for any mutually distinct $i,j,k\in\alpha$. Hence there is
a function $h_i(x_i)$
and non-zero constants $c_{ij}$ such that
$$h'_{ij}=c_{ij}h'_i$$
for any $i,j\in\alpha$ $(i\ne j)$. Then by Proposition
1.7 (6) we have
$$c_{ij}c_{jk}c_{ki}=c_{ji}c_{kj}c_{ik}\qquad (i,j,k\in\alpha)$$
Namely, $\{c_{ij}/c_{ji}\}$ satisfies the cocycle condition. An easy
calculation
shows that it is actually a coboundary, i.e.,
there are non-zero constants $c_i$
$(i\in\alpha)$ such that $c_{ij}/c_{ji}=c_i/c_j$. Then, putting
$$\wt h_{ij}=\frac{c_i}{c_{ij}} h_{ij},\qquad \wt h_i=c_ih_i,$$
we obtain
$$\partial_i\wt h_{ij}=\partial_i\wt h_i,\qquad \wt h_{ij}-\wt h_{ji}
=\frac{c_i}{c_{ij}}(h_{ij}-h_{ji}).$$

Putting $\wt h_{ij}-\wt h_i=d_{ij}$ (a constant), and using
Proposition 1.7 (6)
again, we have
$$d_{ij} +d_{jk}+d_{ki}=d_{ji}+ d_{kj}+ d_{ik}.$$
Put
$$d_i=|\alpha|^{-1}\sum_{j\in\alpha\atop j\ne i}
(d_{ij}-d_{ji})\qquad (i\in\alpha).$$
Then we have
$$d_{ij}-d_{ji}=d_i-d_j,\quad \wt h_{ij}-\wt h_{ji}
=(\wt h_i+d_i)-(\wt h_j+d_j)$$
for any $i,j\in\alpha$ $(i\ne j)$. Hence, by redefining
$h_i$ as
$$h_i=\wt h_i+d_i\qquad (i\in\alpha),$$
the difference
$$\log a_i-\sum_{j\in\alpha\atop j\ne i}\log |h_j-h_i|\qquad (i\in\alpha)$$
becomes a function of the single variable $x_i$.
Thus we can take $a_i$
$(i\in\alpha)$ as
$$a_i=\prod_{j\in\alpha\atop j\ne i}|h_j-h_i|.$$

Now, let $\alpha$ be a non-minimal element of
$\A$, and let $\beta\in\A$
$(\beta\pcc\alpha)$ be a unique element such that
there is no $\gamma\in\A$ satisfying $\beta\pcc\gamma
\pcc\alpha$. We assume that the functions $h_i$ are
defined for all $i\in\gamma$
and $\gamma\pcc\alpha$ ($\gamma\pcc\beta$ if $|\beta|=1$) and the constants
$e_{\gamma\delta}$ are defined for all $\gamma$ and $\delta$
$(\gamma\pcc
\delta\pcc\alpha)$ so that the conditions (1), (2), (3),
and (4) in the proposition
are satisfied. Under this assumption we shall define
suitable functions $h_i$
$(i\in\alpha)$ when $|\alpha|\ge 2$, $h_j$ $(j\in\beta)$ when $|\beta|=1$,
and a
constant $e_{\beta\alpha}$.

Let $i\in\alpha$, $j\in\beta$, and
$k\in\gamma(\pc\beta)$. Then the assumption
implies that
$$\partial_k\log a_j=\cases
-\partial_k\log |h_k-h_j|\quad & (\gamma=\beta,\ k\ne j)\\
-\partial_k\log |h_k+e_{\gamma\beta}|\quad & (\gamma\pcc\beta).
\endcases$$
Since
$$(\partial_j\log a_i)(\partial_k\log a_i)
=(\partial_j\log a_k)(\partial_k\log a_i)
+(\partial_k\log a_j)(\partial_j\log a_i),$$
we have $\partial_k\log a_i=\partial_k\log a_j$
in case $\gamma\pcc\beta$. Hence in this case, putting
$e_{\gamma\alpha}
=e_{\gamma\beta}$, we have
$$\partial_k\log |a_i(h_k+e_{\gamma\alpha})|=0.\tag 1.13 $$
Also, in case $\gamma=\beta$ and $j\ne k$ (hence
$|\beta|\ge 2$), we have
$$-h_k+\frac{h'_k}{h'_{ki}}(h_{ki}-h_{ik})=
-h_j+\frac{h'j}{h'_{ji}}(h_{ji}-h_{ij}).$$
Since the left- and the right-hand side of
the formula above are functions of $x_k$
and $x_j$ respectively, it follows that they are
constants. Let $e_{\beta\alpha}$
be this common constant. Then, in case $|\beta|\ge
2$, we have
$$h_j+e_{\beta\alpha}=c_{ji}(h_{ji}-h_{ij})\tag 1.14$$
for any $j\in\beta$, where $c_{ji}$ are non-zero constants.
Moreover, we put
$$h_j=|h_{ji}-h_{ij}|,\qquad e_{\beta\alpha}=0 \tag 1.15$$
in case $|\beta|=1$. Then, by (1.13), (1.14), and
(1.15), we see that
$$\partial_k\log \left(a_i\prod_{\gamma\pcc\alpha}\prod_{l\in\gamma}
|h_l+e_{\gamma\alpha}|\right)=0$$
for every $k\in\gamma$ and $\gamma\pcc\alpha$.

Now, applying the same argument as for minimal
element, we obtain functions
$h_l$ $(l\in\alpha)$ so that the function
$$\log \left(a_i\prod_{\gamma\pcc\alpha}\prod_{l\in\gamma}|h_l+
e_{\gamma\alpha}|\right)
+\log \prod_{l\in\alpha\atop l\ne i}|h_l-h_i|$$
depends only on $x_i$ for every $i\in\alpha$ (if
$|\alpha|=1$, the second term
does not appear, and $h_i$ $(i\in\alpha)$ remains undetermined).
Hence $a_i$ $(i\in\alpha)$ can be taken so that
this function is equal to zero. This completes
the induction.
The remaining part of the proposition is easy.
\qed
\enddemo

We now assume that each equivalence class $\alpha\in\A$
consists of
successive numbers:
$$\alpha=\{s(\alpha),s(\alpha)+1,\dots,t(\alpha)\}.$$
Let $m(i)$ $(i\in\alpha)$ be the number of functions
$h_j$ $(j\in\alpha)$
such that $h_i>h_j$ at $p$. Also, let $m(\alpha,\beta)$
be the number of negative
functions in
$$\{h_i+e_{\alpha\beta}\ |\ i\in\alpha\}.$$
Let $\n(\alpha)$ be the set of $\beta\in\A$ such
that $\alpha\pcc \beta$
and there is no $\gamma$ between $\alpha$ and
$\beta$. Also, put
$$u_{\alpha}=\prod_{\gamma\pcc\alpha}\prod_{l\in\gamma}(h_l+
e_{\gamma\alpha}).$$

\proclaim{Proposition 1.13} For a suitably chosen basis $F_1,\dots,F_n$
of
$\Cal F$, the functions $b_{ij}=a_{ij}/a_i$ $(i\in\alpha)$ are given
by:
$$b_{ij}=\cases (-1)^{m(i)}(-h_i)^{j-s(\alpha)} & (j\in\alpha)\\
(-1)^{m(i)-1+m(\alpha,\beta)}(h_i+e_{\alpha\beta})^{-1} &
(j=t(\beta),\ \beta\in \n(\alpha))\\
0 & (\text{otherwise}).
\endcases$$
Moreover,
$$\sum_{\beta\succeq\alpha}\sum_{j\in\beta}a_jb_{jk}=
\cases |u_{\alpha}|^{-1} & (k=t(\alpha))\\
0 & (k\ne t(\alpha)).
\endcases$$
\endproclaim
\demo{Proof} Put
$$\wt b_{ij}=\cases (-1)^{m(i)}(-h_i)^{j-s(\alpha)} & (j\in\alpha)\\
(-1)^{m(i)-1+m(\alpha,\beta)}(h_i+e_{\alpha\beta})^{-1} &
(j=t(\beta),\ \beta\in \n(\alpha))\\
0 & (\text{otherwise}).
\endcases$$
Then a direct computation shows that
$$\sum_{\beta\succeq\alpha}\sum_{j\in\beta}a_j\wt b_{jk}=
\cases |u_{\alpha}|^{-1} & (k=t(\alpha))\\
0 & (k\ne t(\alpha)).
\endcases$$
Let $B$ be the $n\times n$ matrix $(\wt
b_{ij})$, and let $B_i$ be the $(n-1)\times
(n-1)$ matrix obtained by deleting $i$-th row and
$t(\alpha)$-th column from $B$,
where $i\in\alpha$.

\proclaim{Lemma 1.14}
\roster
\runinitem"(1)" $\det B\ne 0$.
\item"(2)" $\det B_i\ne 0$.
\endroster
\endproclaim
\demo{Proof} It is easily seen that $\det B$
is equal to $\prod_{\alpha\in\A}
\det B^{\alpha}$, where $B^{\alpha}=(\wt b_{ij})_{i,j\in\alpha}$. Since
$\det B^{\alpha}\ne 0$ by Vandermonde's formula, (1) follows.
(2) is similar.
\qed
\enddemo

Define functions $c_{kj}$ $(1\le k,j\le n)$ by the
formula
$$b_{ij}=\sum_k \wt b_{ik}c_{kj}.\tag 1.16$$
To prove Proposition 1.13 it suffices to show
that $c_{kj}$ are constants.
Let $\alpha\in\A$ and $l\in\alpha$.
We claim that $\partial_lc_{t(\alpha),j}=0$
for any $j$. In fact, since $\sum_ia_ib_{ij}$ are
constants, we obtain
$$0=\partial_l\sum_ia_ib_{ij}=
\partial_l\sum_{\beta\succeq\alpha}\sum_{i\in\beta}a_ib_{ij}.$$
By (1.16) we also have
$$\sum_{\beta\succeq\alpha}\sum_{i\in\beta}a_ib_{ij}=
|u_{\alpha}|^{-1}c_{t(\alpha),j}.$$
Thus, it follows that $\partial_lc_{t(\alpha),j}=0$.

Moreover, the formula (1.16) implies that
$$0=\sum_{k\ne t(\alpha)}\wt b_{ik}\partial_lc_{kj}=0,$$
provided $i\ne l$ and $\l\in\alpha$. It then follows
from Lemma 1.14 that
$$\partial_lc_{kj}=0\qquad (l\in\alpha,\ k\ne t(\alpha)).$$
This completes the proof.
\qed
\enddemo

Correspondingly, $(f_{ij})=(a_{ij})^{-1}$ is given as follows.

\proclaim{Proposition 1.15} Suppose $i\in\alpha$. Then:
$$f_{ij}=\cases
|u_{\alpha}|\Cal S_{t(\alpha)-i}(h_l\,;\, l\in\alpha-\{j\}) & (j\in\alpha)\\
|u_{\alpha}| \sum_{m=0}^{t(\alpha)-i} e_{\alpha\beta}^m
\Cal S_{t(\alpha)-i-m} (h_l\,;\, l\in\alpha) \quad & (j\in\gamma\succeq\beta,
\beta\in \n(\alpha))\\
0 & (\text{otherwise}),
\endcases$$
where $\Cal S_m (h_l\,;\, l\in\alpha)$ stands for the
elementary symmetric
function of degree $m$ $(0\le m\le |\alpha|)$ with
respect to $|\alpha|$
functions $\{h_l\ |\ l\in\alpha\}$.
\endproclaim

The proof is straightforward. Put
$$v_i=u_{\alpha}\Cal S_{t(\alpha)-i+1}
(h_l\,;\,l\in\alpha)\qquad (i\in\alpha,\
\alpha\in\A),$$
and let $\Cal V$ be the vector space
of functions on $U$ spanned by constant
functions and $v_i$ $(1\le i\le n)$. Put
$$\frak k= \{\sgrad(v)\ |\ v\in\Cal V\}.$$

\proclaim{Proposition 1.16}
\roster
\item"(1)" $[Y,\,W_j]=[Y,\,IW_j] =0\quad (Y\in\frak k,\ \text{any }j)$.
\item"(2)" $\{Y,\,F\}=0\qquad (Y\in\frak k,\ F\in\Cal F)$.
\item"(3)" $\frak k$ is the commutative Lie algebra
of infinitesimal
automorphisms of the K\"ahler manifold $M$ on $U$.
\endroster
\endproclaim
\demo{Proof} Put $Y=\sgrad v_m$. Then
$$Y=\sum_ja_j(W_jv_m)IW_j.$$
Hence, by Propositions 1.4 and 1.5 we have
$[Y,IW_j]=0$ for any $j$.
Since
$$[W_i,IW_i]=\left(-\partial_i\log a_i-\frac{h''_i}{h'_i}\right)IW_i+
\sum_{j\ne i}\frac{a_j}{a_i}(\partial_j\log a_i)IW_j,$$
it follows that
$$\align
& [W_i,Y]= a_i\left(\partial^2_iv_m-
(\partial_iv_m)\frac{h''_i}{h'_i}\right)IW_i\\
& +\sum_{j\ne i}a_j(\partial_i\partial_jv_m+
(\partial_i\log a_j)(\partial_jv_m)
+(\partial_j\log a_i)(\partial_iv_m))IW_j.
\endalign$$
Then, it is easily seen that each term
in the right-hand side of the formula above
vanishes. Thus (1) follows.

{}From (1) it follows that
$$0=\{Y,W_i^2+(IW_i)^2\}=\sum_jb_{ij}\{Y,F_j\}.$$
This indicates (2). In particular, we have $\{Y,E\}=0$,
which implies that
$Y$ is an infinitesimal isometry. The property (1)
also implies that
$(\Cal L_YI)W_i= (\Cal L_YI)IW_i=0$ for any $i$, where
$\Cal L_Y$ denotes the Lie
derivative with respect to $Y$. Hence it follows
that $\Cal L_YI=0$. Moreover,
putting $Y'=\sgrad v_l$ $(l\ne m)$, we have
$$i_{[Y,Y']}\omega=-\Cal L_Y(dv_l)=-d(Yv_l)=0.$$
Hence $[Y,Y']=0$, and (3) follows.
\qed
\enddemo

\specialhead 2. Summing up the local data
\endspecialhead
In the previous section we have given, for
each $p\in M^1$, the neighborhood $U$,
the constants $e_{\alpha\beta}$, the functions $h_i$, $a_i$, $b_{ij}$,
$v_i$, and
the basis $F_1,\dots,F_n$ of $\Cal F$, as well
as their numbering.
{}From now on, to clarify the point-dependence, we
shall write $U^{(p)}$, $h_i^{(p)}$,
$F_i^{(p)}$, etc. instead.
We assume that each neighborhood $U^{(p)}$ is taken
to be a small distance ball
centered at $p$ so that it is convex.

Take a point $p_0\in M^1$ and fix it.
Let $M^{1,0}$ be the connected component
of $M^1$ containing $p_0$. Let $p\in M^{1,0}$, and
let $\gamma(t)$ $(0\le t\le 1)$
be a curve in $M^{1,0}$ such that $\gamma(0)=p_0$,
$\gamma(1)=p$. Along the
curve $\gamma$ there is a unique numbering of
$\{D_i\}$
so that $t\mapsto (D_i)_{\gamma(t)}$ is continuous.

Since the relations $i\pc j$ are locally stable,
we have the
following

\proclaim{Lemma 2.1} The relation $\pc$ on $\{1,\dots,n\}$ is
constant along
the curve $\gamma$. In particular, the partially ordered
set $\A$ is constantly
defined along $\gamma$.
\endproclaim

We put
$$F_{\alpha}^{(q)}(\lambda)=\sum_{i\in\alpha}(-\lambda)^{i-s(\alpha)}
F_i^{(q)}.$$
As is easily seen,
$$\align
F_{\alpha}^{(q)}(\lambda) & =|u_{\alpha}^{(q)}|\sum_{j\in\alpha}
\prod_{k\in\alpha\atop k\ne j}(h_k^{(q)}-\lambda)\cdot (V_j^2+(IV_j)^2)\\
& +|u_{\alpha}^{(q)}|\sum_{\beta\in \n(\alpha)}
\frac{\prod_{l\in\alpha}(h_l^{(q)}+
e_{\alpha\beta}^{(q)})-\prod_{l\in\alpha}
(h_l^{(q)}-\lambda)}{e_{\alpha\beta}^{(q)}
+\lambda}\sum_{\gamma\succeq\beta}\sum_{j\in\gamma} (V_j^2+(IV_j)^2).
\endalign$$

\proclaim{Proposition 2.2} There are constants $c_{\alpha}\ne 0$ and
$d_{\alpha}$ such that
$$F_{\alpha}^{(p)}(c_{\alpha}\lambda-d_{\alpha})
=\left(c_{\alpha}^{|\alpha|-1}
\prod_{\beta\pc\alpha}c_{\beta}^{|\beta|}\right)
F_{\alpha}^{(p_0)} (\lambda).$$
If $\{h_i^{(p)}\}$ and $e_{\alpha\beta}^{(p)}$ are suitably chosen, then
those
constants become
$$c_{\alpha}=1,\quad d_{\alpha}=0\qquad (\text{any }\alpha).$$
\endproclaim
\demo{Proof} Let $I_{\gamma(t)}$ be the connected component of
the intersection
of $U^{(\gamma(t))}$ and the image of $\gamma$ containing
the point $\gamma(t)$.
Suppose that $I_{\gamma(t_1)}\cap I_{\gamma(t_2)}\ne\emptyset$. Then, by
Proposition 1.11 there are constants $\wt c_{\alpha}(\ne 0)$
and
$\wt d_{\alpha}$ such that
$$h_i^{(\gamma(t_1))}=\wt c_{\alpha}h_i^{(\gamma(t_2))}-\wt d_{\alpha},\quad
e_{\alpha\beta}^{(\gamma(t_1))}=\wt c_{\alpha}e_{\alpha\beta}^{(\gamma(t_2))}
+\wt d_{\alpha}$$
on $U^{(\gamma(t_1))}\cap U^{(\gamma(t_2))}$. This implies that
$$F_{\alpha}^{(\gamma(t_1))}(\wt c_{\alpha}\lambda-\wt d_{\alpha})
=\left({\wt c}_{\alpha}^{|\alpha|-1}
\prod_{\beta\pc\alpha}{\wt c}_{\beta}^{|\beta|}
\right) F_{\alpha}^{(\gamma(t_2))} (\lambda).$$
Taking a finite number of points on $\gamma$
and iterating this argument
successively, we obtain the former half of the
proposition.

Now, putting
$$\wt h_i^{(p)}=c_{\alpha}^{-1}(h_i^{(p)}+d_{\alpha}),\quad
\wt e_{\alpha\beta}^{(p)}=
c_{\alpha}^{-1}(e_{\alpha\beta}^{(p)}-d_{\alpha}),$$
and denoting by $\wt F_{\alpha}^{(p)}(\lambda)$ the corresponding polynomial,
we have
$$\wt F_{\alpha}^{(p)}(\lambda)=\left(c_{\alpha}^{|\alpha|-1}
\prod_{\beta\pc\alpha}c_{\beta}^{|\beta|}\right)^{-1}F_{\alpha}^{(p)}
(c_{\alpha}\lambda-d_{\alpha})=F_{\alpha}^{(p_0)} (\lambda).$$
\qed
\enddemo

Now, let us suppose that $p=p_0$, i.e., $\gamma$
is a loop.
Let $\nu$ be the permutation of the indices
$1,\dots,n$ defined by
$$(D_i)_{\gamma(1)}=(D_{\nu(i)})_{\gamma(0)}.$$

\proclaim{Proposition 2.3} $\nu$ is the identity.
\endproclaim
\demo{Proof} Since $\nu$ preserves the relation $\pc$, it
induces an
automorphism of the partially ordered set. We also
denote it by $\nu$.
Then, by Proposition 2.2 there are constants $c_{\alpha}\ne
0$ and $d_{\alpha}$
such that
$$F_{\nu(\alpha)}^{(p_0)}(c_{\alpha}\lambda-d_{\alpha})
=\left(c_{\alpha}^{|\alpha|-1}\prod_{\beta\pc\alpha}
c_{\beta}^{|\beta|}\right)
F_{\alpha}^{(p_0)}(\lambda).$$
This formula clearly indicates that $\nu(\alpha)=\alpha$ and $F_i^{(p)}$
$(i\in\alpha)$ is written as a linear combination of
$F_{\nu(j)}^{(p)}$ $(i\le j\le
t(\alpha))$. Thus $\nu(i)=i$ for every $i\in\alpha$.
\qed
\enddemo

This proposition implies that the subbundles $D_i$ are
globally defined
on $M^{1,0}$.

\proclaim{Proposition 2.4} Suppose that $\{h_i^{(q)}\}$
and $e_{\alpha\beta}^{(q)}$ $(q\in M^{1,0})$ are taken so that
$F_{\alpha}^{(q)}(\lambda)= F_{\alpha}^{(p_0)}(\lambda)$ for all
$\alpha\in \A$. Then for any $p,q\in M^{1,0}$ such
that $U^{(p)}\cap U^{(q)}
\ne \emptyset$,
$e_{\alpha\beta}^{(p)}=e_{\alpha\beta}^{(q)}$ and $h_i^{(p)}=h_i^{(q)}$
for any $i$ on the intersection. Hence, there
are functions $\{h_i\}$ on $M^{1,0}$
such that $h_i=h_i^{(p)}$ on $U^{(p)}$.
\endproclaim
\demo{Proof} Since $F_{\alpha}^{(p)}(\lambda)= F_{\alpha}^{(q)}(\lambda)$,
the proposition follows from the proof of Proposition
2.2.
\qed
\enddemo

\specialhead 3. Structure of $M-M^1$
\endspecialhead
In the previous section we have obtained the
constants $e_{\alpha\beta}$, the
functions $h_i$, $a_i$, $b_{ij}$, $f_{ij}$ on $M^{1,0}$
and the basis $\{F_i\}$ of $\Cal F$. Also,
the functions
$v\in\Cal V$ and the vector fields $Y\in \frak
k$ are now defined on $M^{1,0}$, and the
properties described in Proposition 1.16 holds on $M^{1,0}$.
Since for each
$\alpha$ the functions $h_i$ $(i\in\alpha)$ take mutually distinct
values at every
point in $M^{1,0}$, we may assume that the
numbering of $\{D_i\}$ is chosen so that
$$h_{s(\alpha)}> h_{s(\alpha)+1}>\dots >h_{t(\alpha)}\qquad (\alpha\in\A)$$
on $M^{1,0}$. We put
$$\align
v_{\alpha}(\lambda)&=u_{\alpha}\prod_{i\in\alpha}(h_i-\lambda)\\
&=\sum_{m=0}^{|\alpha|-1}(-\lambda)^m v_{s(\alpha)+m}
+(-\lambda)^{|\alpha|}u_{\alpha}.
\endalign$$

The main purpose of this section is to
prove the following

\proclaim{Theorem 3.1} $M-M^1$ is equal with a locally
finite union of closed,
totally geodesic, complex hypersurfaces $L$. In particular $M^1$
is connected and
dense in $M$. Moreover, the functions $h_i$ are
continuously extended to the whole
$M$ and the subbundles $D_i$ are smoothly extended
to $M^0$, and they possess the
following properties:
\roster
\item"(1)" $h_i+e_{\alpha\beta}$
$(i\in\alpha,\ \alpha\pcc \beta)$ are everywhere
nonzero on $M$;
\item"(2)" $h_{s(\alpha)}> h_{s(\alpha)+1}>
\dots >h_{t(\alpha)}\quad\text{on }M^0
\quad (\alpha\in\A)$;
\item"(3)" $h_i$ are $C^{\infty}$ functions on $M^0$;
\item"(4)" $\Cal S_m(h_i\,;\,i\in\alpha)$
$(1\le m\le |\alpha|)$ are $C^{\infty}$
functions on $M$;
\item"(5)" $D_{\alpha}=\sum_{i\in\alpha}D_i$ $(\alpha\in\A)$ are extended to
the whole $M$ as $C^{\infty}$ subbundles of $TM$;
\item"(6)" The vector fields $Y\in\frak k$ are globally
defined and of $C^{\infty}$
on $M$;
\item"(7)" For each hypersurface $L$ there are $\alpha\in\A$
and $c\in\R$
such that a connected component of the set
of zeros of $\sgrad v_{\alpha}(c)\in\frak k$
coincides with $L$, and $v_{\alpha}(c)$ vanishes on $L$.
\endroster
\endproclaim

We shall call $\{h_i\}$ and $\{e_{\alpha\beta}\}$ {\it the
fundamental functions} and
{\it the conjunction constants} of the K\"ahler-Liouville manifold
$(M,\Cal F)$
(of type (A)) respectively. Note that if $\{h'_i\}$
and
$\{e'_{\alpha\beta}\}$ are other choice, then there are constants
$c_{\alpha}\ne 0$
and $d_{\alpha}$ such that $e'_{\alpha\beta}
=c_{\alpha}e_{\alpha\beta}+d_{\alpha}$
and for $i\in\alpha$,
$$h'_i=\cases
c_{\alpha}h_i-d_{\alpha}\quad & (c_{\alpha}>0)\\
c_{\alpha}h_{t(\alpha)+s(\alpha)-i}-d_{\alpha}\quad & (c_{\alpha}<0)
\endcases$$

The following theorem is an immediate consequence of
the theorem above and
Proposition 1.16.

\proclaim{Theorem 3.2} The geodesic flow of a K\"ahler-Liouville
manifold of
type (A) is integrable with respect to the
first integrals in $\Cal F$ and $\frak k$.
\endproclaim

We begin the proof of Theorem 3.1 with
a characterization of the functions $v_i$ in
terms of $F_j$.

\proclaim{Lemma 3.3} Let $i\in\alpha$. If $i\ne s(\alpha)$, the
function $v_i$ is
equal with a linear combination of $\tr F_{i-1},\dots,\tr
F_{t(\alpha)}$ and $\tr F_j$
$(j\in \beta,\ \beta\pcc\alpha)$. If $i=s(\alpha)$ and $\alpha$ is
not a maximal
element, then $v_i$ is equal with a linear
combination of $\tr F_j$ $(j\in
\beta,\ \beta\pc\alpha)$ and $\tr F_{t(\gamma)}$, where $\gamma$ is
any element of $\n(\alpha)$. Finally, if $i=s(\alpha)$ and
$\alpha$ is a maximal
element, then $v_iu_{\alpha}^{|\alpha|-1}$ is a polynomial of $\det
F_{t(\alpha)-1}$
and $\tr F_j$ $(j\in\beta,\ \beta\pc\alpha)$.
\endproclaim

The proof is straightforward. The lemma above implies
that there is a $C^{\infty}$
function on $M$ expressed by the traces of
$F_j$ that coincides
with $v_i$ on $M^{1,0}$, unless $i=s(\alpha)$ and $\alpha$
is maximal. Though
the expression is not unique in general, we
take one and fix it. We denote
the extended function by the same symbol $v_i$.
Note that every $u_{\beta}$
$(\beta\in\A)$ is a linear combination of those functions.
Similarly, if
$i=s(\alpha)$ and $\alpha$ is maximal, then there is
a $C^{\infty}$ function on $M$
that coincides with $v_iu_{\alpha}^{|\alpha|-1}$ on $M^{1,0}$.

Let $M^{0,0}$ be the connected component of $M^0$
that contains $M^{1,0}$. For
technical reason we shall first show that $M^{0,0}\cap
M^1$ is connected and
dense in $M^{0,0}$, and that the functions $h_i$
are smoothly extended to $M^{0,0}$.
Let $c(t)$ be a geodesic such that
$$p=c(0)\in M^{1,0}, \quad c([0,t_0))\subset M^{1,0},\quad q=c(t_0)\in
M^{0,0}-M^{1,0}.$$
Since $q\in M^0$, every simultaneous eigenspace of $\{F_q^e\
|\ F\in\Cal F\}$ is
(complex) one-dimensional. Hence there is a neighborhood $U$
of $q$, and there
are subbundles $D_i$ of $TM$ on $U$ such
that each $D_i$ coincides with that
on $M^{1,0}$ on the connected component of $U\cap
M^{1,0}$ containing a curve segment
of the form $c((t_0-\epsilon,t_0))$, $\epsilon >0$.

Since the polynomial $v_{\alpha}(\lambda)$ has $|\alpha|$ real roots
at each point in
$M^{1,0}$, so does at $q$ if $u_{\alpha}(q)\ne 0$.
In this case,
denoting those roots by $h_i(q)$ $(i\in\alpha)$, $h_{s(\alpha)}
(q)\ge\dots\ge h_{t(\alpha)}(q)$, we have the continuous extension of
$h_i$ up to $q$.

\proclaim{Lemma 3.4}
\roster
\runinitem"(1)" $u_{\alpha}(q)\ne 0$ for any $\alpha$.
\item"(2)" $h_{s(\alpha)}(q)>\dots> h_{t(\alpha)}(q)$.
\endroster
\endproclaim
\demo{Proof} If $u_{\alpha}(q)=0$ for some $\alpha$, then we
have $F_{t(\alpha)}=0$
at $q$, contradicting $q\in M^0$. Hence $u_{\alpha}(q)\ne 0$
for every $\alpha$, and
the functions $h_i$ are well-defined at $q$. Since
the eigenvalues of the endomorphism
$F_i^e$ at each point in $M^{1,0}$ are given
by $f_{ij}$ described in Proposition 1.15, so
are at $q$ by continuity. Therefore, if $h_i(q)=h_{i+1}(q)$
for some $i,i+1\in\alpha$,
then one can easily see that $D_i$ and
$D_{i+1}$ cannot be separated by means of the
eigenvalues, which contradicts the facts that every simultaneous
eigenspace is
one-dimensional.\qed
\enddemo

By virtue of the lemma above, we see
that the polynomial $v_{\alpha}(\lambda)$ has
$|\alpha|$ distinct real roots on $U$, if $U$
is taken small enough. Denoting those roots
again by $h_i$ $(i\in\alpha)$, $h_{s(\alpha)}>
\dots> h_{t(\alpha)}$, we obtain
the smooth
extension of the functions $h_i$ to $M^{1,0}\cup U$.

\proclaim{Lemma 3.5} There is some $h_i$ such that
$dh_i=0$ at $q$. In this case we
also have:
\roster
\item"(1)" The hessian $\hess h_i$ of $h_i$ at
$q$ is given by
$$\hess h_i(X,Y)=ag([X]_{D_i},[Y]_{D_i}),\qquad X,Y\in T_qM,$$
where $a\in\R$, $a\ne 0$, and $[X]_{D_i}$ denotes the
$D_i$-component of $X$;
\item"(2)" $\dot c(t_0)$ is not orthogonal to $D_i$.
\endroster
\endproclaim
\demo{Proof} Put $U'=M^{1,0}\cap U$, and let $\{\wt a_i\}$
be functions around $q$
given in Proposition 1.1. Then we have
$$d\log \wt a_i\equiv -d\log \left(|u_{\alpha}|\prod_{j\in\alpha\atop j\ne i}
|h_j-h_i|\right)\quad \mod D_i$$
on $U'$. Clearly, it is also valid on
the closure of $U'$ in $U$ by continuity. Hence,
if every derivatives $h'_i$ does not vanish at
$q$, then it follows that
$q\in M^1$, a contradiction. Thus there is some
$i$ such that $h'_i= 0$ at $q$.

Put $b=h_i(q)$ and suppose $i\in\alpha$. Then $dv_{\alpha}(b)= 0$
and
$$\hess v_{\alpha}(b)=\wt a\,\hess h_i$$
at $q$, where $\wt a$ is a non-zero
constant. Put $Y=\sgrad v_{\alpha}(b)$. Then we
have
$$\hess v_{\alpha}(b)(X,Z)=g(\nabla_XY,IZ).$$
Moreover, since $Y$ is an infinitesimal automorphism of
the K\"ahler manifold $M$ on
$U'$, we have $\nabla_{IX}Y=I\nabla_XY$. Hence $\hess v_{\alpha}(b)$ is
a hermitian
form on $U'$, and so is at $q$
by continuity. Since $(\hess h_i)(X,Z)=0$ at $q$ if $X$
or $Z$ is orthogonal to $D_i$, it follows
that
$$(\hess h_i)(X,Z)=(\hess h_i)([X]_{D_i},[Z]_{D_i})=ag([X]_{D_i},[Z]_{D_i})$$
at $q$ for some constant $a$.

Now, we show $(\hess v_{\alpha}(b))(\dot c(t_0),\cdot)\ne 0$, which
will prove that
$a\ne 0$, and (2). Since $Y$ is a
Jacobi field along $c(t)$ $(0\le t< t_0)$, it
satisfies the equation of Jacobi field up to
$q=c(t_0)$ by continuity. Hence,
$Y_q$ being $0$, we have $\nabla_{\dot c(t_0)}Y\ne 0$
at $q$. From this it follows
that
$$(\hess v_{\alpha}(b))(\dot c(t_0),\cdot)\ne 0.$$
\qed
\enddemo

\proclaim{Lemma 3.6} There is a constant $\epsilon>0$ such
that
$c((t_0,t_0+\epsilon))\subset M^{1,0}$.
\endproclaim
\demo{Proof} Let $S$ be the subspace of $T_qM$
spanned by $\dot c(t_0)$ and
$I\dot c(t_0)$, and let $B$ be the image
of the $\epsilon$-ball $\{V\in S\ |\ |V|<
\epsilon\}$ in S by the exponential mapping $\Exp_q:T_qM\to
M$. Then,
it follows from the previous lemma that $h'_i\ne
0$ at $q'$ for every $i$ and
$q'\in B-\{q\}$, provided $\epsilon$ small enough. It also
follows from the proof of
the previous lemma that if $q'\in U$ lies
on the boundary of $U'$, then some $h'_i$
vanishes at $q'$. Hence we have $B-\{q\}\subset M^{1,0}$.
\qed
\enddemo

\proclaim{Proposition 3.7} $M^{0,0}\cap M^1=M^{1,0}$, and it is dense
in $M^{0,0}$.
Also, the functions $h_i$ and the subbundles $D_i$
extend smoothly to $M^{0,0}$
and satisfies
$$h_{s(\alpha)}>\dots> h_{t(\alpha)}
\quad (\alpha\in\A),\quad dh_i|_{D_j}=0\quad
\text{if }i\ne j.$$
Moreover, $D_i$ are the simultaneous eigenspaces of the
endomorphisms
$\{F^e\ |\ F\in\Cal F\}$ at each point in
$M^{0,0}$.
\endproclaim
\demo{Proof} Let $N$ be the closure of $M^{1,0}$
in $M^{0,0}$, and assume that
$N\ne M^{0,0}$. Let $q\in N$ be a boundary
point of $N$, and let $U$ be an
open distance ball centered at $q$, which is
small enough so that it is convex and
contained in $M^{0,0}$. Let $p_0$ and $q_0$ be
two points in $U$ such that
$p_0\in M^{1,0}$ and $q_0\not\in N$. Let $c(t)$ $(0\le
t\le T)$ be the minimal
geodesic from $p_0$ to $q_0$, and let $t=t_1$
be the earliest time when $c(t)$
meets the boundary of $N$.

Applying the argument above to the geodesic $c(t)$
and the point $p_1=c(t_1)$,
we see that there is $\epsilon>0$ such that
$c((t_1,t_1+\epsilon))\subset M^{1,0}$.
Iterating this procedure successively, we obtain a sequence
of times
$0=t_0<t_1<t_2<\dots<T$ such that
$$c(t_k)\not\in M^{1,0},\quad c((t_{k-1},t_k))
\subset M^{1,0}\qquad (k\ge 1).$$
We claim that the number of such $t_k$
is finite. In fact, suppose that it is not the
case, and put $t_{\infty}=\lim_{k\to\infty}t_k\le T$. Then, as seen
above, any
vector fields of the form
$$Y=\sgrad v_{\alpha}(b)\qquad (\alpha\in\A,\ b\in\R)$$
are Jacobi fields along $c(t)$ $(0\le t\le t_{\infty})$,
and among them there is
$Y_k$ that vanish at $c(t_k)$ for every $k\ge
1$. Let $\wt{\frak k}$ be the
vector space of Jacobi fields spanned by such
$Y$. Since $\wt{\frak k}$ coincides with
$\frak k$ on $c([0,T])\cap M^{1,0}$, it follows that
it is $n$-dimensional, and
satisfies
$$g(Y,\nabla_{\dot c}Y')=g(\nabla_{\dot c}Y,Y')\qquad Y,Y'\in\wt{\frak k}.$$
{}From this property one can easily conclude that
the set of points $t$ such that
some non-zero $Y\in\wt{\frak k}$ vanishes at $c(t)$ is
discrete. On the other hand,
choosing a subsequence if necessary, we obtain $Y_{\infty}\in\wt{\frak
k}-\{0\}$
as a limit of (constant multiples of) $\{Y_k\}$.
Since $Y_{\infty}$ vanishes at
$c(t_{\infty})$, it is a contradiction. Hence there are
only finitely many
$t_k$; $t_1,\dots,t_l$.

Now, again by Lemma 3.6 we see that
$c((t_l,T))\subset M^{1,0}$. However, this
contradicts the fact that $q_0=c(T)\not\in N$. Thus we
conclude that $N=M^{0,0}$, and
$M^{0,0}\cap M^1=M^{1,0}$. The remaining part is obvious.
\qed
\enddemo

Next, we shall prove that $M^0$ is connected
and dense in $M$. Let $q$ be a boundary
point of $M^{0,0}$, and assume that $u_{\alpha}(q)\ne 0$
for every $\alpha$. Then in
the same way as above the functions $h_i$
are continuously extended to $M^{0,0}
\cup\{q\}$. Also, $T_qM$ is decomposed to the simultaneous
eigenspaces of the endomorphisms $\{F^e\ |\ F\in\Cal F\}$.
Clearly those subspaces
are uniquely extended as $C^{\infty}$ subbundles around $q$
so that they are
sums of simultaneous eigenspaces at each point.

\proclaim{Lemma 3.8} Let $q\in M^s\cap \ol {M^{0,0}}$, and
assume that
$u_{\alpha}(q)\ne 0$ for any $\alpha$. Let $D$ be
one of the simultaneous
eigenspaces of $\{F^e\ |\ F\in\Cal F\}$ at $q$,
and let $\dim D=m$. Then there
is $\alpha\in\A$ and its subset ${\alpha}'$ consisting of
successive numbers
$i,\dots, i+m-1$ such that the extended
subbundle (also denoted by $D$) is equal with
$\sum_{j\in{\alpha}'}D_j$ on $M^{0,0}$
near $q$. Moreover, there is a constant $c$
such that:
\roster
\item"(1)" $h_i(q)=\dots =h_{i+m-1}(q)=c,\qquad h_j(q)\ne c\quad
(j\in\alpha-{\alpha}')$;
\item"(2)" The functions $\Cal S_k(h_j\,;\,j\in{\alpha}')$ on $M^{0,0}$ can
be smoothly
extended around $q$\quad $(1\le k\le m)$;
\item"(3)" $dv_{\alpha}(c)=0$ at $q$ if $m\ge 2$;
\item"(4)" The hessian of $v_{\alpha}(c)$ vanishes at $q$
if $m\ge 3$.
\endroster
Also, there exists at least one simultaneous eigenspace
$D$ with $\dim D\ge 2$.
\endproclaim
\demo{Proof} Let $U$ be a neighborhood of $q$
where the subbundle $D$ is defined.
Since $D$ is a sum of simultaneous eigenspaces
at each point, it is a sum of $D_j$
on $U\cap M^{0,0}$. Let us consider the endomorphisms
$F_{t(\alpha)}^e$. It is
$|u_{\alpha}|$ times the identity on
$\sum_{\beta\succeq\alpha}D_{\beta}$, and
$0$ on the orthogonal complement.
Therefore the subbundles $D_{\alpha}$ are continuously extended to
$q$ so that
$\sum_{\beta\succeq\alpha}D_{\beta}$
is still the eigenspace of $F_{t(\alpha)}^e$
corresponding to the eigenvalue $|u_{\alpha}(q)|$. Since $D_{\alpha}$ is
a sum of
simultaneous eigenspaces at $q$, it is consequently extended
on $U$ as the subbundle
of $TM$. Clearly $D\subset D_{\alpha}$ for some $\alpha$.
Also, the endomorphism
$F_j^e$ at $q$ is a constant multiple of
the identity on $D_{\alpha}$ if $j\notin
\alpha$, and the eigenvalues of $F_j^e$ $(j \in
\alpha)$ on $D_{\alpha}$ are
$$|u_{\alpha}(q)|\Cal S_{t(\alpha)-j}
(h_l(q)\,;\,l\in\alpha-\{k\})\qquad (k\in\alpha)$$
at $q$. This implies that there is a
subset ${\alpha}'$ of $\alpha$
consisting of successive numbers $i,\dots,i+m-1$ such that
$$\align
& D=\sum_{j\in{\alpha}'} D_j\quad\text{on }U\cap M^{0,0}\\
& h_i(q)=\dots =h_{i+m-1}(q)=c\\
& h_j(q)\ne c\qquad \text{for }j\in\alpha-{\alpha}'.
\endalign$$

To prove (2) we note that $\sum_{j\in\alpha} h_j$
is extended as the $C^{\infty}$
function around $q$, because $v_{t(\alpha)}$ and $u_{\alpha}$ are
of $C^{\infty}$, and
$u_{\alpha}(q)\ne 0$. Hence the symmetric functions of the
eigenvalues of the
endomorphisms
$$u_{\alpha}^{-1}\left(\left(\sum_{j\in\alpha}h_j\right)F_{t(\alpha)}^e
-F_{t(\alpha)-1}^e\right)$$
on $D$, which are
$$\Cal S_l\left(\sum_{j\in {\alpha}'-\{k\}}h_j;\ k\in{\alpha}'\right),$$
on $M^{0,0}$, are $C^{\infty}$ functions around $q$. This
proves (2).

To prove the assertions (3) and (4), assume
$m\ge 2$, and let $G(\lambda)$ be the
endomorphism of $D_{\alpha}$ defined to be the restriction
of
$$F_{t(\alpha)-1}^e-\left(\sum_{j\in\alpha}h_j(q)
-\lambda\right)F_{t(\alpha)}^e$$
to $D_{\alpha}$. Then $\det G(\lambda)$ is a $C^{\infty}$
function on $U$. Since
$G(c)=0$ on $D$ at $q$, the order of
the zero $q$ of the function $\det G(c)$ is not
less than $m$. On the other hand, we
have
$$\align
\det G(\lambda) & =|u_{\alpha}|^{|\alpha|} \prod_{j\in\alpha}\left(\sum_
{k\in\alpha}h_k -\sum_{k\in\alpha}h_k(q)-(h_j-\lambda)\right)\\
& = {\epsilon}^{|\alpha|}\sum_{l=0}^{|\alpha|} (-1)^l u_{\alpha}^l\Cal S_l
(h_j-\lambda\,;\,j\in\alpha)
\left(v_{t(\alpha)}-u_{\alpha}\sum_{k\in\alpha}h_k(q)\right)^{|\alpha|-l}
\endalign$$
on $U\cap \ol {M^{0,0}}$, where $\epsilon$ is the
sign of $u_{\alpha}$.
Note that $u_{\alpha}\Cal S_l(h_j-\lambda\,;\,j\in\alpha)$ is
expressed as a linear combination of $v_j$ $(j\in\alpha)$
and $u_{\alpha}$. Hence
the last formula described above also expresses a
$C^{\infty}$ function on $U$
that coincides with $\det G(\lambda)$ on $U\cap \ol
{M^{0,0}}$. In particular
those two functions coincides at $q$ up to
infinite order. Since $d(\det G(c))=0$ at
$q$, and since
$$\Cal S_{|\alpha|-1}(h_j-c\,;\,j\in\alpha)=v_{t(\alpha)}-u_{\alpha}
\sum_{k\in\alpha}h_k(q)=0$$
at $q$, it follows that
$$d(u_{\alpha}\Cal S_{|\alpha|}(h_j-c\,;\,j\in\alpha))=0$$
at $q$. Moreover, if $m\ge 3$, considering the
function $(d/d\lambda)
\det G(\lambda)|_{\lambda=c}$, we also have
$$d(u_{\alpha}\Cal S_{|\alpha|-1}(h_j-c\,;\,j\in\alpha))=0$$
at $q$. Then observing the function $\det G(c)$
again, we see that the hessian of the
function
$$u_{\alpha}\Cal S_{|\alpha|}(h_j-c\,;\,j\in\alpha)=v_{\alpha}(c)$$
vanishes at $q$.

Finally, assume that every simultaneous eigenspace is of
dimension one. Then as is
easily seen, the $|\alpha|$ functions $h_i$ $(i\in\alpha)$ take
mutually different
values at $q$ for any $\alpha$. Hence the
functions $a_{ij}=a_ib_{ij}$ can be
continuously extended to $q$, which implies the linear
independence of $F_1,\dots
F_n$ at $q$.
\qed
\enddemo

\proclaim{Proposition 3.9} $M^0$ is connected and dense in
$M$.
\endproclaim

\demo{Proof} Assume that $M^0$ is not connected, and
let $p,p'\in M^0$ such that
$p\in M^{1,0}$ and $p'$ lies in another component.
Let $c(t)$ $(0\le t\le t_0)$
be a geodesic from $p$ to $p'$. A
slight modification of the geodesic $c$ and the
point $p'$ enables us to assume that $D_i$-component
of $\dot c(0)\in T_pM$ is
not zero for any $i$. Let $t_1$ be
the time such that $q=c(t_1)\in M^s$ and $c(t)
\in M^{0,0}$ for any $t\in [0,t_1)$.

We first claim that every (extended) function $u_{\alpha}$
$(\alpha\in\A)$ does not vanish at $q$. In fact,
assume that $u_{\alpha}(q)=0$.
Since
$$\tr F_{t(\alpha)}=u_{\alpha}\sum_{\beta\succ\alpha} |\beta|$$
and since $F_{t(\alpha)}$ is positive semi-definite on $M^{0,0}$,
it follows that
$F_{t(\alpha)} =0$ at $q$. Let $\zeta_t$ be the
geodesic flow and $\pi:T^*M\to M$
the bundle projection. Then $c(t)=\pi(\zeta_t\lambda_0)$, where $\lambda_0\in
T^*_pM$ is given by
$$\lambda_0(X)=g(\dot c(0),X)\qquad X\in T_pM,$$
and we have
$$F_{t(\alpha)}(\lambda_0) =F_{t(\alpha)}(\zeta_{t_1}\lambda_0)=0.$$
This implies that $\dot c(0)$ is orthogonal to
$D_j$ for any $j\in\beta$, $\beta
\succeq \alpha$, a contradiction. Hence $u_{\alpha}(q)\ne 0$ for
every $\alpha$.

Therefore, as stated above, the functions $h_i$ is
extended up to $q$, and Lemma 3.8 is
applicable. Let $D$ be a simultaneous eigenspace of
the endomorphisms
$\{F^e\ |\ F\in\Cal F\}$ at $q$, and let
$\alpha$ and ${\alpha}'$ be as in Lemma 3.8.
Suppose $m=\dim D\ge 2$, and put $h_j(q)=a$ $(j\in{\alpha}')$.
Then $Y=\sgrad
v_{\alpha}(a)$ is the Jacobi field along the geodesic
$c|_{[0,t_1]}$, and vanishes at
$q=c(t_1)$. Hence $\nabla_{\dot c(t_1)}Y\ne 0$, which implies that
the hessian
of $v_{\alpha}$ at $q$ does not vanish. Thus
we have $m=2$ by Lemma 3.8.

Put ${\Cal F}'=\{F\in\Cal F\ |\ F_q=0\}$, and let
$\dim {\Cal F}'=k$. For each
2-dimensional simultaneous eigenspace $D$, put
$$H_D=F_{\alpha}(a)-\sum_{\beta\in \n(\alpha)}
\epsilon_{\beta}(e_{\alpha\beta}+a)^{-1}
F_{t(\beta)},$$
where $\epsilon_{\beta}$ is the sign of
$\prod_{l\in\alpha}(h_l+e_{\alpha\beta})$.
As is easily seen, the elements $H_D$ form
a basis of ${\Cal F}'$.
For $F\in\Cal F$, let $X_F$ be the vector
field on $T^*M$ defined by
$$i_{X_F}d\theta=-dF,$$
where $\theta$ is the canonical 1-form on the
cotangent bundle. Clearly, the vector
field $X_F$ $(F\in{\Cal F}')$ on $T^*M$ is tangent
to the fibre
$T_q^*M$ at each point $\lambda\in T_q^*M$. Thus we
have the vector field
$X_F|_{T_q^*M}$ on the vector space $T_q^*M$ whose coefficients
are hermitian forms.
Let $S_q^*M$ be the sphere of unit covectors
at $q$, and put
$$\Lambda= \{\lambda\in S_q^*M\ |\ X_F|_{T_q^*M}=0\ \text{at }\lambda \text
{ for some }F\in {\Cal F}'-\{0\}\}.$$

\proclaim{Lemma 3.10}
\roster
\runinitem"(1)" $(X_F)_{\zeta_{t_1}\lambda_0}
\ne 0$ for any $F\in\Cal F-\{0\}$.
\item"(2)" For each 2-dimensional simultaneous eigenspace $D$, there
is
a unit vector $V_D\in D$ such that $\Lambda$
is given by
$$\Lambda=\{\lambda\in S_q^*M\ |\ \lambda(V_D)=\lambda(IV_D)=0
\text{ for some }D\}.$$
In particular the complement of $\Lambda$ in $S_q^*M$
is connected and dense in
$S_q^*M$.
\endroster
\endproclaim
\demo{Proof} Describing
$$F_p=\sum_jb_j(V_j^2+(IV_j)^2),$$
we have
$$\pi_*((X_F)_{\lambda_0})=
2\sum_jb_j(g(V_j,\dot c(0))V_j+g(IV_j,\dot c(0))IV_j).$$
Since the right-hand side does not vanish because
of the assumption on $\dot c(0)$,
it follows that
$$(X_F)_{\zeta_{t_1}\lambda_0}={\zeta_{t_1}}_*((X_F)_{\lambda_0})\ne 0.$$

As is easily seen, the endomorphism $H_D^e$ on
the orthogonal complement $D^{\perp}$
of $D$ vanishes up to order 2 at
$q$. Hence, taking an orthonormal frame $\wt V_i,
I\wt V_i$ $(i=1,2)$ of $D$ around $q$, we
see that $X_{H_D}|_{T^*_qM}$ is of the form
$$X_{H_D}|_{T^*_qM}=f_1\wt V_1^*+f_2\wt V_2^*+f_3(I\wt V_1)^*+
f_4(I\wt V_2)^*,$$
where ${\wt V}_1^*,\dots, (I{\wt V}_2)^*$ are covectors (constant
vector fields on
$T^*_qM$) that vanish on $D^{\perp}$, and dual to
${\wt V}_1,\dots,I{\wt V}_2$; also,
$f_i$ are linear combinations of the hermitian forms
$${\wt V}_i^2+(I{\wt V}_i)^2\quad (i=1,2),\quad {\wt V}_1{\wt V}_2+
(I{\wt V}_1)(I{\wt V}_2),\quad {\wt V}_1I{\wt V}_2-{\wt V}_2I{\wt V}_1.$$
Then, since
$$0=X_{H_D}F_{t(\alpha)}=
|u(\alpha)(q)|X_{H_D}\sum_{i=1}^2({\wt V}_i^2+(I{\wt V}_i)^2),$$
we can choose ${\wt V}_1$ and ${\wt V}_2$
so that
$$\align
aX_{H_D} & =({\wt V}_1^2+(I{\wt V}_1)^2){\wt V}_2^*-({\wt V}_1{\wt V}_2+
(I{\wt V}_1)(I{\wt V}_2)){\wt V}_1^*\\
& +({\wt V}_1I{\wt V}_2-{\wt V}_2I{\wt V}_1)(I{\wt V}_1)^*,
\endalign$$
where $a$ is a non-zero constant, and the
covectors ${\wt V}_i^*,(I{\wt V}_i)^*$ are
identified with the constant vector fields on $T^*_qM$.
Hence, the zero set of
$X_{H_D}$ on $T^*_qM$ is the vector subspace defined
by ${\wt V}_1=I{\wt V}_1=0$.
By putting $V_D=\wt V_1$, (2) follows.
\qed
\enddemo

Let $\Lambda_1$ be the set of points $\lambda\in
S_q^*M$ such that
$$\dim \{(X_F)_{\lambda}\ |\ F\in\Cal F\}<n.$$

\proclaim{Lemma 3.11}
\roster
\runinitem"(1)" $S_q^*M-\Lambda_1$ is connected and dense in $S_q^*M$.
\item"(2)" For $\lambda\in S_q^*M-\Lambda_1$ the set of $t\in\R$
such that $\pi(\zeta_t\lambda)\notin M^0$ is discrete. In particular
there is a
constant $\epsilon>0$ such that $\pi(\zeta_t\lambda)\in M^{0,0}$ for $t$
satisfying
$|t|<\epsilon,t\ne 0$.
\endroster
\endproclaim
\demo{Proof} We define $V_D$ for each 1-dimensional simultaneous
eigenspace $D$ as
a unit vector in $D$. Let $\Lambda_2$ be
the set of $\lambda\in S_q^*M$ such that
$$\lambda(V_D)=\lambda(IV_D)=0$$
for some simultaneous eigenspace $D$ (of dimension 1
or 2). Then the complement of
$\Lambda_2$ in $S_q^*M$ is still connected and dense
in $S_q^*M$. We prove that
$$S_q^*M-\Lambda_2\subset S_q^*M-\Lambda_1,$$
which will indicate (1).

Let $\lambda\in S_q^*M-\Lambda_2$.
Since $\lambda\not\in \Lambda$, the previous lemma implies that
$(X_F)_{\lambda}$ $(F\in\Cal F')$ form a $k$-dimensional subspace of
$T_{\lambda}(T_q^*M)$. On the other hand, we know that
$\{F_q\ |\ F\in\Cal F\}$
is $(n-k)$-dimensional, and spanned by
$$\gather
\wt V_1^2+\wt V_2^2+(I\wt V_1)^2+
(I\wt V_2)^2\quad (D, 2\text{-dimensional}),\\
(V_D)^2+(IV_D)^2\quad (D, 1\text{-dimensional}),
\endgather$$
where $\wt V_1,\dots,I\wt V_2$ is an orthonormal basis
of $D$. Hence, it follows
that
$$\{\pi_*((X_F)_{\lambda})\ |\ F\in\Cal F\}$$
is $(n-k)$-dimensional subspace of $T_qM$,
provided $\lambda\not\in \Lambda_2$.
{}From these two facts we see that
$$\dim \{(X_F)_{\lambda}\ |\ F\in\Cal F\}=
n\qquad (\lambda\in S_q^*M-\Lambda_2).$$
Hence $\lambda\in S_q^*M-\Lambda_1$.

Now, let $\lambda\in S_q^*M-\Lambda_1$, and put $Z_F(t)=
\pi_*((X_F)_{\zeta_t\lambda})$. Then, $Z_F(t)$ are Jacobi fields
along the geodesic $c_1(t)=\pi(\zeta_t\lambda)$, and satisfy
$$g(Z_F,\nabla_{\dot c_1}Z_{\wt F})=g(\nabla_{\dot c_1}Z_F,Z_{\wt F})\qquad
(F,\wt F\in\Cal F).$$
This implies that the set of $t$ at
which $Z_F(t)=0$ for some $F\in\Cal F-\{0\}$
is discrete. Since $F_{c_1(t)}=0$ implies $Z_F(t)=0$ for each
$t$,
it follows that $\pi(\zeta_t\lambda)\in M^0$ except discrete $t$'s.
\qed
\enddemo

We now continue the proof of Proposition 3.9.
The lemmas above imply that there are
only finite number of $t$ on the interval
$(0,t_0)$ such that $c(t)\not\in M^0$.
Let $t_1,\dots,t_l$ $(t_1<\dots<t_l<t_0)$ be those points.
Then the previous lemma also implies that $c((t_1,t_2))\subset
M^{0,0}$.
Since Lemmas 3.10 and 3.11 are still valid
for the point $c(t_2)$, it also follows that
$c((t_2,t_3))\subset M^{0,0}$. Iterating
this procedure successively, we consequently see that $p'=c(t_0)\in
M^{0,0}$, which
shows the connectedness of $M^0$.

The denseness of $M^0$ in $M$ is now
clear, because any geodesic $c(t)$
emanating from the point $p$ meets $M^s$ only
at discrete values of $t$, if every
$D_i$-component of $\dot c(0)$ is nonzero. This completes
the proof of Proposition 3.9.
\qed
\enddemo
We have proved that $M^0$ and $M^1$ are
connected and dense in $M$, and $u_{\alpha}$
does not vanish everywhere on $M^0$ for any
$\alpha$. Now we shall show that
$u_{\alpha}\ne 0$ everywhere on $M$. Assume that $u_{\alpha}(q)=0$
for some
$\alpha\in\A$ and $q\in M^s$. Then $du_{\alpha}=0$ at $q$,
because $u_{\alpha}$ is
either non-positive or non-negative on $M$.
Since the functions $u_{\beta}$ are globally defined and
smooth on $M$, the vector
fields $\sgrad u_{\beta}$ are globally defined infinitesimal automorphisms
of the
K\"ahler manifold $M$. Therefore the connected component $L$
of
the set of zeros of $\sgrad u_{\alpha}$ containing
$q$ is a totally geodesic complex
submanifold of $M$ (see [8], for instance). Also,
the tangent space of $L$ at
$q$ coincides with the kernel of the hessian
of $u_{\alpha}$ at $q$.

Let $U$ be a small neighborhood of $q$,
and let $p\in U\cap M^1$ that does not belong to
any such submanifolds $L$ corresponding to the vector
fields $\sgrad u_{\beta}$
vanishing at $q$. Let $c(t)$ be the minimal
geodesic joining $q$ and $p$; $c(0)=q$,
$c(t_0)=p$. Fix $\alpha$ such that $u_{\alpha}(q)=0$, and put
$$v=v_{t(\alpha)}-\left(\sum_{i\in\alpha}h_i(p)\right) u_{\alpha}.$$
Note that the function $v$ is well-defined and
smooth on the whole $M$.

\proclaim{Lemma 3.12} $\hess v(\dot c(0),\dot c(0))=0$.
\endproclaim
\demo{Proof} Taking $U$ small enough, we may assume
that $u_{\gamma}(c(t))\ne 0$
for any $\gamma\in\A$. Hence, as we have seen
before, there is an open and dense
subset $J$ of the interval $(0,t_0]$ such that
$c(t)\in M^1$ for $t\in J$. Since
$F_{t(\alpha)}$ is semi-definite everywhere, and since $\tr F_{t(\alpha)}$
is a
non-zero multiple of $u_{\alpha}$, we have $F_{t(\alpha)}=0$ at
$q$. Hence
$F_{t(\alpha)}(\zeta_t\lambda) =0$
for any $t\in\R$ $(c(t)=\pi(\zeta_t\lambda))$, and
this implies that
$$\dot c(t)\in \sum_{\gamma\not\succeq\alpha} D_{\gamma}\qquad (t\in J).$$
Hence $\sum_{i\in\alpha}h_i(c(t))$ is constant on each connected component
of $J$.
So, by continuity we have
$$v(c(t))=u_{\alpha}(c(t))\left(\sum_{i\in\alpha}h_i(c(t))-
\sum_{i\in\alpha}h_i(p)
\right)=0$$
for all $t\in [0, t_0]$, and it therefore
follows that
$\hess v(\dot c(0),\dot c(0))=0$.
\qed
\enddemo

It is clear that a slight modification of
the initial vector $\dot c(0)\in T_qM$
does not affect the conclusion of Lemma 3.12.
Thus we have $\hess v=0$ at $q$,
which contradicts the fact that $\sgrad v$ is
the non-trivial Killing vector field on $M$.
Hence it has been shown that $u_{\alpha}$ does
not vanish everywhere on $M$ for
every $\alpha$.

We now prove the remaining part of Theorem
3.1. We have already shown that
the subbundles $D_i$ are well-defined and of $C^{\infty}$
on $M^0$. Also, (1), (2), and
(3) have been verified. Since the functions $u_{\alpha}$
are everywhere non-zero,
(4) and (6) also follows. Let $q$ be
a point in $M-M^1$. If $c(t)$ is a geodesic passing
through $q$ and a point in $M^1$, then
the set of $t\in\R$ such that
$\{Y_{c(t)}\ |\ Y\in\frak k\}$ is not $n$-dimensional is
discrete, as we have already
seen. Hence Lemma 3.5 and the argument in
the proof of Proposition 3.9 is applicable.
Thus if $q\in M^0-M^1$, there are $i(\in\alpha)$ and
$a,b\in\R$, $b\ne 0$
such that $d(v_{\alpha}(a))=0$ at $q$ and
$$\hess v_{\alpha}(a)(X,Y)=bg([x]_{D_i},[Y]_{D_i}),\qquad X,Y\in T_qM.$$
Therefore the connected component $L$ of the zeros
of the infinitesimal automorphism
$\sgrad v_{\alpha}(a)$ passing through $q$ is the complex
submanifold of codimension
one.

Now let $q\in M^s$. Then it follows that
every simultaneous eigenspace of
$\{F_q^e\ |\ F\in\Cal F \}$ is of dimension
one or two, and they are contained in
some $D_{\alpha}$. Hence (5) follows. Let $D$ be
a simultaneous eigenspace
of dimension two, and let $\alpha$, ${\alpha}'=\{i,i+1\}$ and
$c\in\R$ be as in
Lemma 3.8. From the proof of Lemma 3.8,
it easily follows that $d(v_{\alpha}(c))=0$
and $dv\ne 0$ at $q$, where
$$v=u_{\alpha}\Cal S_{|\alpha|-1}(h_j-c\,;\,j\in\alpha)$$
This implies that the exterior derivative of the
function $h_i+h_{i+1}$, which
is of $C^{\infty}$ around $q$, does not vanish
at $q$. Since $d((h_i-c)(h_{i+1}-c))=0$
at $q$, it follows that $\grad v\in D$,
and the orthogonal complement of $D$ is
contained in the kernel of $\hess v_{\alpha}(c)$ at
$q$. Since $dv_{\alpha}(c)(\sgrad v)
=0$ everywhere, it therefore follows that the connected
component $L$ of the set of
zeros of $\sgrad v_{\alpha}(c)$ passing through $q$ is
also $(n-1)$-dimensional.

Thus $M-M^1$ is equal with the union of
such hypersurfaces $L$. The local finiteness
is clear. This complete the proof of Theorem
3.1.

Also, we have just proved the following

\proclaim{Proposition 3.13} Let $q\in M^s$. Then the simultaneous
eigenspaces $D$ of
the linear endomorphisms $\{F_q^e\ |\ F\in\Cal F\}$ are
of dimension one or two. They
are smoothly extended to a neighborhood $U$ of
$q$ as the subbundles of $TM$ in the
following way: If $\dim D=1$, then $D=D_i$ on
$U\cap M^0$ for some $i$; if $\dim D=2$,
then there is $\alpha\in\A$ and $i,i+1\in\alpha$ such that
$D=D_i+D_{i+1}$ on
$U\cap M^0$. In the first case, the function
$h_i$ is smooth around $q$, and if
$i\in\alpha$,
$$h_j(q)\ne h_i(q)\qquad (j\in\alpha,\ j\ne i).$$
In the second case,
$$h_i(q)=h_{i+1}(q),\quad h_j(q)\ne h_i(q)\qquad (j\in\alpha,\ j\ne i,i+1),$$
and the functions $h_i$ and $h_{i+1}$ are not
differentiable at $q$. Moreover,
$h_i+h_{i+1}$ and $h_ih_{i+1}$ are smooth functions around $q$
such that their
exterior derivatives are non-zero and
$d((h_i-h_i(q))(h_{i+1}-h_i(q)))=0$ at $q$.
\endproclaim

\specialhead 4. Torus action and the invariant hypersurfaces
\endspecialhead
In the rest of the paper we shall
assume that the K\"ahler-Liouville manifold $M$
is {\it compact}. Let $\frak k$ be as
before, and put
$$\frak g=\frak k + I\frak k,$$
which is a commutative Lie algebra of infinitesimal
holomorphic transformations
of $M$. Let $K$ and $G$ be the
Lie transformation group of $M$ generated by
$\frak k$ and $\frak g$ respectively. Note that
$\frak g$ is naturally regarded as
a complex Lie algebra. Accordingly, $G$ is regarded
as a complex Lie group so
that the action $G\times M\to M$ is holomorphic.
In this section we shall
investigate the properties of the action of those
groups and the hypersurfaces
contained in $M-M^1$.

We first prove the following

\proclaim{Proposition 4.1} \roster
\runinitem"(1)" $G$ preserves $M^1$ and each hypersurface $L\subset
M-M^1$.
\item"(2)" The action of $G$ on $M^1$ is
simply transitive.
\endroster
\endproclaim
\demo{Proof} To prove (1), recall that $L$ is
a connected component of the set of
zeros of $\sgrad v\in\frak k$
for some $v\in\Cal V$. Therefore, every $Y\in\frak k$
is tangent to $L$. Since
$L$ is a complex submanifold, any $Z\in\frak g$
is also tangent to $L$. Thus $G$
preserves $L$. Since $M^1$ is the complement of
the union of such hypersurfaces,
it is also preserved by $G$. To prove
(2), note that $\{Y_p\ |\ Y\in \frak g\}$
is real $2n$-dimensional for every point $p\in M^1$.
Hence the $G$-action on $M^1$ is
transitive. Suppose $gp=p$ for some $g\in G$ and
$p\in M^1$. Then $gq=q$ for every
$q\in M^1$, because $G$ is abelian. Hence $g$
should be the identity transformation of
$M$ by continuity.
\qed
\enddemo

\proclaim{Theorem 4.2} The Lie group $K$ is isomorphic
to
$$U(1)^n=U(1)\times\dots\times U(1)\quad (n\text{ times}),$$
where $U(1)$ is the group of unit complex
numbers. Also, $G$ is isomorphic to
$(\C^{\times})^n$ as complex Lie group, where $\C^{\times}$
is the multiplicative group of non-zero complex numbers.
\endproclaim

To prove this theorem we need several lemmas.
Let $L$ be a complex hypersurface
contained in $M-M^1$. As observed before, there is
$v_{\alpha}(c)\in\Cal V$
such that $L$ coincides with a connected component
of the set
of zeros of $\sgrad v_{\alpha}(c)$, and $v_{\alpha}(c)$ vanishes
on $L$.
In this case we shall call $v_{\alpha}(c)$ {\it
the function that determines} $L$.

\proclaim{Lemma 4.3} Let $L$ and $v_{\alpha}(c)$ be as
above. Then the vector
field $\sgrad v_{\alpha}(c)$ generates a circle action on
$M$.
\endproclaim
\demo{Proof} Put $Y=\sgrad v_{\alpha}(c)$, and let $\ad Y$
be the linear
endomorphism of $T_pM$ $(p\in L)$ given by the
formula
$$(\ad Y)(X)= [Y,X](=-\nabla_XY),\qquad X\in T_pM.$$
Clearly, $\ad Y$ preserves the normal space $N_pL$.
Since $N_pL$ is complex $1$-dimensional and $Y$ is
the infinitesimal isometry,
there is $a\in\R$ such that
$$(\ad Y)(X)=aIX\quad \text{for any }\ X\in N_pL.$$
We claim that $a$ is independent of $p\in
L$. In fact, let $Z\in T_pL$, and let $X$
be a unit normal vector field along $L$.
Then
$$-\nabla_Z(\nabla_XY)=(Za)IX+ aI\nabla_ZX.$$
Since $Y$ is a Killing vector field, we
have
$$\nabla_Z(\nabla_XY)=\nabla_{\nabla_ZX}Y+
\frac 12 (R(Z,X)Y-R(Y,Z)X+R(X,Y)Z).$$
Hence $Za=0$, which shows that $a$ is constant
on $L$.

The lemma is now clear, because the linear
isotropy action of the 1-parameter
group $\phi_t$ generated by $Y$ has the least
period $2\pi/|a|$ at every point
$p\in L$, and so is $\phi_t$ itself via
the exponential mapping Exp$\,:NL\to M$.
\qed
\enddemo

\proclaim{Lemma 4.4} Fix $\alpha\in\A$, and take a constant
$b$ so that
$h_i+b>0$ on $M$ for any $i\in\alpha$. Put
$$v=u_{\alpha}\prod_{i\in\alpha}(h_i+b)\in\Cal V,$$
and let $p\in M$ be a point where
$|v|$ attains the maximum. Then the functions
$h_j$ are smooth around $p$, and $dh_j=0$ at
$p$ for any $j$ such that
$j\in\beta$, $\beta\pc\alpha$.
\endproclaim
\demo{Proof} If $p\in M^0$, then the assertion is
clear, because $dh_i(D_j)=0$
for any $i,j$ such that $i\ne j$. Similarly,
if there is no 2-dimensional
simultaneous eigenspace of $\{F_p^e\ |\ F\in\Cal F\}$ contained
in
$\sum_{\beta\pc\alpha}D_{\beta}$, the assertion is also clear.
Now, assume that $p\in M^s$ and there are
$\beta\pc\alpha$ and $i, i+1\in\beta$
such that $h_i(p)=h_{i+1}(p)$. Let $D$ be the subbundle
of $TM$ defined on a
neighborhood $U$ of $p$ that coincides with $D_i+D_{i+1}$
on $M^0\cap U$. Then, putting
$$X=\grad (h_i+h_{i+1}),$$
we have $X\in D$ and $X\ne 0$ by
virtue of Proposition 3.13. Since
$$d((h_i-h_i(p))(h_{i+1}-h_i(p)))=0$$
at $p$, it therefore follows that
$$d((h_i+b)(h_{i+1}+b))(X)\ne 0\qquad \text{at }\ p.$$
However, this implies that $dv(X)\ne 0$ at $p$,
a contradiction.
\qed
\enddemo

\demo{Proof of Theorem 4.2} Fix $\alpha\in\A$, and let
$p\in M$ be as in
Lemma 4.4. Clearly the function $v_{\alpha}(h_i(p))$
$(i\in\alpha)$ determines
a hypersurface $L_i$ contained in $M-M^1$ that pass
through $p$. Hence the vector fields
$$Y_i=\sgrad v_{\alpha}(h_i(p))\qquad (i\in\alpha)$$
generate circle actions on $M$ by Lemma 4.3.
Executing the same procedure for all
$\alpha$ we thus obtain $n$ vector fields $Y_i\in
\frak k$ each of which
generates a circle action on $M$. As is
easily seen, those vector fields form a
basis of $\frak k$. Hence $K$ is compact,
and is isomorphic to $U(1)^n$.

Now, fix a point $q\in M^1$, and let
$\phi_i:\frak g\to \R$ be the mapping defined by
$$\phi_i(Z)= v_{\alpha}(h_i(p))((\exp Z)q).$$
Put $\Phi=(\phi_i):\frak g\to \R^n$. Then we have $\Phi(Z+Y)=\Phi(Z)$
for any $Y\in\frak k$
and $Z\in\frak g$, and
$$\frac{\partial}{\partial t_j}\phi_i(\sum_k t_kIY_k)
=-g(Y_i,Y_j).$$
{}From this formula it easily follows that the
inner product of two vectors
$$\Phi(\sum_k t_kIY_k)-
\Phi(\sum_k s_kIY_k)\quad\text{and}\quad t-s=(t_k-s_k)$$
in $\R^n$ is negative, provided $t\ne s$. Hence
$\Phi|_{I\frak k}$ is injective.

These facts indicate that the kernel of the
homomorphism $\exp:\frak g\to G$ is
equal to that of $\exp:\frak k\to K$. This
proves the latter half of the theorem.
\qed
\enddemo

Our next goal is to determine all the
hypersurfaces contained in $M-M^1$. For this
purpose we need deeper information on the fundamental
functions $h_i$. Let $L$ be a
hypersurface in $M-M^1$ determined by the function $v_{\alpha}(c)\in\Cal
V$.
The following lemma is easy.

\proclaim{Lemma 4.5} Let $p\in L$. Then there are
neighborhoods $W$ and $U$ of $p$ in
$L$ and $M$ respectively, a neighborhood $V$ of
$0$ in $\C=\{(z)\}$, and a holomorphic
diffeomorphism $\phi: W\times V\to U$ such that
\roster
\item"(1)" $\phi(q,0)=q$ for any $q\in W$,
\item"(2)" $\phi_*z\frac{\partial}{\partial z}=(2a)^{-1}(-IY+\sqrt{-1}Y)$,
\endroster
where $Y=\sgrad v_{\alpha}(c)$ and $a$ is the eigenvalue
of $I\circ\ad Y$ on $NL$.
\endproclaim

Note that the hessian of $v_{\alpha}(c)$ on the
normal bundle $NL$ is equal to $a$ times the
metric $g$, $a$ being the constant in Lemma
4.5. Hence there is a neighborhood
$U$ of $L$ such that $v_{\alpha}(c)> 0$ (resp.
$< 0$) on $U-L$ if $a>0$ (resp. $a<0$).

\proclaim{Lemma 4.6} $v_{\alpha}(c)> 0$ (resp. $< 0$) on
$M-L$ if
$a>0$ (resp. $a<0$).
\endproclaim
\demo{Proof} Suppose $a>0$. Let $\psi_t$ be the one-parameter
group of transformations
generated by $-\grad v_{\alpha}(c)$. Then the previous lemma
implies that $\psi_t(q)$
converges to a point in $L$ as $t\to
\infty$ for any $q\in U$, provided $U\ (\supset L)$
small enough. Now, fix $q_0\in U$, and let
$q\in M^1$ be an arbitrary point. Then there
is $g\in G$ such that $q=gq_0$. Since $\psi_t(q)=g\psi_t(q_0)$
and $gL=L$, it follows that
$\psi_t(q)\in U$ for sufficiently large $t$. Hence we
have
$$v_{\alpha}(c)(q)>v_{\alpha}(c)(\psi_t(q))>0.$$
Let $b>0$ be the minimal value of the
function $v_{\alpha}(c)$ on the boundary of $U$, and
let $q'$ be a point in $M-(M^1\cup U)$.
Then one can take a point $q\in M^1- U$ arbitrary
near $q'$. Since $v_{\alpha}(c)(q)>b$, it follows that $v_{\alpha}(c)(q')\ge
b$, proving
the lemma. The case $a<0$ is similar.
\qed
\enddemo

\proclaim{Lemma 4.7} Let $\alpha\in\A$, and let $i,i+1\in\alpha$. Then
$$\min h_i\ge \max h_{i+1},$$
where the minimum and the maximum are taken
on $M$.
\endproclaim
\demo{Proof} Suppose that $h_i$ takes its minimal value
$c$ at $p\in M$. In this case we
have $h_{i-1}(p)>c$ if $i-1\in\alpha$. In fact, if $i-1\in\alpha$
and $h_{i-1}(p)=c$, then
$c$ is also the minimal value of $h_{i-1}$.
Hence $d(h_{i-1}+h_i)=0$ at $p$, contradicting
Proposition 3.13. Then, there are two cases: (1)
$h_{i+1}(p)<c$, or (2) $h_i(p)
=h_{i+1}(p)=c$. If the case (1) occurs, then the
function $h_i$ is smooth around $p$. Hence
the function $v_{\alpha}(c)$ determines a hypersurface $L\subset M-M^1$.
Since $h_{i+1}
\ne c$ on $M-L$ by the previous lemma,
we have $\max h_{i+1}\le c$.

Now, suppose that the case (2) occurs. Then
by virtue of Proposition 3.13 the function
$v_{\alpha}(c)$ again determines a hypersurface $L\subset M-M^1$. Proposition
3.13 also
implies that there is a point $q\in L$
near $p$ such that $h_{i+1}(q)<c$ and $h_i(q)=c$.
Hence we again conclude that $\max h_{i+1}\le c$.
\qed
\enddemo

\proclaim{Proposition 4.8} A function in $\Cal V$ of
the form $v_{\alpha}(c)$ determines
a hypersurface $L$ in $M-M^1$ if and only
if $c$ is the maximal or minimal value of $h_i$
for some $i\in\alpha$.
\endproclaim
\demo{Proof} The ``\,if\,'' part has been indicated in
the proof of the previous proposition.
Now, suppose that a function $v_{\alpha}(c)$ determines a
hypersurface $L$ in $M-M^1$, and
let $p\in L$. Then there is $i\in\alpha$ such
that $h_i(p)=c$. Since $h_i\ne c$ on $M-L$,
$c$ should be the maximal value or the
minimal value of $h_i$.
\qed
\enddemo

\proclaim{Corollary 4.9} Suppose that the function $h_i$ is
smooth around a point $p\in M$,
and $dh_i=0$ at $p$. Then $h_i$ takes its
maximum or minimum at $p$.
\endproclaim
\demo{Proof} The assumption implies that the function $v_{\alpha}(h_i(p))$
determines a
hypersurface in $M-M^1$ $(i\in\alpha)$. Thus the corollary follows
from Proposition 4.8.
\qed
\enddemo

We now prove the following theorem.

\proclaim{Theorem 4.10} For any $\alpha\in\A$ and $i,i+1\in\alpha$,
$$\min h_i=\max h_{i+1}.$$
\endproclaim

The following corollary is an immediate consequence of
Theorem 4.10 and Proposition 4.8.

\proclaim{Corollary 4.11} Put
$$\align
& c_{\alpha,0}=\max h_{s(\alpha)},
\qquad c_{\alpha,|\alpha|}= \min h_{t(\alpha)},\\
& c_{\alpha,\nu}=\min h_{s(\alpha)+\nu-1}=\max h_{s(\alpha)+\nu}\quad
(1\le\nu\le |\alpha|-1),
\endalign$$
and let $L_{\alpha,\nu}$ be the hypersurface in $M-M^1$
determined by $v_{\alpha}
(c_{\alpha,\nu})$. Then the hypersurfaces
$L_{\alpha,\nu}$ are mutually distinct,
and the set
$$\{L_{\alpha,\nu}\ |\ \alpha\in\A,\ 0\le\nu\le |\alpha|\}$$
coincides with the set of all closed hypersurfaces
contained in $M-M^1$.
\endproclaim

Let us recall that the fundamental functions $h_i$
$(i\in\alpha)$ and the conjunction
constants $e_{\alpha\beta}$ $(\alpha\pcc\beta)$ may be replaced with
$$k_{\alpha}h_i-l_{\alpha}\quad (k_{\alpha}>0)\quad\text{or}\quad
k_{\alpha}h_{t(\alpha)+s(\alpha)-i}-l_{\alpha}\quad (k_{\alpha}<0)$$
and $k_{\alpha}e_{\alpha\beta}+l_{\alpha}$
respectively, where $k_{\alpha}\ne 0$ and $l_{\alpha}$ are constants.
Hence it is always
possible to choose $h_i$ and $e_{\alpha\beta}$ so that
$$1=c_{\alpha,0}>c_{\alpha,1}>\dots>c_{\alpha,|\alpha|}=0. \tag 4.1$$
In this case we also have
$$e_{\alpha\beta}>0\quad \text{or}\quad e_{\alpha\beta}<-1. \tag 4.2$$
Under this condition the only possible alternative choice
of $h_i$ $(i\in\alpha)$ and
$e_{\alpha\beta}$ are given by
$$h'_i=1-h_{t(\alpha)+s(\alpha)-i}\quad (i\in\alpha),\qquad e'_{\alpha\beta}=
-1-e_{\alpha\beta}.$$
In the rest of the paper we shall
always assume that {\it the fundamental functions $\{h_i\}$
and the conjunction constants $\{e_{\alpha\beta}\}$ are chosen so
that the condition (4.1)
is satisfied}. Also, we shall call $\{c_{\alpha,\nu}\}$ {\it
the fundamental constants}.

\demo{Proof of Theorem 4.10} Assume that $i, i+1\in
\alpha$ and
$$\min h_i=c_2> c_3 =\max h_{i+1},$$
and put $c_1=\max h_i$, $c_4=\min h_{i+1}$. Let $L_{\mu}$
be the hypersurface in $M-M^1$
determined by $v_{\alpha}(c_{\mu})$, and put
$$Y_{\mu}=\sgrad v_{\alpha}(c_{\mu})\qquad (\mu=1,\dots,4).$$
Also, put
$$b_j=\cases \max h_j\quad (j<i)\\
\min h_j\quad (j>i+1)
\endcases$$
for $j\in\alpha$, $j\ne i, i+1$.

\proclaim{Lemma 4.12} There are four points $p_{13}$, $p_{14}$,
$p_{23}$, $p_{24}$ such
that
\roster
\item"(1)" $p_{\mu\nu}\in L_{\mu}\cap L_{\nu}\qquad (\mu=1,2,\ \nu=3,4)$,
\item"(2)" $h_j$ is smooth and $dh_j=0$ at the
four points for every $\beta\pc\alpha$
and every $j\in\beta$,
\item"(3)" $h_j=b_j$ at the four points for every
$j\in\alpha$, $j\ne i, i+1$.
\endroster
\endproclaim
\demo{Proof} First we show that $L_{\mu}\cap L_{\nu}\ne\emptyset$ $(\mu=1,2,\
\nu=3,4)$. Let $b$ the maximal value of the
function $h_i$ on the hypersurface $L_3$,
and suppose $h_i(q)=b$, $q\in L_3$. Then the similar
argument as the proof of
Lemma 4.4 implies that the function $h_i$ is
smooth and $dh_i=0$ at $q$. Hence
$b=c_1$ or $c_2$ by Corollary 4.9. Since the
case $b=c_2$ contradicts the choice of $q$,
we have $b=c_1$. Thus $L_1\cap L_3\ne \emptyset$. The
other cases are similar.

Now, choose a constant $d$ such that $c_2>d>c_3$,
and let $p_{\mu\nu}\in L_{\mu}
\cap L_{\nu}$ be a point where the function
$|v_{\alpha}(d)|$, restricted to $L_{\mu}
\cap L_{\nu}$, takes the maximal value $(\mu=1,2,\ \nu=3,4)$.
Then the similar argument
as above clearly indicates that the conditions (2)
and (3) are satisfied.
\qed
\enddemo

Let us consider the following identity (the decomposition
to linear fractions):
$$\align
\frac{v_{\alpha}(\lambda)}{\prod_{j\in{\alpha}'}
(b_j-\lambda)\prod_{\mu=1}^3(c_{\mu}-
\lambda)} = \sum_{j\in{\alpha}'}\frac{1}{b_j-\lambda}\frac{v_{\alpha}(b_j)}
{\prod_{k\in{\alpha}'\atop k\ne j}(b_k-b_j)\prod_{\mu=1}^3(c_{\mu}-b_j)}\\
+\sum_{\mu=1}^3\frac{1}{c_{\mu}-\lambda}\frac{v_{\alpha}(c_{\mu})}
{\prod_{j\in{\alpha}'}(b_j-c_{\mu})
\prod_{1\le\nu\le 3\atop \nu\ne\mu}(c_{\nu}-c_{\mu})},
\tag 4.3 \endalign$$
where ${\alpha}'=\alpha-\{i,i+1\}$. Multiplying both sides
by $-\lambda$ and
taking the
limit $\lambda\to\infty$, we obtain the identity:
$$\align
u_{\alpha} & = \sum_{j\in{\alpha}'}\frac{v_{\alpha}(b_j)}
{\prod_{k\in{\alpha}'\atop k\ne j}(b_k-b_j)
\prod_{\mu=1}^3(c_{\mu}-b_j)}\\
& +\sum_{\mu=1}^3\frac{v_{\alpha}(c_{\mu})}
{\prod_{j\in{\alpha}'}(b_j-c_{\mu})
\prod_{1\le\nu\le 3\atop \nu\ne\mu}(c_{\nu}-c_{\mu})}.
\endalign$$

This formula gives the linear dependence of the
skew-gradient vector fields $Y_1$, $Y_2$,
$Y_3$, $\sgrad v_{\alpha}(b_j)$ $(j\in{\alpha}')$,
and $\sgrad u_{\alpha}$. Let
$d_{\mu}$ be
the positive number so that the least period
of the one-parameter group generated by
$Y_{\mu}$ is $2\pi/d_{\mu}$ $(\mu=1,\dots,4)$. We now consider the
linear isotropy
action of the one-parameter groups generated by those
vector fields at the point $p_{13}$
(Note that those vector fields vanish at $p_{13}$).
There we have the decomposition
$$T_{p_{13}}M= N_{p_{13}}L_1\oplus N_{p_{13}}
L_3\oplus T_{p_{13}}(L_1\cap L_3).$$
Clearly, the linear isotropy action of the one-parameter
groups generated by $Y_3$,
$\sgrad u_{\alpha}$, and $\sgrad v_{\alpha}
(b_j)$ $(j\in{\alpha}')$ are trivial
on
$N_{p_{13}}L_1$. Hence the linear endomorphism
$$\exp\left(\frac{tY_2}{\prod_{j\in{\alpha}'}
(b_j-c_2)\prod_{1\le\mu\le 3\atop \mu\ne 2}
(c_{\mu}-c_2)}\right)
\exp\left(\frac{tY_1}{\prod_{j\in{\alpha}'}
(b_j-c_1)\prod_{1\le\mu\le 3\atop \mu\ne 1}
(c_{\mu}-c_1)}\right)\tag 4.4$$
of $N_{p_{13}}L_1$ is trivial for all $t\in\R$. Substituting
$$t=2\pi d_2^{-1}(c_1-c_2)(c_3-c_2)\prod_{j\in{\alpha}'}(b_j-c_2)$$
in (4.4), we conclude that
$$2\pi d_2^{-1}(c_1-c_2)(c_3-c_2)\prod_{j\in{\alpha}'}(b_j-c_2)=
m\cdot 2\pi d_1^{-1}(c_2-c_1)(c_3-c_1)
\prod_{j\in{\alpha}'}(b_j-c_1).\tag 4.5$$
The similar argument at the point $p_{23}$ gives
the formula
$$2\pi d_1^{-1}(c_2-c_1)(c_3-c_1)\prod_{j\in{\alpha}'}(b_j-c_1)=
m'\cdot 2\pi d_2^{-1}(c_1-c_2)(c_3-c_2)\prod_{j\in{\alpha}'}(b_j-c_2)$$
for some $n'\in\Z$. Therefore $m=\pm 1$. Replacing $c_3$
by $c_4$ in the
formula (4.3), and considering the linear isotropy actions
at $p_{14}$ and $p_{24}$, we
also obtain the formula
$$2\pi d_2^{-1}(c_1-c_2)(c_4-c_2)\prod_{j\in{\alpha}'}(b_j-c_2)=
\pm 2\pi d_1^{-1}(c_2-c_1)(c_4-c_1)\prod_{j\in{\alpha}'}(b_j-c_1).\tag 4.6$$
{}From (4.5) and (4.6) we have
$$(c_3-c_2)(c_4-c_1)=\pm (c_4-c_2)(c_3-c_1),$$
which contradicts the inequality; $c_1>c_2>c_3>c_4$. This completes the
proof of
Theorem 4.10.
\qed
\enddemo

In the rest of the section we shall
observe the detail of the action of $G$ on $M$. In
particular we shall show that $M$ is a
toric variety. Let $c_{\alpha,\nu}$ and
$L_{\alpha,\nu}$ be as in Corollary 4.11, and put
$$\J=\{(\alpha,\nu)\ |\ \alpha\in\A,\ 0\le \nu\le |\alpha|\}.$$
Let $d_{\alpha,\nu}$ be the non-zero constant such that
$$\hess v_{\alpha}(c_{\alpha,\nu})(X,X)=d_{\alpha,\nu}g(X,X)\quad
(X\in NL_{\alpha,\nu}),$$
and put
$$Y_{\alpha,\nu}=d_{\alpha,\nu}^{-1}
\sgrad v_{\alpha}(c_{\alpha,\nu})\in\frak k.$$
Also, let $\I$ be the set of sections
of the mapping $\J\to\A$ $((\alpha,\nu)
\mapsto \alpha)$, that is, $\iota\in\I$ is an assignment
of an index $\iota(\alpha)$
$(0\le \iota(\alpha)\le |\alpha|)$ to each $\alpha\in\A$. Put
$$\J(\iota)=\{(\alpha,\nu)\in \J\ |\ \nu\ne\iota(\alpha)\}
\qquad (\iota\in\I).$$

\proclaim{Lemma 4.13} \roster
\runinitem"(1)" For any $\alpha$,
$\cap_{\nu=0}^{|\alpha|} L_{\alpha,\nu}=\emptyset$.
\item"(2)" For any $\iota\in\I$, $\cap_{(\alpha,\nu)\in\J(\iota)}
L_{\alpha,\nu}\ne\emptyset$.
\endroster
\endproclaim
\demo{Proof} (1) Let $p\in\cap_{\nu=0}^{|\alpha|}
L_{\alpha,\nu}$. First, since
$p\in L_{\alpha,0}$, it follows that $h_{s(\alpha)}(p)
=c_{\alpha,0}$. Next, the
condition that
$p\in L_{\alpha,1}$ implies that $h_{s(\alpha)+1}(p)
=c_{\alpha,1}$, and so on.
Consequently,
we have $h_{s(\alpha)+\nu}(p)=c_{\alpha,\nu}$
$(0\le \nu\le |\alpha|-1)$ from the
condition
that $p\in\cap_{\nu=0}^{|\alpha|-1} L_{\alpha,\nu}$. However, since $p\in
L_{\alpha,|\alpha|}$, we also have $h_{t(\alpha)}(p)
=c_{\alpha,|\alpha|}$, a contradiction.

(2) Put
$$v=\prod_{\alpha\in\A}\prod_{i\in\alpha}(h_i-c_{\alpha,\iota(\alpha)}),$$
and let $p\in M$ be a point where
the function $|v|$ takes the maximum. Then, as in the
proof of Lemma 4.4, we see that $p\in
M^0$ and $dh_i=0$ at $p$ for every $i$. Hence
$$\{h_i(p)\ |\ i\in\alpha\}=
\{c_{\alpha,\nu}\ |\ 0\le\nu\le |\alpha|,\ \nu\ne\iota(\alpha)\}.$$
This indicates that $p\in L_{\alpha,\nu}$ for every
$(\alpha,\nu)\in\J(\iota)$.
\qed
\enddemo

Let $U$ be the neighborhood of $L_{\alpha,\nu}$ given
in the proof of Lemma 4.6, and put
$$U_{\alpha,\nu}=\cup_{t\in\R}\psi_t(U),$$
where $\psi_t$ is the one-parameter group generated by
$-IY_{\alpha,\nu}$. It is clear
from Lemma 4.6 that $U_{\alpha,\nu}$ is $G$-invariant, and
$\psi_t(q)$ converges to
a point in $L_{\alpha,\nu}$ as $t\to -\infty$ for
any $q\in U_{\alpha,\nu}$. Define the
mapping $\rho_{\alpha,\nu}: U_{\alpha,\nu}\to L_{\alpha,\nu}$ by
$$\rho_{\alpha,\nu}(q)=\lim_{t\to -\infty}\psi_t(q).$$

\proclaim{Proposition 4.14} \roster
\runinitem"(1)" $U_{\alpha,\nu}=M-\cap_{0\le\mu\le |\alpha|\atop \mu\ne\nu}
L_{\alpha,\nu}$.
\item"(2)" $\rho_{\alpha,\nu}: U_{\alpha,\nu}
\to L_{\alpha,\nu}$ is the holomorphic
fibre bundle with typical fibre $\C$. Also, $\rho_{\alpha,\nu}:
U_{\alpha,\nu}
-L_{\alpha,\nu}\to L_{\alpha,\nu}$ is the principal $\C^{\times}$-bundle.
\endroster
\endproclaim
\demo{Proof} Put
$$S=\cap_{0\le\mu\le |\alpha|\atop \mu\ne\nu}L_{\alpha,\nu}.$$
Since $S$ is $G$-invariant, and $S\cap L_{\alpha,\nu}=\emptyset$, it
follows that
$U_{\alpha,\nu}\subset M-S$. Now, we shall show the reversed
inclusion. Put
$$v=\prod_{i\in\alpha}(h_i-c_{\alpha,\nu})
=u_{\alpha}^{-1}v_{\alpha}(c_{\alpha,\nu}).$$
Then, the function $|v|$ is positive on $M-L_{\alpha,\nu}$
and takes its maximal value
$$\prod_{0\le\mu\le |\alpha|\atop \mu\ne\nu}|c_{\alpha,\mu}-c_{\alpha,\nu}|$$
on $S$. Also, we have
$$\frac{d}{dt}\, |v(\psi_t(q))|=|\grad v|^2\cdot \frac{d_{\alpha,\nu}^{-1}
v_{\alpha}(c_{\alpha,\nu})}{|v|} >0$$
for every $q\in U_{\alpha,\nu}-L_{\alpha,\nu}$ and $t\in\R$.

Let $q\in M-S$. Let $p_j=\psi_{t_j}(q)$ $(0\ge t_j\to -\infty)$
be a
converging sequence, and let $p\in M$ be its
limit point. Since $|v(p)|<|v(q)|$, it follows that
$p\in M-S$. Also, we have $(d|v|)_p=0$. As is
easily seen, the set of critical
points of the function $|v|$ is equal to
$L_{\alpha,\nu}\cup S$. Hence it follows that
$p\in L_{\alpha,\nu}$. This implies that $\psi_{-t}(q)\in U$ for
sufficiently large
$t$. Thus $q\in U_{\alpha,\nu}$, completing the proof of
(1). (2) is the immediate consequence of
Lemma 4.5.
\qed
\enddemo

\proclaim{Proposition 4.15} Let $\J_0$ be a subset of
$\J$ such that
$\cap_{(\alpha,\nu)\in\J_0} L_{\alpha,\nu}\ne \emptyset$. Let $\rho_0$ be the
composition of all $\rho_{\alpha,\nu}$, $(\alpha,\nu)\in\J_0$, and put
$$S_0=\cup_{\alpha\in \A}
\cap_{0\le\nu\le |\alpha|\atop (\alpha,\nu)\not\in\J_0}
L_{\alpha,\nu}.$$
Then the mapping $\rho_0: M-S_0\to
\cap_{(\alpha,\nu)\in\J_0} L_{\alpha,\nu}$ is
well-defined, and is a fibre bundle with typical
fibre $\C^k$, where $k=\#\J_0$.
In particular, $\cap_{(\alpha,\nu)\in\J_0} L_{\alpha,\nu}$ is connected.
\endproclaim
\demo{Proof} We shall prove the proposition by induction
on $k$. Let $\J_0$, $\rho_0$,
and $S_0$ be as in the statement. Suppose
that $(\beta,\mu)\not\in\J_0$, and put
$$\J_1=\J_0\cup \{(\beta,\mu)\}.$$
We assume that $\cap_{(\alpha,\nu)\in\J_1}
L_{\alpha,\nu}\ne\emptyset$. $\rho_1$
and $S_1$ are similarly defined. We then have
the following commutative diagram:
$$\CD
M-S_1 @>{\rho_0}>> \cap_{(\alpha,\nu)\in\J_0}L_{\alpha,\nu}-
\cap_{0\le\nu \le |\beta|\atop \nu\ne\mu}L_{\beta,\nu}\\
@V{\rho_{\beta,\mu}}VV @VV{\rho_{\beta,\mu}}V\\
L_{\beta,\mu}-S_0 @>>{\rho_0} > \cap_{(\alpha,\nu)\in \J_1}L_{\alpha,\nu}
\endCD $$
Hence $\rho_1=\rho_0\circ \rho_{\beta,\mu}$ is well-defined.
{}From the induction assumption and the previous proposition,
the rows and the columns in
the diagram are fibre bundles with fibre $\C^k$
and $\C$ respectively. Let $q\in M-S_1$,
and let $p=\rho_1(q)$. Then, there is a neighborhood
$U$ of $p$ in $\cap_{(\alpha,\nu)\in
\J_1}L_{\alpha,\nu}$ such that $\rho_{\beta,\mu}^{-1}(U)$
and $\rho_0^{-1}(U)$ are
isomorphic to $U\times\C$ and $U\times \C^k$ respectively. Moreover,
for each $r\in
\rho_0^{-1}(p)$ the mapping
$$\rho_0: \rho_{\beta,\mu}^{-1}(r)\to\rho_{\beta,\mu}^{-1}(p)$$
is an isomorphism, because it commutes with the
$\C^{\times}$-action generated by
$Y_{\beta,\mu}$ and $IY_{\beta,\mu}$. Therefore the mapping
$$\rho_1^{-1}(U)\to U\times \C\times\C^k$$
defined by $\rho_0$ and $\rho_{\beta,\mu}$ is isomorphic, and
it gives the local
triviality of $\rho_1$.
\qed
\enddemo

The following corollary is an immediate consequence of
Proposition 4.15.

\proclaim{Corollary 4.16} Let $\iota\in\I$, and fix a point
$p_0\in M^1$. Then there is
a $G$-equivariant holomorphic isomorphism from $M-\cup_{\alpha\in\A}
L_{\alpha,\iota(\alpha)}$ to
$\C^n=\{(z_{\alpha,\nu}^{(\iota)};\ (\alpha,\nu)\in
\J(\iota))\}$ such that $p_0$ corresponds to the point
given by $z_{\alpha,\nu}^{(\iota)}=1$
for every $(\alpha,\nu)\in\J(\iota)$. Here the $G$-action on $\C^n$
is given by
$$\exp\left(\sum_{(\beta,\mu)\in J(\iota)}(-t_{\beta,\mu}IY_{\beta,\mu}+
s_{\beta,\mu}Y_{\beta,\mu})\right)
(z_{\alpha,\nu}^{(\iota)})=\left(e^{t_{\alpha,\nu}+\sqrt{-1}
s_{\alpha,\nu}}z_{\alpha,\nu}^{(\iota)}\right).$$
\endproclaim

Let $\Gamma$ be the lattice in $\frak k$
such that $2\pi\Gamma$ is equal to the kernel
of the homomorphism $\exp:\frak k\to K$ of the
abelian groups. Clearly $Y_{\alpha,\nu}
\in\Gamma$ for every $(\alpha,\nu)\in\J$.

\proclaim{Proposition 4.17} For any $\iota\in\I$, the elements
$Y_{\alpha,\nu}\ (\,(\alpha,\nu)\in\J(\iota)\,)$
form a $\Z$-basis of $\Gamma$.
\endproclaim
\demo{Proof} Fix $\iota\in\I$. Then, by virtue of Lemma
4.13 (2) there is a point $p$
such that $p\in L_{\alpha,\nu}$ for every $(\alpha,\nu)\in\J(\iota)$.
Since the associated endomorphisms $\ad Y_{\alpha,\nu}$ of $T_pM$
are linearly
independent, it follows that $Y_{\alpha,\nu}$
$(\,(\alpha,\nu)\in\J(\iota)\,)$ form a
$\R$-basis of $\frak k$. Let $Y$ be any
element of $\Gamma$, and let
$$Y=\sum_{(\alpha,\nu)\in\J(\iota)} a_{\alpha,\nu}Y_{\alpha,\nu},
\quad a_{\alpha,\nu}\in\R.$$
We recall that the linear isotropy action of
the one-parameter
group $\exp(tY_{\alpha,\nu})$ on
$N_pL_{\alpha',\nu'}$ is trivial if $(\alpha,\nu)\ne
(\alpha',\nu')$, and has the least period $2\pi$ if
$(\alpha,\nu)=(\alpha',\nu')$.
Since $\exp(2\pi Y)$ is the identity, it thus
follows that $a_{\alpha,\nu}\in\Z$ for all
$(\alpha,\nu)\in\J(\iota)$.
\qed
\enddemo

\proclaim{Theorem 4.18} $M$ is a toric variety with
respect to the action of $G$.
\endproclaim
\demo{Proof} By virtue of Corollary 4.16, $M$ is
covered by the open sets
$U_{\iota}=M-\cup_{\alpha\in\A}
L_{\alpha,\iota(\alpha)}$ $(\iota\in\I)$ each of which is
holomorphically isomorphic to $\C^n$. As is easily seen,
the coordinate change on $U_{\iota}
\cap U_{\iota'}$ is given by Laurent monomials whose
exponents are equal to the
coefficients of the base change of $\Gamma$: $Y_{\alpha,\nu}$
$((\alpha,\nu)\in\J(\iota))$
to $Y_{\alpha,\nu}$ $((\alpha,\nu)\in\J(\iota'))$. Hence $M$ is an algebraic
variety.
Corollary 4.16 also indicates that the action of
the ``algebraic torus'' $G$ on $M$ is
algebraic.
\qed
\enddemo

The next several propositions will give the information
on the structure of the toric variety $M$.

\proclaim{Proposition 4.19} For each $\alpha\in\A$ the value
$$\frac{d_{\alpha,\nu}}{\prod_{0\le\mu\le |\alpha|\atop \mu\ne\nu}
(c_{\alpha,\mu}-c_{\alpha,\nu})}$$
does not depend on $\nu$ $(0\le\nu\le |\alpha|)$.
\endproclaim
\demo{Proof} In the same way as the proof
of Lemma 4.12 we have
$$u_{\alpha}=\sum_{\nu=0}^{|\alpha|}\frac{v_{\alpha}(c_{\alpha,\nu})}
{\prod_{0\le\mu\le |\alpha|\atop \mu\ne\nu}
(c_{\alpha,\mu}-c_{\alpha,\nu})}.$$
Taking the skew gradient vector fields of both
sides, we then have
$$\sgrad u_{\alpha}=\sum_{\nu=0}^{|\alpha|}
\frac{d_{\alpha,\nu}}{\prod_{0\le\mu
\le |\alpha|\atop \mu\ne\nu} (c_{\alpha,\mu}
-c_{\alpha,\nu})} Y_{\alpha,\nu}.\tag 4.7$$

We claim here that $\sgrad u_{\alpha}$ is written
as a linear combination of $Y_{\gamma,
\nu}$ $(\gamma\pcc\alpha,\ 1\le\nu\le |\gamma|)$. In fact, it is
clear from the definition
of $u_{\alpha}$ that $\sgrad u_{\alpha}$ is written as
a linear combination of
$Y_{\beta,\nu}$ $(0\le\nu\le |\beta|)$, where $\beta$ is the maximal
element of the
totally ordered set $\{\gamma\in\A\ |\ \gamma\pcc\alpha\}$. Hence, by
the formula (4.7)
(replaced $\alpha$ with $\beta$) $\sgrad u_{\alpha}$ is written
as a linear combination of
$\sgrad u_{\beta}$ and $Y_{\beta,\nu}$ $(1\le\nu\le |\beta|)$. Thus the
claim follows by
induction on $\beta$.

Hence, it has been shown that each $Y_{\alpha,\nu}$
is written as a linear combination
of
$$Y_{\alpha,\mu}\quad (\mu\ne\nu),
\quad Y_{\gamma,\mu}\quad (\gamma\pcc\alpha,\
1\le\mu\le |\gamma|),$$
which are part of a basis of $\Gamma$.
Since $Y_{\alpha,\nu}$ is a primitive element, i.e.,
there is no integer $m>1$ such that $m^{-1}Y_{\alpha,\nu}\in\Gamma$,
the coefficients
are integers. This being true for every $\nu$,
we have
$$\frac{d_{\alpha,\nu}}{\prod_{0\le\mu\le |\alpha|\atop \mu\ne\nu}
(c_{\alpha,\mu}-c_{\alpha,\nu})}
=\pm \frac{d_{\alpha,\nu'}}{\prod_{0\le\mu\le |\alpha|\atop \mu\ne\nu'}
(c_{\alpha,\mu}-c_{\alpha,\nu'})}$$
for any $\nu$ and $\nu'$.

Note the sign of $d_{\alpha,\nu}$ is equal to
the sign of the function $v_{\alpha}
(c_{\alpha,\nu})$ on $M-L_{\alpha,\nu}$. This implies that the sign
of
$$\frac{d_{\alpha,\nu}}{\prod_{0\le\mu\le |\alpha|\atop \mu\ne\nu}
(c_{\alpha,\mu}-c_{\alpha,\nu})}$$
is equal to the sign of $u_{\alpha}$. In
particular it does not depend on $\nu$. This
completes the proof of the proposition.
\qed
\enddemo

We put
$$d_{\alpha}=\frac{d_{\alpha,\nu}}{\prod_{0\le\mu\le |\alpha|\atop \mu\ne\nu}
(c_{\alpha,\mu}-c_{\alpha,\nu})}\quad (0\le\nu\le |\alpha|).$$
Then, $d_{\alpha}^{-1} \sgrad u_{\alpha}=
\sum_{\nu=0}^{|\alpha|} Y_{\alpha,\nu}\in\Gamma$.
Put
$$Z_{\alpha}=\sum_{\nu=0}^{|\alpha|} Y_{\alpha,\nu}\in\Gamma$$
Note that $Z_{\alpha}=0$ if $\alpha$ is minimal. We
shall call $d_{\alpha}$ $(\alpha\in\A)$
the {\it scaling constants}.

For convenience we shall use two symbols $\p(\alpha)$
and $\n(\alpha)$: For each
non-minimal $\alpha\in \A$, $\p(\alpha)$ denotes the maximal element
of the totally
ordered subset $\{\gamma\in\A\ |\ \gamma\pcc\alpha\}$; for each non-maximal
$\alpha\in \A$, $\n(\alpha)$ denotes the subset of $\A$
defined by
$$\n(\alpha)=\{\gamma\in\A\ |\ \p(\gamma)=\alpha\}$$
($\n(\alpha)$ is identical with the one defined before).
For non-minimal elements $\alpha\in\A$ we define constants $m_{\alpha,\nu}$
$(0\le \nu
\le |\p(\alpha)|)$; putting $\beta=\p(\alpha)$,
$$
m_{\alpha,\nu} =\frac{d_{\beta}}{d_{\alpha}}
\prod_{0\le\mu\le |\beta|\atop \mu\ne\nu}
(c_{\beta,\mu}+e_{\beta\alpha}).\tag 4.8$$

\proclaim{Proposition 4.20}
$$\align
Z_{\alpha} & =m_{\alpha,0}Z_{\p(\alpha)}+\sum_{\nu=1}^{|\p(\alpha)|}
(m_{\alpha,\nu}-m_{\alpha,0})Y_{\p(\alpha),\nu}\\
& =\sum_{\beta; \text{non-minimal}\atop \beta\pc\alpha}
\left(\prod_{\beta\pcc\gamma\pc\alpha} m_{\gamma,0}\right)
\sum_{\nu=1}^{|\p(\beta)|} (m_{\beta,\nu}-m_{\beta,0})Y_{\p(\beta),\nu}.
\endalign$$
\endproclaim
\demo{Proof} Putting $\lambda=-e_{\beta\alpha}$ in the identity
$$\frac{v_{\beta}(\lambda)}{\prod_{\nu=0}^{|\beta|}(c_{\beta,\nu}-\lambda)}
=\sum_{\nu=0}^{|\beta|}\frac{1}{c_{\beta,\nu}-\lambda}
\frac{v_{\beta}(c_{\beta,\nu})}{\prod_{\mu\ne\nu}
(c_{\beta,\mu}-c_{\beta,\nu})},$$
and taking the skew gradient of both sides,
we have
$$\align
& Z_{\alpha}= \left(\frac{d_{\beta}}{d_{\alpha}}
\prod_{1\le\mu\le |\beta|}(c_{\beta,\mu}+e_{\beta\alpha})\right)Z_{\beta}\\
+ & \sum_{\nu=1}^{|\beta|}\left(\frac{d_{\beta}}{d_{\alpha}}
\prod_{0\le\mu\le |\beta|\atop \mu\ne\nu}
(c_{\beta,\mu}+e_{\beta\alpha})-\frac{d_{\beta}}{d_{\alpha}}
\prod_{1\le\mu\le |\beta|}(c_{\beta,\mu}+
e_{\beta\alpha})\right)Y_{\beta,\nu}.
\endalign$$
This proves the first equality. The second one
is immediate.
\qed
\enddemo

\proclaim{Proposition 4.21} The constants $m_{\alpha,\nu}$
$(\alpha$, not minimal, $0\le\nu\le\p(\alpha))$
possess the following properties.
\roster
\item"(1)" $m_{\alpha,\nu}-m_{\alpha,0}\in \Z$.
\item"(2)" $m_{\alpha,\nu}\in\Q$ if $\p(\alpha)$ is not minimal.
\item"(3)"
$$\left(\prod_{\beta\pcc\gamma\pc\alpha} m_{\gamma,0}\right)
(m_{\beta,\nu}-m_{\beta,0})\in\Z$$
for any non-minimal $\beta$ and $\alpha$ ($\succ\beta$).
\item"(4)" Either
$$0<m_{\alpha,0}<\dots<m_{\alpha,|\p(\alpha)|}$$
or
$$m_{\alpha,0}>\dots>m_{\alpha,|\p(\alpha)|}> 0.$$
\item"(5)" if $\alpha,\alpha'\in\n(\beta)$,
$\alpha\ne\alpha'$, then for any $\nu$,
$1\le\nu\le |\beta|$,
$$\frac{m_{\alpha,0}}{m_{\alpha,\nu}}\frac
{m_{\alpha,|\beta|}-m_{\alpha,\nu}}{m_{\alpha,|\beta|}-m_{\alpha,0}}
=\frac{m_{\alpha',0}}{m_{\alpha',\nu}}\frac
{m_{\alpha',|\beta|}-m_{\alpha',\nu}}{m_{\alpha',|\beta|}-m_{\alpha',0}}.$$
\endroster
\endproclaim
\demo{Proof} (1), (2), and (3) are immediately obtained
from the second equality in
Proposition 4.20. To prove $m_{\alpha,\nu}>0$, note that the
sign of $d_{\alpha}$ is equal to
the sign of $u_{\alpha}$. This implies that the
sign of $d_{\p(\alpha)}d_{\alpha}^{-1}$ is
equal to the sign of everywhere non-zero function
$\prod_{i\in\p(\alpha)}
(h_i+e_{\p(\alpha)\alpha})$. Thus we have
$m_{\alpha,\nu}>0$. The remaining inequalities
in (4) follows from the inequality
$$1=c_{\p(\alpha),0}>\dots>c_{\p(\alpha),|\p(\alpha)|}=0$$
and the fact that either $e_{\p(\alpha)\alpha}>0$
or $e_{\p(\alpha)\alpha}<-1$.
To prove (5), note that $d_{\beta}d_{\alpha}^{-1}$,
$c_{\beta,\nu}$ $(1\le\nu\le
|\beta|-1)$,
and $e_{\beta\alpha}$ are uniquely determined from
$m_{\alpha,\nu}$ $(0\le\nu\le
|\beta|)$,
where $\alpha\in\n(\beta)$. In particular we have
$$e_{\beta\alpha}=\frac{m_{\alpha,0}}{m_{\alpha,|\beta|}-m_{\alpha,0}},\qquad
c_{\beta,\nu}=\frac{m_{\alpha,0}}{m_{\alpha,\nu}}\frac
{m_{\alpha,|\beta|}-m_{\alpha,\nu}}
{m_{\alpha,|\beta|}-m_{\alpha,0}}.\tag 4.9$$
Hence (5) follows.
\qed
\enddemo

\remark{Remark} If $h_i$ $(i\in\beta)$ and
$e_{\beta\alpha}$ $(\beta\pcc\alpha)$ are
replaced with $1-h_{s(\beta)+t(\beta)-i}$
and $-1-e_{\beta\alpha}$ respectively for
a non-maximal $\beta$, then, (1) the order of
$Y_{\beta,0},\dots,Y_{\beta,|\beta|}$ and
$m_{\alpha,0},\dots,m_{\alpha,|\beta|}$ are reversed ($\alpha\in\n(\beta)$),
and (2) $d_{\gamma}$ is replaced with
$(-1)^{|\beta|}d_{\gamma}$ $(\beta\pcc\gamma)$.

If $\A$ is totally ordered, it is therefore
possible to choose $\{h_i\}$ and
$\{e_{\beta\alpha}\}$ so that every $e_{\beta\alpha}$ is positive and
$m_{\alpha,0}<\dots<m_{\alpha,|\p(\alpha)|}$ for every non-minimal $\alpha$.
But, in general it is impossible.
\endremark

\specialhead 5. Properties as a toric variety
\endspecialhead
In the previous section we have proved that
$M$ is a toric variety with respect
to the action of $G$. In this section
we shall specify the fan of the toric variety
$M$, and describe some properties that are useful
for the ``existence problem''.
Throughout this section we shall refer to Fulton
[3] and Oda [11] on the general
theory for toric varieties.

\head The fan of $M$
\endhead

As a toric variety, $M$ is constructed from
the lattice
$\Gamma\subset\frak k$ and a set $\D$ of polyhedral
cones in the Lie algebra
$\frak k$. The pair $(\Gamma,\D)$ (or the set
$\D$ if $\Gamma$ is known) is called the
fan of $M$. In our case, the invariant
hypersurfaces $L_{\alpha,\nu}$ $((\alpha,\nu)\in\J)$
and the information on their intersections will determine
$\D$. We first describe
$(\Gamma,\D)$ in terms of the partially ordered set
$\A$ and the numbers
$|\alpha|$, $m_{\alpha,\nu}$, and then prove that it is
the fan of $M$.

Let $\wt{\frak k}$ be the real vector space
of dimension $n+\#\A$ equipped with
the basis $\wt Y_{\alpha,\nu}$ $((\alpha,\nu)\in\J)$, and let $\wt\Gamma$
be the
lattice in $\wt{\frak k}$ generated by $\wt Y_{\alpha,\nu}$
$((\alpha,\nu)\in\J)$.
Define $\wt Z_{\alpha}\in\wt\Gamma$ $(\alpha\in\A)$ by
$$\wt Z_{\alpha}=\sum_{\beta; \text{non-minimal}\atop \beta\pc\alpha}
\left(\prod_{\beta\pcc\gamma\pc\alpha} m_{\gamma,0}\right)
\sum_{\mu=1}^{|\p(\beta)|} (m_{\beta,\mu}-m_{\beta,0})\wt Y_{\p(\beta),\mu}$$
if $\alpha$ is non-minimal, and by $\wt Z_{\alpha}=0$
if $\alpha$ is minimal.
Let $\Gamma_0$ be the subgroup of $\wt\Gamma$ generated
by
$$R_{\alpha}=\sum_{\nu=0}^{|\alpha|} \wt Y_{\alpha,\nu}- \wt Z_{\alpha}
\qquad (\alpha\in\A),$$
and let $\frak k_0$ be the subspace of
$\wt{\frak k}$ spanned by $R_{\alpha}$
$(\alpha\in\A)$. Then we have the exact sequences
$$\aligned
& 0\to\frak k_0\to\wt{\frak k}\overset\rho\to\to\frak k\to 0,\\
& 0\to\Gamma_0\to\wt\Gamma\overset\rho\to\to\Gamma\to 0,
\endaligned \tag 5.1$$
where $\rho$ is the homomorphism defined by $\rho(\wt
Y_{\alpha,\nu})=Y_{\alpha,\nu}$
$((\alpha,\nu)\in\J)$. Namely, $\Gamma$ and $\frak k$ are isomorphic
to
$\wt\Gamma/\Gamma_0$ and $\wt{\frak k}/{\frak k}_0$ respectively.

For each $\iota\in\I$, let $\sigma_{\iota}$ be the $n$-dimensional
cone in $\frak k$
generated by $Y_{\alpha,\nu}$ $((\alpha,\nu)\in\J(\iota))$, i.e.,
$$\sigma_{\iota}=\{\sum_{(\alpha,\nu)\in\J(\iota)} a_{\alpha,\nu}
Y_{\alpha,\nu}\ |\ a_{\alpha,\nu}\ge 0\}.$$
Let $\D$ be the set of the cones
$\sigma_{\iota}$ $(\iota\in\I)$ and all the faces of
them. Here, a face of the cone $\sigma_{\iota}$
means a cone $\sigma$ generated by
$Y_{\alpha,\nu}$ $((\alpha,\nu)\in\J_0)$ for some subset $\J_0\subset
\J(\iota)$. Hence the assignment $\sigma\to\J_0$ gives the one-to-one
correspondence between the cones in $\D$ and the
subsets of $\J$ contained in some
$\J(\iota)$. The $0$-cone $\{0\}$ is supposed to be
the face of every cone, which
corresponds to the empty subset of $\J$. It
is easily seen that $(\Gamma,\D)$
satisfies the conditions that a fan should satisfy,
and the resulting toric variety,
denoted by $X(\D)$, is compact and non-singular (see
[3] Sections 2.4 and 2.5).

Notice that the construction above only needs a
partially ordered set $\A$, integers
$|\alpha|$ $(\alpha\in\A)$, $m_{\alpha,\nu}-m_{\alpha,0}$
$(\alpha; \text{non-minimal},\ 0\le\nu\le |\p(\alpha)|)$,
and rational numbers
$m_{\alpha,0}$ $(\alpha,\,\p(\alpha); \text{non-minimal})$
satisfying Proposition 1.8 (2)
and Proposition 4.21. Namely, only the differences
$m_{\alpha,\nu}-m_{\alpha,0}$
are used
for such $\alpha$ that $\p(\alpha)$ is minimal, and
the condition (3) of Proposition 1.8 on
the integers $|\alpha|$ for maximal $\alpha$ are not
necessary for the construction.
Generally, fans and toric varieties obtained in such
a way from those data will be called
{\it of KL type}. If Proposition 1,8 (3)
is satisfied, then they
will be called {\it of KL-A type}. If
not, then they will be called {\it of KL-B type}.

\remark{Remark} Since the fan of a toric variety
is unique (cf. [12] Theorem 4.1), and
since $\D$ determines the sets of elements $\{Y_{\alpha,\nu}\}$
and $\{Z_{\alpha}\}$ of
$\Gamma$, it follows that the partially ordered set
$\A$ and the integers $|\alpha|$
$(\alpha\in\A)$ are uniquely associated with a toric variety
of KL type. Also, for each
$\alpha$ there are only two possibilities for the
ordering of $Y_{\alpha,0},\dots,
Y_{\alpha,|\alpha|}$ so that the corresponding
numbers $m_{\alpha,\nu}$ satisfy
Proposition 4.21 (4); the alternative is the reversed
order.
\endremark

\proclaim{Proposition 5.1} $(\Gamma,\D)$ is the fan of $M$.
\endproclaim

Proposition 5.1 will be proved by giving an
explicit identification of the toric variety
$X(\D)$ (of KL-A type) with $M$. So, let
us first review the construction of
$X(\D)$.

For each $\sigma\in\D$ we define the semigroup $\Cal
S_{\sigma}$ by
$$\Cal S_{\sigma}=\{\eta\in\Gamma^*\ |\ <\eta,Y>\ge 0\ \ \text{for any}\
Y\in\sigma
\cap\Gamma\},$$
where $\Gamma^*$ denotes the dual lattice of $\Gamma$.
Regarding $\C$ as the
multiplicative semigroup, we put
$$U_{\sigma}=\{u: \Cal S_{\sigma}\to
\C,\ \ \text{a semigroup homomorphism}\}.$$
Here, homomorphisms are assumed to map unit to
unit. If $\tau\subset\sigma$, then it is
easily seen that
$$\Cal S_{\sigma}\subset\Cal S_{\tau},\quad U_{\tau}\subset U_{\sigma}.$$
For the $0$-cone, $\Cal S_{\{0\}}=\Gamma^*$ and
$$U_{\{0\}}=\{u:\Gamma^*\to \C^{\times},\ \ \text{a group homomorphism}\}=
\Gamma\otimes_{\Z}\C^{\times},$$
which is an algebraic torus isomorphic to $(\C^{\times})^n$.
We shall denote it by
$\Cal T_{\Gamma}$. The group $\Cal T_{\Gamma}$
naturally acts on each $U_{\sigma}$ by
$$(u_0u)(\eta)=u_0(\eta)u(\eta),
\qquad u_0\in \Cal T_{\Gamma},\ u\in U_{\sigma},\ \eta\in
\Cal S_{\sigma}.$$
It follows from the definition that $U_{\sigma}$ is
an affine variety with coordinate
ring $\C[\Cal S_{\sigma}]$ (the semigroup ring). If the
cone $\sigma$ is $k$-dimensional,
then $U_{\sigma}$ is isomorphic to
$\C^k\times(\C^{\times})^{n-k}$. Then $X(\D)$
is
obtained by gluing all $U_{\sigma}$ with the relations
$$U_{\sigma}\supset U_{\sigma\cap\tau}\subset U_{\tau}.$$
The action of $\Cal T_{\Gamma}$ on $X(\D)$ is
also well-defined.

\demo{Proof of Proposition 5.1} Fix a point $p_0\in
M^1$. Let $\iota\in\I$, and let
$(z_{\alpha,\nu}^{(\iota)};\ (\alpha,\nu)\in\J(\iota))$
be the coordinate system on
$M-\cup_{\alpha\in\A} L_{\alpha,\iota(\alpha)}$ given
in
Corollary 4.16. Also, let
$\eta_{\alpha,\nu}^{(\iota)}$ $((\alpha,\nu)\in\J(\iota))$
be the basis of $\Gamma^*$ dual to $Y_{\alpha,\nu}$
$((\alpha,\nu)\in\J(\iota))$ (note
that they are generators for $\Cal S_{\sigma_{\iota}}$). Then,
there is a natural
holomorphic isomorphism
$$U_{\sigma_{\iota}}\to M-\cup_{\alpha\in\A} L_{\alpha,\iota(\alpha)}\quad
(u\mapsto p)\tag 5.2$$
given by
$$u(\eta_{\alpha,\nu}^{(\iota)})=
z_{\alpha,\nu}^{(\iota)}(p)\quad ((\alpha,\nu)
\in\J(\iota)).\tag 5.3$$
If $\sigma\in\D$ is a face of $\sigma_{\iota}$, then
the mapping (5.2) gives the
holomorphic isomorphism
$$U_{\sigma}\to M-\cup_{(\alpha,\nu)\in\J-\J_0} L_{\alpha,\nu},\tag 5.4$$
where $\J_0$ is the subset of $\J$ corresponding
to $\sigma$;
$$\J_0=\{(\alpha,\nu)\in\J\ |\ Y_{\alpha,\nu}\in\sigma\}\subset \J(\iota).$$
It is easily seen that the isomorphism (5.4)
is independent of the choice of
$\sigma_{\iota}$ containing $\sigma$. Hence we obtain the holomorphic
isomorphism
$X(\D) \to M$.

Defining the isomorphism $\Cal T_{\Gamma}\to G$ by
$$Y_{\alpha,\nu}\otimes e^{t+\sqrt{-1}s}
\mapsto \exp(-tIY_{\alpha,\nu}+sY_{\alpha,\nu}),$$
we can easily see that the isomorphism $X(\D)
\to M$ commutes with the actions of the
groups $\Cal T_{\Gamma}$ and $G$. This completes the
proof.
\qed
\enddemo

{}From now on, we shall fix a point
$p_0\in M^1$ (the base point) and identify each
$U_{\sigma}$ with the subset
$M-\cup_{(\alpha,\nu)\in\J-\J_0} L_{\alpha,\nu}$ of $M$
by the isomorphism given in the proof of
Proposition 5.1, where $\J_0\subset \J$
corresponds to $\sigma$. Also, the group $\Cal T_{\Gamma}$
will be identified with $G$.
For $\sigma\in\D$ corresponding to $\J_0$ we put
$$O_{\sigma}=U_{\sigma}\cap \cap_{(\alpha,\nu)\in\J_0} L_{\alpha,\nu}.$$
Then, $O_{\sigma}$ is a $G$-orbit isomorphic to $(\C^{\times})^{n-k}$
$(k=\dim\sigma)$, and
$$M=\cup_{\sigma\in\D}, O_{\sigma}\quad
U_{\sigma}=\cup_{\tau\subset\sigma} O_{\tau}$$
Note that $O_{\{0\}}=M^1$. Let $V(\sigma)$ be the closure
of
$O_{\sigma}$ in $M$. Then we have
$$V(\sigma)=\cap_{(\alpha,\nu)\in\J_0} L_{\alpha,\nu}=
\cup_{\tau\supset\sigma} O_{\tau}$$
(see [3] Section 3.1).

Since $\D$ contains $n$-dimensional cones, we have:

\proclaim{Corollary 5.2} $M$ is simply connected.
\endproclaim

For the proof, see [3] p.56, Proposition. Actually,
one can see more about the
topology of $M$: There is a cell-decomposition of
$M$ consisting of $\prod_{\alpha}
(|\alpha|+1)$ even-dimensional cells; the number of $2k$-dimensional cells
is equal to the
number of $\iota\in\I$ such that $\sum_{\alpha} \iota(\alpha)=k$. This
is made in a
similar way as [3] pp.101-103. Since the result
is not used in this paper, we
omit the detail.

\head Fibre bundles associated with $M$
\endhead
In general, let $(\Gamma_i,\D_i)$ $(i=1,2)$ be two fans,
and let $\phi:\Gamma_1\to
\Gamma_2$ be a homomorphism such that the induced
linear homomorphism
$$\phi:\Gamma_1\otimes_{\Z}\R\to\Gamma_2\otimes_{\Z}\R$$
maps each $\sigma\in\D_1$ into some $\sigma'\in\D_2$. Denoting the
resulting toric varieties by $X(\D_i)$, we have

\proclaim{Proposition 5.3} There is a natural holomorphic mapping
$\phi_{\sharp}:X(\D_1)
\to X(\D_2)$ that is equivariant with respect to
the naturally induced homomorphism
$$\phi_{\sharp}: \Cal T_{\Gamma_1}\to\Cal T_{\Gamma_2}$$
of algebraic tori.
\endproclaim

For the proof, see [11] p.19, Theorem 1.13
(see also [3] p.41, Exercise). As
applications of this general result, we shall obtain
two kinds of fibre bundles associated
with $M$. We now explain the first one.
Let $\wt\Gamma$ and $\wt{\frak k}$ be as in
the previous subsection. For each $\sigma\in\D$ corresponding to
$J_0\subset\J$, define
the cone $\wt\sigma$ in $\wt{\frak k}$ by
$$\wt\sigma=\{\sum_{(\alpha,\nu)\in\J_0} a_{\alpha,\nu}
\wt Y_{\alpha,\nu}\ |\ a_{\alpha,\nu}\ge 0\},$$
and put $\wt\D=\{\wt\sigma\ |\ \sigma\in\D\}$. As is easily
verified, $(\wt\Gamma,
\wt\D)$ is a fan. Then, by Proposition 5.3
the homomorphism $\rho$ induces the
equivariant holomorphic mapping $\rho_{\sharp}: X(\wt\D)\to M$.

\proclaim{Proposition 5.4} \roster
\runinitem"(1)" The toric variety $X(\wt\D)$ and the algebraic
torus
$\Cal T_{\wt\Gamma}$ are naturally isomorphic to
$$\prod_{\alpha\in\A}(\C^{|\alpha|+1}-\{0\})=
\{(z_{\alpha};\ \alpha\in\A)\ |\ z_{\alpha}
=(z_{\alpha,0},\dots,z_{\alpha,|\alpha|})\in\C^{|\alpha|+1}-\{0\}\}$$
and
$$(\C^{\times})^{n+\#\A}=
\{(\lambda_{\alpha,\nu};\ (\alpha,\nu)\in\J)\ |\ \lambda_
{\alpha,\nu}\in\C^{\times}\}$$
respectively.
\item"(2)" $\rho_{\sharp}: X(\wt\D)\to M$ is a principal fibre
bundle with structure group
$\Cal T_{\Gamma_0}$.
\endroster
\endproclaim
\demo{Proof} (1) Let $\wt Y_{\alpha,\nu}^*$ $((\alpha,\nu)\in\J)$ be the
basis of $\wt
\Gamma^*$ dual to $\wt Y_{\alpha,\nu}$ $((\alpha,\nu)\in\J)$. Then, all
the semigroups
$\Cal S_{\wt\sigma}$ contain the semigroup generated by
$\wt Y_{\alpha,\nu}^*$ $((\alpha,\nu)\in\J)$. This implies that all the
affine varieties
$U_{\wt\sigma}$ are realized in $\C^{n+\#\A}$;
$$\gather
\C^{n+\#\A}=\{(z_{\alpha,\nu};\ (\alpha,\nu)\in\J)\},\\
U_{\wt\sigma}= \{(z_{\alpha,\nu})\ |\ z_{\alpha,\nu}
\ne 0\text{ for } (\alpha,\nu)
\not\in\J_0\},
\endgather$$
$\J_0$ corresponding to $\sigma$. From this the assertion
easily follows.

(2) We first review how the mapping $\rho_{\sharp}$
is constructed: The surjective
homomorphism $\rho:\wt\Gamma\to\Gamma$ induces the inclusion $\rho^*:\Gamma^*
\to \wt\Gamma^*$. This gives the inclusion $\Cal S_{\sigma}\to
\Cal S_{\wt\sigma}$
for any $\sigma\in\D$. Thus the mapping
$\rho_{\sharp}:U_{\wt\sigma}\to U_{\sigma}$
is defined by
$$u\mapsto u|_{\Cal S_{\sigma}}\qquad (u\in U_{\wt\sigma}).$$

Now, fix $\iota\in\I$ and define a splitting $j_{\iota}:\Gamma\to\wt\Gamma$
of the
exact sequence (5.1) by $j_{\iota}(Y_{\alpha,\nu})=\wt Y_{\alpha,\nu}$
$((\alpha,\nu)\in\J(\iota))$. Let $\Gamma_{\iota}$ be its image. Then we
have
the direct sum decompositions
$$\wt\Gamma=\Gamma_0 +\Gamma_{\iota},\qquad
\wt\Gamma^*=\Gamma_{\iota}^{\perp} + \Gamma^*,$$
and accordingly,
$$\Cal S_{\wt\sigma_{\iota}}=
\Gamma_{\iota}^{\perp}+\Cal S_{\sigma_{\iota}}.$$
Since $\Gamma_{\iota}^{\perp}$ is identified with $\Gamma_0^*$ by the
natural
homomorphism $\wt\Gamma^*\to\Gamma_0^*$, we thus obtain the holomorphic
isomorphism
$$U_{\wt\sigma_{\iota}}\to \Cal T_{\Gamma_0}\times U_{\sigma_{\iota}}\qquad
(u\mapsto (u|_{\Gamma_{\iota}^{\perp}},
u|_{\Cal S_{\sigma_{\iota}}})).\tag 5.5$$
Clearly, the mapping above also gives the isomorphism
$\Cal T_{\wt\Gamma}\to
\Cal T_{\Gamma_0}\times G$ of
algebraic tori, with which the isomorphism (5.5) is
equivariant. This proves (2).
\qed
\enddemo

Now, let us explain the other kind of
fibre bundles that are naturally associated with $M$.
For convenience we shall introduce a topology on
the set $\A$: A subset $\Cal B$ of $\A$ is
open if it possesses the property;
$$\text{if\ }\beta\in\Cal B\ \text{and } \gamma\pc\beta,\ \text{then }
\gamma\in\Cal B.\tag 5.6$$
Let $\A'$ be an open subset of $\A$,
and put $\A''=\A-\A'$. Let $\Gamma'$ be a subgroup
of $\Gamma$ generated by $Y_{\alpha,\nu}$
$((\alpha,\nu)\in\J,\ \alpha\in\A')$, and
let
$\frak k'$ be the subspace of $\frak k$
spanned by those vectors. Put
$$\Gamma''=\Gamma/\Gamma',\quad \frak k''=\frak k/\frak k',$$
and let $\pi:\frak k\to \frak k''$ be the
natural projection. Also, put
$$\D' =\{\sigma\in\D\ |\ \sigma\subset \frak k'\},\quad
\D'' =\{\pi(\sigma)\ |\ \sigma\in\D\}.$$
Then the pairs $(\Gamma',\D')$ and $(\Gamma'',\D'')$ become fans.
It is clear that
$X(\D'')$ is a toric variety of KL-A type.
For $X(\D')$, it can be only said that
it is of KL type. Since the homomorphism
$\pi$ satisfies the assumption of
Proposition 5.3, we have the equivariant mapping
$$\pi_{\sharp}: M\to X(\D'').$$

\proclaim{Proposition 5.5} $\pi_{\sharp}: M\to X(\D'')$ is a fibre
bundle with typical
fibre $X(\D')$. More precisely, for each $\iota\in\I$ there
is an isomorphism $G\to
\Cal T_{\Gamma'}\times \Cal T_{\Gamma''}$ of algebraic tori and
an equivariant
holomorphic isomorphism
$$\pi_{\sharp}^{-1}(U_{\pi(\sigma_{\iota})})
\to X(\D')\times U_{\pi(\sigma_{\iota})}.$$
\endproclaim

\demo{Proof} Fix $\iota\in\I$, and let $\Gamma_1$ be the
subgroup of $\Gamma$
generated by $Y_{\alpha,\nu}$
$((\alpha,\nu)\in\J(\iota),\ \alpha\in\A'')$. We then have
the direct sum decompositions
$$\Gamma=\Gamma'+\Gamma_1,\quad \Gamma^*=\Gamma_1^\perp +(\Gamma'')^*,$$
and
$$\Cal S_{\sigma_{\iota}}=
\Cal S_{\sigma'_{\iota}}+ \Cal S_{\pi(\sigma_{\iota})},$$
where $\sigma'_{\iota}=
\sigma_{\iota}\cap \frak k'\in\D'$, and $\Gamma_1^\perp$ is
identified with $(\Gamma')^*$ by the projection $\Gamma^*\to (\Gamma')^*$.
This induces
the holomorphic isomorphism
$$U_{\sigma_{\iota}}\to U_{\sigma'_{\iota}}
\times U_{\pi(\sigma_{\iota})},\tag 5.7$$
and the isomorphism
$$\Cal T_{\Gamma}\to\Cal T_{\Gamma'}\times\Cal T_{\Gamma''}$$
of algebraic tori so that the mapping (5.7)
is equivariant. Then, varying $\iota\in\I$
so that $\iota(\alpha)$ remains unchanged for any $\alpha\in\A''$,
and taking the union
of both sides of (5.7) with respect to
all such $\iota$, we have the equivariant isomorphism
$$\pi_{\sharp}:\pi^{-1}(U_{\pi(\sigma_{\iota})})
\to X(\D')\times U_{\pi(\sigma_{\iota})}.$$
\qed
\enddemo
Let $\A=\cup_{i=1}^k \A_i$ (disjoint union) be the decomposition
of $\A$ into the
connected components. Let $\Gamma_i$ be the subgroup of
$\Gamma$ generated by
$Y_{\alpha,\nu}$ $(\alpha\in\A_i)$, and let $\frak k_i$ be the
subspace of $\frak k$
spanned by $\Gamma_i$. Clearly, $\Gamma$ is the direct
sum of those subgroups, and
correspondingly,
$$G=\Cal T_{\Gamma_1}\times\dots\times \Cal T_{\Gamma_k}.$$
Putting
$$\D_i=\{\sigma\cap \frak k_i\ |\ \sigma\in\D\},$$
we obtain fans $(\Gamma_i,\D_i)$.

\proclaim{Corollary 5.6} There is a natural holomorphic isomorphism
$$M\to X(\D_1)\times\dots\times X(\D_k)$$
that is equivariant with respect to the identification
$$G=\Cal T_{\Gamma_1}\times\dots\times \Cal T_{\Gamma_k}.$$
Moreover, each $X(\D_i)$ naturally becomes a K\"ahler-Liouville manifold
so that
$M$ becomes the product manifold as K\"ahler-Liouville manifold.
\endproclaim
\demo{Proof} The former half is an immediate consequence
of Proposition 5.5.
The latter half is then obvious.
\qed
\enddemo

We now go back to the situation of
Proposition 5.5 and observe the fibre bundle
$\pi_{\sharp}: M\to X(\D'')$ from another point of view.
Let $\wt\Gamma'$ and
$\wt\Gamma''$ be the subgroups of $\wt\Gamma$ generated by
$\wt Y_{\alpha,\nu}$
$((\alpha,\nu)\in\J,\ \alpha\in\A')$ and
$\wt Y_{\alpha,\nu}$ $((\alpha,\nu)\in
\J,\ \alpha\in\A'')$ respectively. We then have
$$\wt\Gamma=\wt\Gamma'+ \wt\Gamma''\quad
(\text{direct sum}).\tag 5.8$$
Let $\wt\pi: \wt\Gamma\to\wt\Gamma''$ be the projection.
The homomorphism
$\rho:\wt\Gamma\to\Gamma$ induces the homomorphisms
$\wt\Gamma'\to\Gamma'$ and
$\rho'':\wt\Gamma''\to\Gamma''$ (the latter is given by
$\pi\circ\rho$). Let $\Gamma'_0$ and $\Gamma''_0$ be the kernel
of those
homomorphisms. Thus we have the following commutative diagram
whose rows
and columns are exact:
$$\CD
@. 0 @. 0 @. 0 @.\\
@. @VVV @VVV @VVV @.\\
0 @>>> \Gamma'_0 @>>> \wt\Gamma' @> \rho >>
\Gamma' @>>> 0\\
@. @VVV @VVV @VVV @.\\
0 @>>> \Gamma_0 @>>> \wt\Gamma @> \rho >>
\Gamma @>>> 0\\
@. @VV \wt\pi V @VV \wt\pi V @VV
\pi V @.\\
0 @>>> \Gamma''_0 @>>> \wt\Gamma'' @> \rho''>> \Gamma''
@>>> 0\\
@. @VVV @VVV @VVV @.\\
@. 0 @. 0 @. 0 @.
\endCD $$

Note that the splitting in the mid column
does not induce the splitting in the left one,
i.e., $\Gamma_0\cap\wt\Gamma''\ne
\Gamma''_0$ in general. Let $\wt{\frak k}''$
be
the subspace of $\wt{\frak k}$ spanned by $\wt\Gamma''$,
and put
$$\wt\D''=\{\sigma\in\wt\D\ |\ \sigma\subset\wt{\frak k}''\}.$$
Then, $(\wt\Gamma'',\wt\D'')$ becomes a fan, and the homomorphism
$\rho''$ induces
the principal $\Cal T_{\Gamma''_0}$-bundle
$\rho''_{\sharp}:X(\wt\D'')\to X(\D'')$.

We now define the homomorphism $\psi:\Gamma''_0\to \Gamma'$ as
follows:
Let $\wt\Gamma\to \wt\Gamma'$ be the projection with respect
to the decomposition
(5.8). Restricting it to $\Gamma_0$, we have the
homomorphism $\Gamma_0
\to\wt\Gamma'$. The restriction of this mapping to $\Gamma'_0$
being the identity,
we thus obtain the homomorphism
$$\psi_1:\Gamma''_0=\Gamma_0/\Gamma'_0\to\wt\Gamma'/\Gamma'_0=\Gamma'.$$
We put $\psi=-\psi_1$. The following formula is easily
obtained:
$$\psi(\wt\pi(R_{\alpha}))=
\left(\prod_{\alpha_0\pcc\beta\pc\alpha}m_{\beta,0}\right)
Z_{\alpha_0}\qquad (\alpha\in\A''),$$
where $\alpha_0$ is the minimal element of $\A''$
satisfying $\alpha_0\pc\alpha$.
The induced homomorphism $\Cal T_{\Gamma''_0}\to
\Cal T_{\Gamma'}$ of algebraic
tori
is also denoted by $\psi$. Through this homomorphism
$\Cal T_{\Gamma''_0}$ acts on $X(\D')$.

\proclaim{Proposition 5.7} The fibre bundle $\pi_{\sharp}:M\to X(\D'')$ is
isomorphic
to the fibre product
$$X(\wt\D'')\times_{\Cal T_{\Gamma''_0}} X(\D')\to X(\D'').$$
\endproclaim

\demo{Proof} Let $\wt{\frak k}'$ be the subspace of
$\wt{\frak k}$ spanned by
$\wt\Gamma'$, and put $\wt\D'=
\{\sigma\in\wt\D\ |\ \sigma\subset\wt{\frak k}'\}$.
Then $(\wt\Gamma',\wt\D')$ is also a fan, and we
have
$$X(\wt\D)=X(\wt\D'')\times X(\wt\D').$$
Hence $M=X(\wt\D)/\Cal T_{\Gamma_0}$ is equal to
$$((X(\wt\D'')\times X(\wt\D'))/\Cal T_{\Gamma'_0})/\Cal T_{\Gamma''_0}
=(X(\wt\D'')\times X(\D'))/\Cal T_{\Gamma''_0}.$$
Since the action of $\Cal T_{\Gamma''_0}$ on $X(\wt\D'')\times
X(\D')$ is given by
$$(p,q)g=(pg, \psi_1(g)q)=(pg,\psi(g^{-1})q),$$
the proposition follows.
\qed
\enddemo

\head Line bundles
\endhead

Let Pic($M$) denote the group of the isomorphism
classes of holomorphic line bundles
over $M$. To each $\xi\in\wt\Gamma^*$ we associate a
divisor of $M$;
$$\xi\mapsto \sum_{(\alpha,\nu)\in\J}
<\xi,\wt Y_{\alpha,\nu}>L_{\alpha,\nu},$$
where $<,>$ denotes the natural pairing of $\wt\Gamma^*$
and $\wt\Gamma$.
Let $Q_{\xi}$ be the line bundle over $M$
associated with this divisor. Then we have the
homomorphism $\wt\Gamma^*\to \text{Pic}(M)$ $(\xi\mapsto Q_{\xi})$. Also, the
homomorphism $-\xi:\wt\Gamma\to \Z$ induces the homomorphism $\chi_{-\xi}:
\Cal T_{\wt\Gamma}\to\C^{\times}$. Restricting it to $\Cal T_{\Gamma_0}$, we
have
another line bundle over $M$ associated with $\pi_{\sharp}:X(\wt\D)\to
M$.

\proclaim{Proposition 5.8}
\roster
\runinitem"(1)" The following sequence is exact:
$$0\to \Gamma^*\overset{\rho^*}\to\to\wt\Gamma^*\to\text{\rm Pic}(M)\to 0.$$
\item"(2)" $Q_{\xi}$ is isomorphic to the fibre product
$X(\wt\D)\times_{\chi_{-\xi}}\C$,
where $\chi_{-\xi}$ is regarded as the homomorphism $\Cal
T_{\Gamma_0}\to\C^{\times}$
by restriction.
\item"(3)" The assignment of the first Chern class
$c_1(Q)$ to each $Q\in
\text{\rm Pic}(M)$ gives the isomorphism $\text{\rm Pic}(M)\to H^2(M,\Z)$.
\endroster
\endproclaim
\demo{Proof} For (1) and (3), see [3] pp.63-64
and [11] Corollary 2.5. (2) is easy.
\qed
\enddemo

Put
$$\zeta_{\alpha}=c_1(Q_{\wt Y_{\alpha,0}^*})
\in H^2(M,\Z)\qquad (\alpha\in\A).$$
The proposition above implies that the elements $\zeta_{\alpha}$
$(\alpha\in\A)$ form
a basis of $H^2(M,\Z)$. Its dual basis is
given as follows: Let $\tau(\alpha)\in\D$ be the
$(n-1)$-dimensional cone generated by
$$\{Y_{\beta,\nu}\ |\ 1\le\nu\le |\beta|
\text{ if }\beta\ne\alpha ;\ 2\le\nu\le
|\beta|
\text{ if }\beta=\alpha\}.$$
Then $V(\tau(\alpha))$ is 1-dimensional (isomorphic to $\C P^1$).
Let $[V(\tau(\alpha))]$
denote its fundamental class in $H_2(M,\Z)$.

\proclaim{Proposition 5.9} $<\zeta_{\alpha},\ [V(\tau(\beta))]>
=\delta_{\alpha\beta}$,
where $<,>$ denotes the natural pairing of $H^2(M,\Z)$
and $H_2(M,\Z)$.
\endproclaim
\demo{Proof} The assertion easily follows from the fact:
$$L_{\alpha,0}\cap V(\tau(\beta))=\cases
0\quad & \text{if } \alpha\ne\beta\\
V(\sigma_{\iota})=\{\text{a point}\}\quad & \text{if } \alpha=\beta,
\endcases$$
where $\iota\in\I$ is given by $\iota(\gamma)=0$ for every
$\gamma\in\A$.
\qed
\enddemo

The next theorem specifies the cohomology class $[\omega]\in
H^2(M,\R)$ of the K\"ahler
form $\omega$. Let $\A_i$ $(1\le i\le k)$ be
the connected components of $\A$, and let
$\alpha_i$ denote the unique minimal element of $\A_i$.

\proclaim{Theorem 5.10}
\roster
\item"(1)"
$$\int_{V(\tau(\alpha))} \omega =
\frac{2\pi}{d_{\alpha}}\prod_{\beta\pcc\alpha}
\prod_{\nu=1}^{|\beta|} (c_{\beta,\nu}+e_{\beta\alpha}).$$
\item"(2)"
$$[\omega]= \sum_{i=1}^k \frac{2\pi}
{d_{\alpha_i}}\left(\sum_{\alpha\succ\alpha_i}
\left(\prod_{\alpha_i\pcc\beta\pc\alpha} m_{\beta,0}\right)\, \zeta_{\alpha}
+ \zeta_{\alpha_i}\right).$$
\endroster
\endproclaim
\demo{Proof} (1) As is easily seen, the action
of the circle group $\{\exp(tY_{\alpha,0})\}$
$(t\in\R/2\pi\Z)$ on $V(\tau(\alpha))$ has the two fixed points
$q_0$ and $q_1$;
$$\{q_0\}=V(\tau(\alpha))\cap V(\R_{\ge 0}Y_{\alpha,0}),\qquad
\{q_1\}=V(\tau(\alpha))\cap V(\R_{\ge 0}Y_{\alpha,1}),$$
where $\R_{\ge 0}$ denotes the set of non-negative
real numbers.
Let $\gamma(s)$ $(0\le s\le l)$ be a geodesic
of unit speed on $V(\tau(\alpha))$ from
$q_0$ to $q_1$. Then, parametrizing $V(\tau(\alpha))$ by $(s,t)$,
we have
$$\int_{V(\tau(\alpha))} \omega=
\int_0^l\int_0^{2\pi} \omega(\frac{\partial}{\partial s},\,
\frac{\partial}{\partial t})dtds,$$
and
$$\omega(\frac{\partial}{\partial s},\,\frac{\partial}{\partial t})
= d_{\alpha,0}^{-1} \frac{d}{ds} v_{\alpha}(c_{\alpha,0})(\gamma(s)).$$
Hence
$$\int_{V(\tau(\alpha))} \omega=
\frac{2\pi}{d_{\alpha,0}} \left(v_{\alpha}(c_{\alpha,0})
(q_1) - v_{\alpha}(c_{\alpha,0})(q_0)\right).$$
Since $v_{\alpha}(c_{\alpha,0})(q_0)=0$ and
$$v_{\alpha}(c_{\alpha,0})(q_1)=
\left(\prod_{\beta\pcc\alpha}\prod_{\nu=1}^{|\beta|}
(c_{\beta,\nu}+e_{\beta\alpha})\right)\,
\prod_{\mu=1}^{|\alpha|} (c_{\alpha,\mu}-c_{\alpha,0}),$$
the assertion follows. (2) is immediately obtained from
(1).
\qed
\enddemo

\specialhead 6. Bundle structure associated with a subset
of $\A$
\endspecialhead
In the previous section we have proved that
an open subset $\A'$ induces the
fibre bundle $\pi_{\sharp}:M\to X(\D'')$ whose typical fibre is
$X(\D')$.
In this section we shall show that the
toric varieties $X(\D')$ and $X(\D'')$
naturally possess structures of K\"ahler-Liouville manifold inherited
from $M$. Since the numbering $i=1,\dots,n$ is inconvenient
for the purpose of
this section, we shall use $(\alpha,\nu)$ $(\alpha\in\A,1\le \nu\le|\alpha|)$
instead.
The correspondence is given by
$$(\alpha,\nu)\leftrightarrow i=s(\alpha)+\nu-1.$$

Fix an open subset $\A'$, and let $(\Gamma',\D')$
and $(\Gamma'',\D'')$
be as in the previous section. $TM$ is
naturally decomposed to the sum of
(mutually orthogonal) two subbundles; $TM=D'+D''$, where
$$D'=\sum_{\alpha\in\A'}D_{\alpha},
\quad D''=\sum_{\alpha\in\A''}D_{\alpha}.$$
Clearly, $D'$ is integrable, and the maximal integral
submanifolds are the fibres
of $\pi_{\sharp}:M\to X(\D'')$. Let $\A''_1,\dots,\A''_r$ be the connected
components of $\A''$, and let $\alpha_s$ be the
(unique) minimal element of
$A''_s$ $(1\le s\le r)$. Put
$$D''_s=\sum_{\alpha\in\A''_s}D_{\alpha}\qquad (1\le s\le r).$$
Recalling the orthonormal frame $V_{\alpha,\nu},IV_{\alpha,\nu}$
$(\alpha\in\A,1\le \nu\le|\alpha|)$ over $M^1$, we put
$$V''_{\alpha,\nu}=\sqrt{|u_{\alpha_s}|}
V_{\alpha,\nu}\qquad (\alpha\in\A''_s ).$$
The following lemma is easily obtained by using
the properties of the vector
fields $W_{\alpha,\nu}$ (cf. Proposition 1.4).

\proclaim{Lemma 6.1} For $\alpha\in\A''$ and $\beta\in\A'$,
$$[V''_{\alpha,\nu},V_{\beta,\mu}]=[V''_{\alpha,\nu},IV_{\beta,\mu}]=
[IV''_{\alpha,\nu},V_{\beta,\mu}]=[IV''_{\alpha,\nu},IV_{\beta,\mu}]=0.$$
\endproclaim

Let us recall the polynomial $F_{\alpha}(\lambda)$ in the
indeterminate $\lambda$ whose
coefficients are elements of $\Cal F$ (cf. Section
2):
$$\align
& F_{\alpha}(\lambda) =|u_{\alpha}|\sum_{\nu=1}^{|\alpha|}
\prod_{1\le\mu\le |\alpha|\atop \mu\ne\nu}(h_{\alpha,\mu}-\lambda)\cdot
(V_{\alpha,\nu}^2+(IV_{\alpha,\nu})^2)\tag 6.1\\
& +|u_{\alpha}|\sum_{\beta\in \n(\alpha)}
\frac{\prod_{\nu=1}^{|\alpha|}(h_{\alpha,\nu}+
e_{\alpha\beta})-\prod_{\nu=1}^{|\alpha|}
(h_{\alpha,\nu}-\lambda)}{e_{\alpha\beta}
+\lambda}\sum_{\gamma\succeq\beta}\sum_{\mu=1}^{|\gamma|}
(V_{\gamma,\mu}^2+(IV_{\gamma,\mu})^2).
\endalign$$
This polynomial is uniquely determined if the fundamental
functions and the conjunction
constants are specified. We shall call it the
{\it generating polynomial}. The next
proposition is an immediate consequence of the lemma
above.

\proclaim{Proposition 6.2} \roster
\runinitem"(1)" The vector fields $V''_{\alpha,\nu}$
are $\Cal T_{\Gamma'}$-invariant.
\item"(2)" The horizontal subbundle $D''$ is $\Cal T_{\Gamma'}$-invariant.
\item"(3)" For any $\alpha\in\A''$, the coefficients of $F_{\alpha}(\lambda)$
are
$\Cal T_{\Gamma'}$-invariant, and are sections of $S^2D''$.
\endroster
\endproclaim

By virtue of Proposition 6.2 (3), the coefficients
of $(\pi_{\sharp})_*F_{\alpha}(\lambda)$
$(\alpha\in\A'')$ are well-defined sections of $S^2TX(\D'')$. Let $\Cal
F''$ be the vector
space spanned by those sections.
Also, the riemannian metric $g''$ on $X(\D'')$ is
defined by the conditions: The
subbundles
$D''_s$ $(1\le s\le r)$ are mutually orthogonal with
respect to
$\pi_{\sharp}^*g''$, and
$$\pi_{\sharp}^*g''=|u_{\alpha_s}|^{-1}g
\quad \text{on } D''_s\qquad (1\le s\le r).$$
It is easily seen that $g''$ is a
K\"ahler metric. We denote by $M''$ the K\"ahler
manifold $(X(\D''),g'')$.

\proclaim{Theorem 6.3} $(M'',\Cal F'')$ is a K\"ahler-Liouville manifold
of
type (A). It possesses the following properties:
\roster
\item"(1)" The associated partially ordered set is naturally
identified with $\A''$;
\item"(2)" the underlying toric variety is identical with
$X(\D'')$;
\item"(3)" the fundamental functions
$\{h''_{\alpha,\nu}\}$ $(\alpha\in\A'',\ 1\le
\nu\le |\alpha|)$ are given by
$\pi_{\sharp}^*(h''_{\alpha,\nu})=h_{\alpha,\nu}$;
\item"(4)" the conjunction constants
$e''_{\alpha\beta}$ $(\alpha\in\A'',\ \alpha
\pc\beta)$ are given by $e''_{\alpha\beta}=e_{\alpha\beta}$;
\item"(5)" the scaling constants $d''_{\alpha}$ $(\alpha\in\A'')$ are given
by
$d''_{\alpha}=\epsilon(\alpha_s)d_{\alpha}$ $(\alpha\in\A''_s)$,
where $\epsilon(\alpha_s)$ is the sign of $d_{\alpha_s}$.
\endroster
\endproclaim
\demo{Proof} The commutativity of $\Cal F''$ with respect
to the Poisson bracket
follows from that of $\Cal F$. Since maximal
elements of $\A''$ are also maximal
in $\A$, it follows that $|\alpha|\ge 2$ for
any maximal element $\alpha$ of $\A''$.
This implies that $(M'',\Cal F'')$ is of type
(A). The properties $(1),\dots,(5)$ are
easily verified.
\qed
\enddemo

Next, let us consider the fibre. Define the
K\"ahler metric $g'(q)$ on the fibre
$\pi_{\sharp}^{-1}(q)$ $q\in M''$ by restricting $g$. With this
metric we regard
$\pi_{\sharp}^{-1}(q)$
as a K\"ahler manifold. Also, we define $\Cal
F'(q)$ as follows: Each
$F\in\Cal F$ is a section of $S^2D'+S^2D''$; so,
taking the $S^2D'$-component $F'$
of $F$, we put
$$\Cal F'(q)= \{F'|_{\pi_{\sharp}^{-1}(q)}\ |\ F\in\Cal F\}.$$

\proclaim{Theorem 6.4}
\roster
\runinitem"(1)" $(\pi_{\sharp}^{-1}(q),
\Cal F'(q))$ is a K\"ahler-Liouville manifold
for any $q\in M''$.
\item"(2)" Let $\wt X$ be the horizontal lift
(i.e., the lift as a section of $D''$) of a
vector field $X$ on $M''$. Then the one-parameter
group $\{\phi_t\}$ of transformations
of $M$ generated by $\wt X$ gives the
automorphisms
$\pi_{\sharp}^{-1}(q)\to \pi_{\sharp}^{-1}
(\phi_t(q))$ of K\"ahler manifolds, and
preserves $F'$ for each $F\in\Cal F$.
\endroster
\endproclaim
\demo{Proof} (1) Let $F'_{\alpha}(\lambda)$ $(\alpha\in\A')$ be the
$S^2D'$-component of $F_{\alpha}(\lambda)$. We have
$$F_{\alpha}(\lambda)= F'_{\alpha}(\lambda)+
\epsilon(\alpha)\sum_{\beta\in\n(\alpha)}
\frac{u_{\beta}-v_{\alpha}(\lambda)}{e_{\alpha\beta}+\lambda}
\sum_{1\le s\le r\atop \alpha_s\succeq\beta}
2|u_{\alpha_s}|^{-1}\wt E''_s,\tag 6.2$$
where $\wt E''_s$ is the $S^2D''_s$-components of the
horizontal lift of the energy
function $E''$ of $M''$, and $\epsilon(\alpha)$ denotes the
sign of $u_{\alpha}$.
Then, taking the $S^3D'$-components of
$$0=\{F_{\alpha}(\lambda), F_{\alpha}(\mu)\},$$
we obtain
$$0=\{F'_{\alpha}(\lambda), F'_{\alpha}(\mu)\}.$$
Hence $\Cal F'$ is commutative.

(2) is an immediate consequence of Lemma 6.1.
\qed
\enddemo

The typical fibre $X(\D')$ is naturally identified with
the fibre
$\pi_{\sharp}^{-1}(\pi_{\sharp}(p_0))$
passing through the base point $p_0\in M^1$. Denoting
the K\"ahler manifold
$\pi_{\sharp}^{-1}(\pi_{\sharp}(p_0))$ by $M'$,
and $\Cal F'(\pi_{\sharp}(p_0))$
by $\Cal F'$, we obtain a K\"ahler-Liouville manifold
$(M',\Cal F')$.
Note that it is of type (A)
if and only if every maximal element $\alpha$
of $\A'$ satisfies $|\alpha|\ge 2$.
The following theorem is immediate.

\proclaim{Theorem 6.5} If $(M',\Cal F')$ is of type
(A), then it possesses the
following properties:
\roster
\item"(1)" The associated partially ordered set is naturally
identified with $\A'$;
\item"(2)" the underlying toric variety is isomorphic to
$X(\D')$;
\item"(3)" the fundamental functions
$\{h'_{\alpha,\nu}\}$ $(\alpha\in\A',\ 1\le
\nu\le |\alpha|)$ are given by the restriction of
$h_{\alpha,\nu}$ to $M'$;
\item"(4)" the conjunction constants
$e'_{\alpha\beta}$ $(\beta\in\A',\ \alpha
\pc\beta)$ are given by $e'_{\alpha\beta}=e_{\alpha\beta}$;
\item"(5)" the scaling constants $d'_{\alpha}$ $(\alpha\in\A')$ are given
by
$d'_{\alpha}=d_{\alpha}$.
\endroster
\endproclaim

In case $(M',\Cal F')$ is not of type
(A), then the structure of toric variety on $M'$ may be
external, i.e., not determined by $(M',\Cal F')$ itself.
Nevertheless, we have the following

\proclaim{Proposition 6.6}
\roster
\runinitem"(1)" The maximal compact subgroup $K'$ of the
algebraic torus
$\Cal T_{\Gamma'}$ acts on the K\"ahler manifold $M'$
as automorphisms and preserves
each element of $\Cal F'$.
\item"(2)" The geodesic flow of $M'$ is integrable
by means of $\Cal F'$ and the Lie
algebra of $K'$.
\endroster
\endproclaim

The proof is clear. We shall say that
a compact K\"ahler-Liouville manifold
is {\it of type (B)} if it can
be realized as the fibre of a fibre bundle obtained from a
compact K\"ahler-Liouville manifold of type (A) and an
open subset of the associated
partially ordered set, and if it is not
of type (A). By the definition, it possesses a
structure of toric variety of KL-B type (not
necessarily unique). It is another type of
K\"ahler-Liouville manifold whose geodesic flow is integrable. In
this paper we shall not
mention further about such a manifold except the
1-dimensional case (see Section 8).

In the rest of this section, we shall
show that the K\"ahler-Liouville manifold
$(M,\Cal F)$ can be reconstructed from the structure
of toric variety on $M$ and the
K\"ahler-Liouville manifolds $(M'',\Cal F'')$ and $(M',\Cal F')$, provided
$(M',\Cal F')$
is of type (A). By virtue of Corollary
5.6, $M''$ is described
as the product $M''_1\times\dots\times M''_r$
of K\"ahler-Liouville manifolds,
corresponding to the decomposition of $\A''$ into the
connected components. Let
$\omega''_s$ denote the K\"ahler form of $M''_s$. Also,
let $\alpha_s$ denote the
unique minimal element of $\A''_s$.

Put
$$Q=\cup_{q\in M''} (\text{the unique open orbit of }\Cal
T_{\Gamma'}\text{ in
the fibre }\pi_{\sharp}^{-1}(q)).$$
Then $Q$ is open and dense in $M$,
and $\pi_{\sharp}: Q\to M''$ is a principal
$\Cal T_{\Gamma'}$-bundle. Proposition 5.7 implies that this bundle
is isomorphic
to
$$X(\wt\D'')\times_{\psi}\Cal T_{\Gamma'}\to M''.$$
Since the horizontal subbundle $D''$ is $\Cal T_{\Gamma'}$-invariant,
it defines
the connection on this principal bundle. Let $\theta$
be the connection form, and
$\Theta$ the $\frak g'$-valued 2-form on $M''$ so
that $\pi_{\sharp}^*\Theta$ is the
curvature form ($\frak g'$ is the Lie algebra
of $\Cal T_{\Gamma'}$).

\proclaim{Proposition 6.7}
$$\Theta=\sum_{s=1}^r d''_{\alpha_s} \omega''_s\otimes Z_{\alpha_s}.$$
\endproclaim
\demo{Proof} By Propositions 1.2 and 1.10 we have
$$[V_i,IV_i]\equiv -\sgrad(\log |u_{\alpha}|)\quad \mod (D_i)$$
for $\alpha\in\A''$ and $i\in\alpha$. This implies
$$[V''_i,IV''_i]\equiv -\sgrad |u_{\alpha_s}| \quad \mod (D'')$$
for any $\alpha\in\A''_s$ and $i\in\alpha$.
Since $|d_{\alpha_s}|=d''_{\alpha_s}$,
we have
$$d\theta (V''_i,IV''_i)= d''_{\alpha_s} Z_{\alpha_s}\qquad
(\alpha\in\A''_s,\, i\in\alpha).$$
Also, it is easily seen that
$$d\theta (V''_i,V''_j)=d\theta (V''_i,IV''_j)=d\theta (IV''_i,IV''_j)=0$$
for any $\alpha,\beta\in\A''$ and $i\in\alpha,\,j\in\beta$ $(i\ne j)$.
Hence the proposition follows.
\qed
\enddemo

{}From now on, we forget the structure of
K\"ahler-Liouville manifold, and only assume
that $M=X(\D)$ is a toric variety of KL-A
type. Let $\A$ be the associated partially
ordered set, and let $\A'$ be an open
subset of it. Put $\A''=\A-\A'$. Then we have the
fibre bundle $\pi_{\sharp}:M\to X(\D'')$ with typical fibre $X(\D')$
as before, and the
principal $\Cal T_{\Gamma'}$-bundle $\pi_{\sharp}:Q\to X(\D'')$ as above.
Let $(M'',\Cal F'')$ be a K\"ahler-Liouville manifold of
type (A) whose underlying toric
variety is isomorphic to $X(\D'')$. We shall identify
$M''$ with $X(\D'')$.
Let $\A_s$, $\alpha_s$, and $\omega''_s$ $(1\le s\le r)$
be as above.
For each non-minimal $\alpha\in\A$, let $l_{\alpha}$ be the
largest positive
integer satisfying $l_{\alpha}^{-1} Z_{\alpha}\in\Gamma$.

\proclaim{Lemma 6.8}
\roster
\runinitem"(1)" For any $\alpha\in\A''_s$,
$$\left(\prod_{\alpha_s\pcc\beta\pc\alpha}
m_{\beta,0}\right)l_{\alpha_s}\in\Z.$$
\item"(2)" $$\left[\frac{d''_{\alpha_s}}{2\pi}l_{\alpha_s}\omega''_s\right]
\in H^2(M''_s,\Z).$$
\endroster
\endproclaim
\demo{Proof} (1) By Proposition 4.20 we have
$$Z_{\alpha}\equiv \left(\prod_{\alpha_s\pcc\beta\pc\alpha}m_{\beta,0}\right)
Z_{\alpha_s}\qquad \mod (\sum_{\beta\in\A''}\sum_{\nu=1}^{|\beta|}\Z
Y_{\beta,\nu}).$$
Thus the assertion follows. (2) follows from (1)
and Theorem 5.10 (2).
\qed
\enddemo

\proclaim{Proposition 6.9} There is a unique connection form
$\theta$ on the principal
bundle $\pi_{\sharp}:Q\to M''$ such that
\roster
\item"(1)" the associated curvature form is given by
$\pi_{\sharp}^*\Theta$,
where
$$\Theta=\sum_{s=1}^r d''_{\alpha_s} \omega''_s\otimes Z_{\alpha_s};$$
\item"(2)" the horizontal distribution defined by the kernel
of $\theta$ is invariant
with respect to the complex structure $I$.
\endroster

Accordingly, the $\Cal T_{\Gamma'}$-invariant horizontal subbundle $D''$ of
$TM$
with respect to $\pi_{\sharp}:M\to M''$ is uniquely determined
by $\Theta$ so that
the connection $D''|_Q$ on $Q$ satisfies the conditions
above. The connection form
$\theta$ and the subbundle $D''$ are invariant under
the action of the maximal compact
subgroup $K$ of $\Cal T_{\Gamma}$.
\endproclaim
\demo{Proof} First, we shall prove the uniqueness. Let
$\theta$ be a connection form
satisfying the condition (1) and (2). $\theta$ is
$\frak g'$-valued, and here $\frak g'$
is regarded as a real Lie algebra with
the complex structure $I$. Now, we regard it
as a complex Lie algebra by replacing $I$
with $\sqrt{-1}$. We shall write $\wt \theta$
(resp. $\wt\Theta$) instead of $\theta$ (resp. $\Theta$) when
$\frak g'$ are regarded as
the complex Lie algebra. By extending it $\C$-linearly
to $TQ\otimes \C$, $\wt\theta$
becomes a $(1,0)$-form. Since $\wt\Theta$ is a $(1,1)$-from,
we have
$$\partial\wt\theta=0, \quad \bar\partial\wt\theta=\pi_{\sharp}^*\wt\Theta.$$
This implies that if $\theta_1$ is another connection
form satisfying the conditions (1)
and (2), then $\wt\theta-\wt\theta_1$ is a holomorphic 1-form,
and is projectable.
Hence there is a holomorphic 1-form $\mu$ on
$M''$ such that $\wt\theta-\wt\theta_1=
\pi_{\sharp}^*\mu$. However, since $M''$ is a compact, simply
connected K\"ahler
manifold, we have $\mu=0$. Thus it follows that
$\theta=\theta_1$.

Next, we shall prove the existence. Let $P_s$
be a hermitian line bundle over $M''_s$
with the canonical hermitian connection form $\wt\theta_s$ whose
first Chern form is
equal to
$$-\frac{d''_{\alpha_s}}{2\pi} l_{\alpha_s}\omega''_s$$
(for the existence of such a hermitian line
bundle, see [7] p.41, Proposition).
Put $U_s=\{v\in P_s\ |\ |v|=1\}$ and $U=\prod_{s=1}^r U_s$.
Then
$U$ is a principal $U(1)^r$-bundle over $M''=\prod_{s=1}^r M''_s$.
Let
$\phi:U(1)^r\to \Cal T_{\Gamma'}$ be the homomorphism given by
$$(\lambda_1,\dots,\lambda_r)\mapsto \prod_{s=1}^r(l_{\alpha_s}^{-1}
Z_{\alpha_s}\otimes \lambda_s).$$
Then we obtain the associated $\Cal T_{\Gamma'}$-bundle
$U\times_{\phi} \Cal T_{\Gamma'}\to M''$.

\proclaim{Lemma 6.10} The principal bundle
$U\times_{\phi} \Cal T_{\Gamma'}\to
M''$
is naturally identified with the bundle $\pi_{\sharp}: Q\to
M''$.
\endproclaim
\demo{Proof} Let $\chi_s:\Gamma''_0\to\Z$ $(1\le s\le r)$ be the
homomorphism
given by
$$\chi_s(\wt\pi(R_{\alpha}))=\left(\prod_{\alpha_s\pcc\beta\pc\alpha}
m_{\beta,0}\right)l_{\alpha_s}.$$
The associated homomorphism $\Cal T_{\Gamma''_0}\to \C^{\times}$ is also
denoted by
$\chi_s$. Then, by Proposition 5.8 and Theorem 5.10
we see that the line bundle $P_s$
is isomorphic to the fibre product
$$X(\wt\D'')\times_{\chi_s}\C\to M''.$$
Moreover, denoting by $K''_0$ the maximal compact subgroup
of $\Cal T_{\Gamma''_0}$,
we have
$$\psi|_{K''_0}=\phi\circ (\chi_1,\dots,\chi_r)|_{K''_0}.$$
Therefore the lemma follows.
\qed
\enddemo

We now continue the proof of Proposition 6.9.
The direct sum of the connection forms
$\wt\theta_s$, restricted to $U_s$, is a connection form
on $U$. Composing this with the
Lie algebra homomorphism associated with $\phi$, we obtain
a connection form
$\theta$ on the principal bundle $\pi_{\sharp}: Q\to M''$.
Then we clearly have
$d\theta=\pi_{\sharp}^*\Theta$.

Finally, we prove the $K$-invariance. Let $k\in K$.
Then the pull-back $k^*\theta$ is a
connection form with the same curvature, because $\Theta$
is preserved by the
transformation of $M''$ induced from $k$. Hence by
the uniqueness we have $k^*\theta
=\theta$. This completes the proof.
\qed
\enddemo

Now, we moreover assume that there is a
K\"ahler-Liouville manifold $(M',\Cal F')$
of type (A) whose underlying toric variety is
$X(\D')$. Then we have the following

\proclaim{Theorem 6.11} There is a unique K\"ahler-Liouville manifold
$(M,\Cal F)$
of type (A) satisfying the following conditions:
\roster
\item"(1)" The underlying structure of toric variety is
identical with the given one;
\item"(2)" the given K\"ahler-Liouville manifolds $(M',\Cal F')$ and
$(M'',\Cal F'')$
are isomorphic with the ones induced from $(M,\Cal
F)$.
\endroster
\endproclaim
\demo{Proof} We first define functions $h_{\alpha,\nu}$ $(\alpha\in\A,
1\le\nu\le |\alpha|)$ on $M$. Let $K'$ be the
maximal compact subgroup of
$\Cal T_{\Gamma'}$. Since the fundamental functions $\{h'_{\alpha,\nu}\}$ of
$(M',\Cal F')$ are $K'$-invariant, they are supposed to
be defined on
$M=U\times_{U(1)^r}M'$. So, we put
$$h_{\alpha,\nu}=\cases h'_{\alpha,\nu}\quad & (\alpha\in\A')\\
\pi_{\sharp}^*h''_{\alpha,\nu}\quad & (\alpha\in\A'')
\endcases$$
Accordingly, we also put $c_{\alpha,\nu}=
c'_{\alpha,\nu}$ if $\alpha\in\A'$ and
$=c''_{\alpha,\nu}$ if $\alpha\in\A''$,
where $c'_{\alpha,\nu}$ and $c''_{\alpha,\nu}$
are the fundamental constants.

We choose the ordering of
$Y_{\alpha,0},\dots, Y_{\alpha,|\alpha|}$ $(\alpha\in\A)$
so that the ordering of $Y''_{\alpha,\nu}=
(\pi_{\sharp})_*Y_{\alpha,\nu}$ $(\alpha\in\A'')$
and $Y'_{\alpha,\nu}=Y_{\alpha,\nu}$ $(\alpha\in\A')$ are equal to the ones
induced from the fundamental functions
$\{h''_{\alpha,\nu}\}$ and $\{h'_{\alpha,\nu}\}$ respectively
(cf. the remark before Proposition 5.1). Hence the
numbers $m_{\alpha,\nu}$ $(\alpha,
\p(\alpha)$, non-minimal) and
$m_{\alpha,\nu}-m_{\alpha,0}$ $(\p(\alpha)$, minimal) are
uniquely determined. We define the number $m_{\alpha,0}$ for
$\alpha$ such that
$\p(\alpha)$ is minimal by the formula (4.9) and
the value of the constants
$c_{\alpha,\nu}$. Then again by (4.9) the value of
$e_{\beta\alpha}$ is determined
for every non-maximal $\beta$ and $\alpha\in\n(\beta)$. It is
also defined for all
$\alpha$ and $\beta\pcc\alpha$ so that Proposition 1.10 (2)
is satisfied.

Now, let us define the function $u_{\alpha}$ on
$M$ by
$$u_{\alpha}=\prod_{\beta\pcc\alpha}\prod_{\nu=1}^{|\beta|}
(h_{\beta,\nu}+e_{\beta\alpha}).$$
Let $D''$ be the (horizontal) subbundle of $TM$
given by Proposition 6.9.
Since the K\"ahler form $\omega'$ on $M'$ is
supposed to be
defined on $M$ so that the kernel coincides
with the horizontal subbundle $D''$,
we can define a 2-form $\omega$ on $M$
by
$$\omega= \omega' + \sum_{s=1}^r |u_{\alpha_s}| \pi_{\sharp}^*\omega''_s.$$
Also, we define $F_{\alpha}(\lambda)$ $(\alpha\in\A)$ by the horizontal
lift of
the generating polynomial $F''_{\alpha}(\lambda)$ of $(M''\Cal F'')$ if
$\alpha\in\A''$,
and by the formula (6.2) if $\alpha\in\A'$, where
$\wt E''_s$ is the horizontal lift of
the energy function $E''_s$ of $M''_s$, and $F'_{\alpha}(\lambda)$
is the generating
polynomial of $(M',\Cal F')$. Let $\Cal F$ be
the vector space spanned by all the
coefficients of $F_{\alpha}(\lambda)$ $(\alpha\in\A)$.

Define the orthonormal frame $V_{\alpha,\nu}, IV_{\alpha,\nu}$ on the
open
$\Cal T_{\Gamma}$-orbit $M^1$ by using the corresponding frames
on $M'$ and $M''$ in the
obvious manner. Then we have the formula (6.1)
and the relations
$$\gather
[W_{\alpha,\nu},W_{\beta,\mu}]=[W_{\alpha,\nu},IW_{\beta,\mu}]
=[IW_{\alpha,\nu},IW_{\beta,\mu}]=0\quad ((\alpha,\nu)\ne(\beta,\mu))\\
[V_{\alpha,\nu},IV_{\alpha,\nu}]
\equiv \sgrad a_{\alpha,\nu}\quad \mod D_{\alpha,\nu},
\endgather$$
where $W_{\alpha,\nu}= a_{\alpha,\nu}^{-1/2} V_{\alpha,\nu}$ and
$$a_{\alpha,\nu}^{-1}=|u_{\alpha}|
\prod_{\mu\ne\nu}|h_{\alpha,\mu}-h_{\alpha,\nu}|.$$
Hence the arguments in Section 1 imply that
$\omega$ is a K\"ahler form, and with
this K\"ahler metric $(M,\Cal F)$ becomes a K\"ahler-Liouville
manifold of type (A).
The uniqueness and the property (2) obviously follow
from the construction above.
To prove (1) we need the following lemma.

\proclaim{Lemma 6.12}
$$Y_{\alpha,\nu}=\wt Y''_{\alpha,\nu}+\frac{d''_{\alpha_s}}{d''_{\alpha,\nu}}
\frac{v_{\alpha}(c_{\alpha,\nu})}
{u_{\alpha_s}} Z_{\alpha_s}\qquad (\alpha\in\A''_s),$$
where $\wt Y''_{\alpha,\nu}$ is the horizontal lift of
$Y''_{\alpha,\nu}$.
\endproclaim
\demo{Proof} Let $\theta$ be the connection form given
by Proposition 6.9. Then we have
$$i_{Y_{\alpha,\nu}}d\theta=d''_{\alpha_s}
\pi_{\sharp}^*(i_{Y''_{\alpha,\nu}}\omega''_s)
\otimes Z_{\alpha_s}=-\frac{d''_{\alpha_s}}{d''_{\alpha,\nu}}
d\left(\frac{v_{\alpha}(c_{\alpha,\nu})}
{u_{\alpha_s}}\right)\otimes Z_{\alpha_s}.$$
Since the left-hand side is equal to $-d(\theta(Y_{\alpha,\nu}))$
by Proposition 6.9, we have
$$Y_{\alpha,\nu}=\wt Y''_{\alpha,\nu}+
\frac{d''_{\alpha_s}}{d''_{\alpha,\nu}}
\frac{v_{\alpha}(c_{\alpha,\nu})}{u_{\alpha_s}}
Z_{\alpha_s}+(\text{constant term}).$$
Then, by comparing both sides at points on
$L_{\alpha,\nu}=V(\R_{\ge 0}Y_{\alpha,\nu})$,
the lemma is proved.
\qed
\enddemo

The lemma above implies that $Y_{\alpha,\nu}=
d_{\alpha,\nu}^{-1}\sgrad v_{\alpha}
(c_{\alpha,\nu})$
for any $\alpha\in\A''_s$, where
$d_{\alpha,\nu}=\epsilon(\alpha_s)d''_{\alpha,\nu}$ and
$\epsilon(\alpha_s)$ is the sign of $u_{\alpha_s}$. Since it
is also true for $\alpha\in\A'$
($d_{\alpha,\nu}=d'_{\alpha,\nu}$ in this case),
the condition (1) is therefore
satisfied.
This completes the proof of Theorem 6.11.
\qed
\enddemo

\specialhead 7. The case where $\#\A=1$
\endspecialhead
In this section we shall classify compact K\"ahler-Liouville
manifolds (of type (A)) such that $\#\A=1$. Note
that
such manifolds are isomorphic to the complex projective
space $\C P^n$ (with the
standard $(\C^{\times})^n=(\C^{\times})^{n+1}/\C^{\times}$ action) as toric
variety.

Let $(M,\Cal F)$ be a compact K\"ahler-Liouville manifold
of type
(A) such that the associated partially ordered set
$\A$ consists of one element.
In this case we write $\nu$ instead of
$(\alpha,\nu)$, and $d_*$ instead of
$d_{\alpha}$. Put
$$S=\cap_{\nu=1}^{n-1} L_{\nu},\quad \{q_1\}=L_0\cap S,\quad
\{q_2\}=L_n\cap S.$$
$S$ is holomorphically isomorphic to $\C P^1$. We
regard $S$ as a K\"ahler
manifold with the induced metric. Clearly, $Y_n$ is
tangent to $S$, and its
zeros are $q_1$ and $q_2$. Let $\gamma(t)$ $(0\le
t\le l/2)$ be a minimal
geodesic of unit speed such that $\gamma(0)=q_1$ and
$\gamma(l/2)=q_2$.
Since $S$ is a surface of revolution, $\gamma(t)$
is extended to a closed
geodesic of the least period $l$. Recalling the
function $v(\lambda)=
\prod_{\nu=1}^n (h_{\nu}-\lambda)$, put
$$h(t)=\frac{v(c_n)(\gamma(t))}{\prod_{\nu=0}^{n-1}(c_{\nu}-c_n)}
\qquad t\in\R/l\Z.$$

\proclaim{Proposition 7.1} $h\in C^{\infty}(\R/l\Z)$ possesses the following
properties:
\roster
\item"(1)" $h(-t)=h(t)$ for any $t$;
\item"(2)" $h(0)=1$, $h(l/2)=0$;
\item"(3)" $h'(t)<0$ if $0<t<l/2$;
\item"(4)" $h'(T_{\nu})=-\sqrt{2d_* c_{\nu}(1-c_{\nu})}$ $(1\le\nu\le n-1)$,
where $T_{\nu}$ $(0<T_{\nu}<l/2)$ is defined by $h(T_{\nu})=c_{\nu}$;
\item"(5)" $-h''(0)=h''(l/2)=d_*$.
\endroster
\endproclaim
\demo{Proof (except (4))} (1), (2), and (5) are
clear. Since $Y_n\ne 0$ at
$\gamma(t)$ $(0<t<l/2)$, (3) is also obvious.
\qed
\enddemo

Let $\Cal C=\Cal C_n$ be the set of
elements $(\{c_0,\dots,c_n\},d_*,l,h)$
such that $d_*$ and $l$ are positive constants,
$\{c_{\nu}\}$ are constants
satisfying
$$1=c_0>c_1>\dots>c_n=0,\tag 7.1$$
and $h\in C^{\infty}(\R/l\Z)$ satisfies the conditions $(1),\dots,
(5)$ in Lemma 7.1. We say that two
elements
$(\{c_{\nu}\},d_*,l,h)$ and $(\{\wt c_{\nu}\},\wt d_*, \wt l,\wt h)$
are
{\it equivalent} if $\wt d_*=d_*$, $\wt l=l$ and
either $\wt c_{\nu}=c_{\nu}$,
$\wt h=h$, or
$$\wt c_{\nu}=1-c_{n-\nu},\qquad \wt h(t)= 1-h(\frac{l}2 -t).$$

\proclaim{Theorem 7.2} The assignment of $(\{c_{\nu}\},d_*, l,h)\in\Cal C$
to
$(M,\Cal F)$ described above gives the one-to-one correspondence
between
the set of the isomorphism classes of compact
K\"ahler-Liouville manifolds
of type (A) satisfying $\#\A=1$ and the equivalence
classes of elements
of $\Cal C$.
\endproclaim

To prove Proposition 7.1 (4) and Theorem 7.2
we shall use the results for (real)
Liouville manifolds obtained by the author [6]. First,
we prove the following

\proclaim{Proposition 7.3} Let $p_0\in M^1$ be the base
point so that $M$ is
identified with the toric variety $X(\D)$. Put
$$N=\Exp_{p_0}(D^+).$$
Then, $N$ is a well-defined real submanifold of
$M$, which is totally geodesic
and diffeomorphic to $\R P^n$.
Moreover, take the $S^2TN$-component $F'$ of each element
$F\in\Cal F$ and put
$$\Cal F'=\{F'\ |\ F\in\Cal F\}.$$
Then $(N,\Cal F')$ is a proper Liouville manifold
of rank one and type {\rm (B)},
and its core is isomorphic to
$$(\R/l\Z,\{[h-c_1],\dots, [h-c_{n-1}]\}).$$
\endproclaim
\demo{Proof} Using the real number field $\R$ instead
of $\C$ in the construction
of the toric variety $X(\D)$ (cf. Section 5),
one obtains a submanifold $X(\D)(\R)$
diffeomorphic to the $n$-dimensional real projective space $\R
P^n$. Since
its tangent space is spanned by $I\frak k=D^+$
at each point on $\R P^n\cap M^1$,
Proposition 1.6 implies that $X(\D)(\R)$ is totally geodesic
and $N=X(\D)(\R)$.

It is easy to verify that $(N,\Cal F')$
is a Liouville manifold. Put
$$G_{\nu}=\sum_{\xi=1}^n \left(\prod_{\mu\ne\xi}(h_{\mu}-c_{\nu})\right)
(V_{\xi}^2 +(IV_{\xi})^2)\in\Cal F \qquad (1\le\nu\le n-1).$$
Then we have
$$\{p\in N\ |\ (G'_{\nu})_p=0\}=M^s\cap L_{\nu}\cap N,$$
$$\{p\in N\ |\ \text{rank\,}(G'_{\nu})_p\le1\}=L_{\nu}\cap N.$$
Also, we have $(dG'_{\nu})_{\lambda}\ne 0$ at some $\lambda\in
T^*_pN$ for every
$p\in M^s\cap L_{\nu}\cap N$, because $d(h_{\nu}+h_{\nu+1})$ does not
vanish
at $p$. Hence the Liouville manifold $(N,\Cal F')$
is proper and of rank one. Since
$N$ is diffeomorphic to $\R P^n$, it is
of type (B) (cf. [6] Theorems $3.3.1$ and
$3.4.1$).

By definition the core of the Liouville manifold
$(N,\Cal F')$ consists of the
1 dimensional riemannian submanifold
$$C=\cap_{\nu=1}^{n-1} (M^s\cap L_{\nu}\cap N)$$
(called the core submanifold) and the equivalence classes
$[\wt h_{\nu}]$ of the
functions $\wt h_{\nu}$ on it defined by
$$(G'_{\nu})_p=\wt h_{\nu}(p) V^2,\qquad p\in C$$
where $V$ is the unit normal vector to
$L_{\nu}\cap N$. The equivalence classes
of functions are defined to be the orbits
of the affine transformation group on
the target space $\R$. It is easily seen
that $C$ is
equal to the image of a closed geodesic
passing through $q_1$ and $q_2$.
Hence we can take $\gamma$ so that its
image coincides with $C$.
Thus $C$ is isometric to $\R/l\Z$, and $\wt
h_{\nu}(\gamma(t))=h(t)-c_{\nu}$.
This completes the proof.
\qed
\enddemo

Note that another choice of the base point
$p_0$ gives another submanifold, but
they are mutually transferred with the action of
$K$. Hence the isomorphism
class of the Liouville manifold $(N,\Cal F')$ is
uniquely determined.
As was shown in [6], the isomorphism classes
of proper Liouville manifolds of
rank one are completely classified by means of
the isomorphism classes of
the cores. In the present case, two cores
$(\R/l\Z,\{[h-c_1],\dots, [h-c_{n-1}]\})$
and $(\R/\wt l\Z,\{[\wt h-\wt c_1],\dots, [\wt h-\wt c_{n-1}]\})$
are isomorphic
if and only if $l=\wt l$ and either
$h(t)=\wt h(t)$, $c_{\nu}=\wt c_{\nu}$, or
$h(t)=1-\wt h(-t+l/2)$, $c_{\nu}=1-\wt c_{n-\nu}$. Hence those isomorphism
classes corresponds to the equivalence classes of elements
of $\Cal C$.

By the proof of Theorem 3.3.1 and Theorem
3.4.1 in [6] we obtain a branched
covering of $N$ whose covering space is a
torus. We now explain it: Put
$$\alpha_{\nu}=4\int_{T_{\nu-1}}^{T_{\nu}}\frac{dt}{\sqrt{(-1)^{\nu-1}
\prod_{\mu=1}^{n-1} (h(t)-c_{\mu})}}\qquad (1\le \nu\le n),\tag 7.1$$
where $T_{\nu}\in [0,l/2]$ is defined by $h(T_{\nu})=c_{\nu}$ $(0\le\nu\le
n)$.
Put
$$R=\prod_{\nu=1}^n \left(\R/\alpha_{\nu}\Z\right),$$
and let $x_{\nu}$ (mod$\, \alpha_{\nu}\Z)$ be the natural
coordinate of
$\R/\alpha_{\nu}\Z$. Let $H(\simeq (\Z/2\Z)^{n})$ be the transformation
group of $R$ generated by $\tau_{2\nu}\circ \tau_{2\nu+1}$ $(1\le\nu\le
n-1)$
and $\tau_1\circ\prod_{\nu=1}^n \tau_{2\nu}$, where
$$\aligned
& \tau_{2\nu-1}(x_1,\dots,x_n)=
(x_1,\dots,x_{\nu-1},\frac{\alpha_{\nu}}2-x_{\nu},
x_{\nu+1},\dots,x_n)\\
& \tau_{2\nu}(x_1,\dots,x_n)=
(x_1,\dots,x_{\nu-1},-x_{\nu},x_{\nu+1},\dots,x_n).
\endaligned\tag 7.2$$
Then we have

\proclaim{Proposition 7.4 ([6])} There is a surjective mapping
$\Phi:R\to N$
possessing the following properties:
\roster
\item"(1)" For any $p\in N$, $\Phi^{-1}(p)$ is an
$H$-orbit;
\item"(2)" $\Phi_*(\partial/\partial x_{\nu})=\pm W_{\nu}$;
\item"(3)" $h_{\nu}\circ\Phi$ are $C^{\infty}$ functions;
\item"(4)" $M^s\cap N=\{p\in N\ |\ \#\Phi^{-1}(p)<2^n\}$;
\item"(5)" $L_{\nu}\cap N=\{\Phi(x)\ |\ x_{\nu}=0,\alpha_{\nu}/2
\text{ or } x_{\nu+1}=\pm \alpha_{\nu+1}/4\}\quad (1\le\nu\le n-1)$,
$L_0\cap N=\{\Phi(x)\ |\ x_1=\pm \frac{\alpha_1}4\}$,
$L_n\cap N=\{\Phi(x)\ |\ x_n=0,\frac{\alpha_n}2\}$;
\item"(6)" $\Phi\circ\tau_{2\nu-1}=\exp(\pi Y_{\nu-1})\circ\Phi,\quad
\Phi\circ\tau_{2\nu}=\exp(\pi Y_{\nu})\circ\Phi$.
\endroster
\endproclaim

\demo{Proof of Proposition 7.1 (4)} Since the function
$h_{\nu}\circ\Phi$ depends
only on the variable $x_{\nu}$, we write it
$\wt h_{\nu}(x_{\nu})$. Observing
the formula
$$\hess v(c_{\nu})=d_*\prod_{0\le\mu\le n\atop\mu\ne \nu}(c_{\mu}-c_{\nu})$$
at a point $p$ such that $h_{\nu}(p)=c_{\nu}$ and
$h_{\nu+1}(p)\ne c_{\nu}$,
we have $\wt h'_{\nu}(0)=0$ and
$$\wt h''_{\nu}(0)=(-1)^{n-\nu}d_*\prod_{\mu\ne \nu}
(c_{\mu}-c_{\nu}).$$
Note that the vector fields $V_{\nu}$ and $W_{\nu}$
are locally well-defined (up to
sign) and smooth as vector fields on $M^0\cap
N$, though they are not determined
around $p\not\in M^1$ as vector fields on $M$.
Since $h(t)=h_{\nu}(\gamma(t))$ on $[T_{\nu-1},T_{\nu}]$, we have
$$\align
h'(T_{\nu})^2 & =\lim_{t\nearrow T_{\nu}}(V_{\nu}h_{\nu})^2(\gamma(t))\\
& =\lim_{t\nearrow T_{\nu}}\frac{(W_{\nu}h_{\nu})^2(\gamma(t))}
{(-1)^{n-\nu}\prod_{\mu\ne\nu}(h_{\mu}-h_{\nu})}\\
& =\lim_{x_{\nu}\to 0}\frac{(\wt h'_{\nu}(x_{\nu}))^2}
{(-1)^{n-\nu}\prod_{\mu=1}^{n-1}(c_{\mu}-\wt h_{\nu}(x_{\nu}))}\\
& =2d_*(1-c_{\nu})c_{\nu}.
\endalign$$
\qed
\enddemo

\demo{Proof of Theorem 7.2} Let $(\{c_{\nu}\},d_*, l,h)\in\Cal C$
be an arbitrary
element, and let $(N, \Cal F')$ be a
proper Liouville manifold of rank one
whose core is isomorphic to
$$(\R/l\Z,\{[h-c_1],\dots,[h-c_{n-1}]\}).$$
To prove Theorem 7.2 it suffices to show
that there is a unique K\"ahler-Liouville
manifold $(M,\Cal F)$ up to isomorphism such that
the associated Liouville manifold
is isomorphic to $(N,\Cal F')$. To do so,
we first review how to construct
$(N, \Cal F')$.

Let $\alpha_{\nu}$, $R$, $\tau_{2\nu-1},\, \tau_{2\nu}$, and $H$ be
as above. It is
not hard to see that R/H is homeomorphic
to $\R P^n$ with the quotient topology.
Put $N=R/H$, and let $\Phi:R\to N$ be the
quotient mapping. To regard $N$ as
differentiable manifold diffeomorphic to $\R P^n$, it is
necessary to specify
coordinate systems around branch points (i.e., points $p\in
N$ such that
$\#\Phi^{-1}(p)<2^n$). Let $N^s$ denote the branch locus. Put
$$\align
I_{\nu} & =\{\Phi(x)\ |\ \tau_{2\nu}(x)=
x \text{ and }\tau_{2\nu+1}(x)=x\}\quad
(1\le\nu \le n-1)\\
J_{\nu} & =\{\Phi(x)\ |\ \tau_{2\nu}(x)=
x \text{ or }\tau_{2\nu+1}(x)=x\}\quad
(0\le\nu\le n)
\endalign$$
Then $N^s=\cup_{\nu=1}^{n-1}I_{\nu}$. Let $p=\Phi(a)\in N^s$. Then there is
a unique
subset $K$ of $\{1,\dots,n-1\}$ such that $p\in I_{\nu}$
if and only if $\nu\in K$.
Writing $K=\{\nu_1,\dots,\nu_k\}$, $\nu_1<\dots<\nu_k$, we have
$\nu_{i+1}-\nu_i\ge 2$. Define functions $y_1,\dots,y_n$ by
$$\align
& y_{\nu_i} =(x_{\nu_i}-a_{\nu_i})^2+(x_{\nu_i+1}-a_{\nu_i+1})^2
\quad (1\le i\le k)\\
& y_{\nu_i+1} =2(x_{\nu_i}-a_{\nu_i})(x_{\nu_i+1}-a_{\nu_i+1})
\quad (1\le i\le k)\\
& y_{\nu}=x_{\nu}\qquad (\nu,\nu-1\not\in K)
\endalign$$
The system of functions $(y_{\nu})$ is then projectable,
and becomes a coordinate
system around $p$ mentioned above (cf. [6] Proposition
3.3.2).

Define $C^{\infty}$ mappings
$$\R/\alpha_{\nu}\Z\to\cases
[-T_1,T_1]\quad & (\nu=1)\\
[T_{\nu-1},T_{\nu}]\quad & (2\le\nu\le n-1)\\
[T_{n-1},l-T_{n-1}]\quad & (\nu=n)
\endcases ,
\qquad (x_{\nu}\mapsto t=t_{\nu}(x_{\nu}))$$
by the differential equations
$$t'_{\nu}(x_{\nu})^2=(-1)^{\nu-1}\prod_{\mu=1}^{n-1}
(h(t_{\nu})-c_{\mu})$$
and the initial conditions
$$t_{\nu}(0)=T_{\nu}\quad (1\le\nu\le n),
\quad \cases t'(0)=0,\ t''_{\nu}(0)<0\quad & (1\le\nu\le n-1)\\
t'_n(0)<0\quad & (\nu=n).\endcases$$
Put $\wt h_{\nu}(x_{\nu})=h(t_{\nu}(x_{\nu}))$ and
$$\align
g' & =\sum_{\nu=1}^n (-1)^{n-\nu}
\left(\prod_{\mu\ne\nu}(\wt h_{\mu}-\wt h_{\nu})
\right) (dx_{\nu})^2\\
F'_{\nu} & =\sum_{\mu=1}^n \frac{\prod_{\xi\ne\mu}(\wt h_{\xi}-c_{\nu})}
{(-1)^{n-\mu}\prod_{\xi\ne\mu}(\wt h_{\xi}-\wt h_{\mu})} \left(\frac
{\partial}{\partial x_{\mu}}\right)^2\qquad (1\le\nu\le n-1).
\endalign$$
Then $g'$ and $F'_{\nu}$ are projectable, and define
the riemannian metric on $N$
and the sections of $S^2TN$ respectively. Denoting by
$E'$ the energy function
associated with $g'$ and by $\Cal F'$ the
vector space spanned by $F'_{\nu}$
$(1\le\nu\le n-1)$ and $E'$, we obtain the proper
Liouville manifold $(N,\Cal F')$
of rank one and type (B), whose core
is isomorphic to the given one.

The functions $\wt h_{\nu}$ are also projectable, and
define the continuous
functions $h_{\nu}$ on $N$. The function $h_{\nu}$ is
smooth outside
$I_{\nu}\cup I_{\nu-1}$. Also, it is easily seen that
the symmetric polynomials
of $h_1,\dots,h_n$ are smooth on the whole $N$.
Put $v(\lambda)=\prod_{\nu}
(h_{\nu}-\lambda)$, and
$$X_{\nu}=\frac{1}{d_*\prod_{0\le\mu\le n\atop \mu\ne\nu}(c_{\mu}-c_{\nu})}
\grad v(c_{\nu})\qquad (0\le\nu\le n).$$
The following lemma is immediate.

\proclaim{Lemma 7.5}
\roster
\runinitem"(1)" $[X_{\mu}, X_{\nu}]=0$ for any $\mu,\nu$.
\item"(2)" $v(c_{\nu})(p)=0$, $(X_{\nu})_p=0$ for $p\in J_{\nu}$.
\item"(3)" $\hess v(c_{\nu})$ on the normal bundle $NJ_{\nu}$
is equal to
$d_*\prod_{\mu\ne\nu}(c_{\mu}-c_{\nu})g'$.
\endroster
\endproclaim

Let $\pi:\R^{n+1}-\{0\}\to \R P^n$ be the natural projection,
and let
$(w_0,\dots,w_n)$ be the natural coordinate system of $\R^{n+1}$.

\proclaim{Proposition 7.6} There is a diffeomorphism $\phi: N\to
\R P^n$ such that
$$\phi_*(X_{\nu})=\pi_*\left(w_{\nu}\frac{\partial}{\partial w_{\nu}}\right)
\qquad (0\le\nu\le n).$$
\endproclaim
\demo{Proof} We first construct $\phi$ on $N-J_0$. Noting
that $X_1,\dots,X_n$
are linearly independent at every point on $N^1=N-\cup_{\nu=0}^n
J_{\nu}$, we define
(closed) 1-forms $\omega_1,\dots,\omega_n$
on $N^1$ by $\omega_{\nu}(X_{\mu})=
\delta_{\nu\mu}$. It is easily seen that $\omega_{\nu}$ is
smoothly extended to
$N-(J_0\cup J_{\nu})$. Fix a point $p_0\in N^1$

Let $\sigma_{\nu}$ $(0\le\nu\le n)$ be the involution on
$N$ defined by
$\Phi\circ\tau_{2\nu}=\sigma_{\nu}\circ\Phi$ or $\Phi\circ\tau_{2\nu+1}=
\sigma_{\nu}\circ\Phi$. Then, $\sigma_{\nu}$ is the reflection with
respect to $J_{\nu}$, and preserves each element of
$\Cal F'$. Also, we see that
$N-(J_0\cup J_{\nu})$ has two connected components; one contains
$p_0$ and the
other contains $\sigma_{\nu}(p_0)$.

Let $x_{\nu}$ $(1\le\nu\le n)$ be the function on
$N-(J_0\cup J_{\nu})$ defined by
$$x_{\nu}(p)=\cases \exp(\int_0^1 \omega_{\nu}(\dot c(t))dt)
\quad & (p\simeq p_0)\\
-\exp(\int_0^1 \omega_{\nu}(\dot c(t))dt)\quad & (p\simeq \sigma_{\nu}(p_0)),
\endcases$$
where $p\simeq p_0$ means $p$ and $p_0$ are
on the same component, and $c(t)$
$(0\le t\le 1)$ is a curve in $N-(J_0\cup
J_{\nu})$ from $p_0$ or
$\sigma_{\nu}(p_0)$ to $p$. Clearly we have
$$\gather
x_{\nu}(\sigma_{\nu}(p))=-x_{\nu}(p),\\
x_{\nu}(p_0)=1,\quad x_{\nu}(\sigma_{\nu}(p_0))=-1,\\
\omega_{\nu}=\frac{dx_{\nu}}{x_{\nu}}.
\endgather$$

We now claim that $x_{\nu}$ is smoothly extended
to $N-J_0$ by putting $x_{\nu}=0$
on $J_{\nu}-J_0$. In fact it is an easy
consequence of the following lemma.

\proclaim{Lemma 7.7} For each $p\in J_{\nu}-J_0$, there is
a neighborhood $U$ of
$p$ and a $C^{\infty}$ function $u_{\nu}$ on $U$
such that $u_{\nu}^2=|v(c_{\nu})|$.
\endproclaim

The lemma above follows from Lemma 7.5. Thus
we have obtained the
diffeomorphism
$$N-J_0\to \R^n\quad (p\mapsto (x_1(p),\dots,x_n(p)),\tag 7.3$$
which maps $X_{\nu}$ to $x_{\nu}\partial/\partial x_{\nu}$ and $p_0$
to
$(1,\dots,1)$ (the surjectivity follows from the completeness of
$X_{\nu}$
on $N^1$). Now, making the coordinate functions on
$N-J_{\nu}$ in the same
way, and gluing them together, we consequently obtain
the desired diffeomorphism
$\phi:N\to \R P^n$.
\qed
\enddemo

We now continue the proof of Theorem 7.2.
By virtue of Proposition 7.6 we may
identify $N$ with $\R P^n$. Hence $X_{\nu}=\pi_*(\partial/\partial w_{\nu})$,
and
$J_{\nu}$ is given by $w_{\nu}=0$. Also, we regard
$\R P^n$ as a submanifold of
$\C P^n$ in the natural manner. The projection
$\C^{n+1}-\{0\}\to \C P^n$ is also
denoted by $\pi$. Let $K=U(1)^n$ be the torus
acting on $\C P^n$ by
$$((\lambda_1,\dots,\lambda_n),\pi(w_0,\dots,w_n))\to
\pi(w_0,\lambda_1w_1,\dots,\lambda_nw_n)\quad (|\lambda_{\nu}|=1).$$
The following lemma is immediate.

\proclaim{Lemma 7.8} Let $H$ be a symmetric 2-form
on $\R P^n$ invariant with
respect to $\sigma_{\nu}$ $(0\le\nu\le n)$. Then there is
a unique hermitian
form $\wt H$ on $\C P^n$ satisfying the
following conditions:
\roster
\item"(1)" $\wt H|_{T\R P^n}=H$;
\item"(2)" $\wt H(X,IY)=0$ for any $X,Y\in T_p\R P^n$,
$p\in\R P^n$;
\item"(3)" $\wt H$ is $K$-invariant.
\endroster
\endproclaim

By the lemma above the riemannian metric $g'$
extends to a hermitian metric $g$
on $\C P^n$. Let $F'\in\Cal F'$ be an
arbitrary element. By using the bundle
isomorphism $T\R P^n\to T^*\R P^n$ induced from $g'$,
$F'$ is translated to a
symmetric 2-form $F'_{\flat}$. Extending it to a hermitian
form $F_{\flat}$
on $\C P^n$, and again translating with respect
to $g$, we obtain a section $F$ of
$S^2T\C P^n$. Let $\Cal F$ be the vector
space (over $\R$) spanned by such $F$.
Then direct calculations show that $g$ is a
K\"ahler metric, and with this metric
$(\C P^n,\Cal F)$ becomes a K\"ahler-Liouville manifold of
type (A) that satisfies
$\#\A=1$. This completes the proof of Theorem 7.2.
\qed
\enddemo

\specialhead 8. Existence theorem
\endspecialhead
Let $(M, \Cal F)$ be a compact K\"ahler-Liouville
manifold of type (A). For each
$\alpha\in\A$, we define a K\"ahler-Liouville
manifold $(M_{\alpha},\Cal F_{\alpha})$ as follows: Define the closed
subset
$\A_{\alpha}$ of $\A$ by
$$\A_{\alpha}=\{\beta\in\A\ |\ \alpha\pc\beta\}.$$
Let $(M'', \Cal F'')$ be the K\"ahler-Liouville manifold
that is the base space of
the fibre bundle determined by the open subset
$\A-\A_{\alpha}$.
Next, regard $(M'', \Cal F'')$ as the total
space, and let
$(M_{\alpha},\Cal F_{\alpha})$
be the K\"ahler-Liouville manifold that is the typical
fibre of the fibre bundle
determined by the open subset $\{\alpha\}$ of $\A_{\alpha}$.

If $|\alpha|\ge 2$, then $(M_{\alpha},\Cal F_{\alpha})$ is of
type (A), and it
possesses the structure of toric variety that is
given by the structure of
K\"ahler-Liouville manifold. In case $|\alpha|=1$, we also regard
$M_{\alpha}$
as a toric variety, whose structure is inherited
from that of $M''$. So, in any case
$M_{\alpha}$ is isomorphic to $\C P^{|\alpha|}$ as toric
variety.

Let $N$ be a 1-dimensional compact
K\"ahler manifold which is also a toric variety
such that the associated
$U(1)$-action preserves the metric. We shall simply call
it a compact toric
K\"ahler manifold (of dimension 1). To such a
manifold we assign positive constants
$d_*$, $l$, and a function $h$ on $\R/l\Z$
as follows: Let $Y$ be an infinitesimal
generator of the $U(1)$-action so that the least
period of $\exp(sY)$ is $2\pi$.
The set of zeros of $Y$ consists of
two points, say $q_0$ and $q_1$. We may assume
that the endomorphism $\nabla Y$ of $T_{q_1}N$ is
equal to the complex structure $I$
(then it is equal to $-I$ at $q_0$).
Let $l/2$ be the distance
between these two points. Then a minimal geodesic
$\gamma(t)$ from $q_0$
to $q_1$ extends to a closed geodesic of
least period $l$. Since the 1-form
$i_Y\omega$ is closed ($\omega$ is the K\"ahler form),
there is a unique
function $\wt h$ on $M$ such that
$$i_Y\omega=-d\wt h,\qquad \quad \wt h(q_1)=0.$$
Put $d_*=\wt h(q_0)^{-1}$ and $h(t)=d_*\wt h(\gamma(t))$.

The following lemma is immediate.

\proclaim{Lemma 8.1} $(d_*,l,h)$ defined above possesses the following
properties:
\roster
\item"(1)" $h(-t)=h(t)$ for any $t$;
\item"(2)" $h(0)=1,\ h(l/2)=0$;
\item"(3)" $h'(t)<0$ if $0< t< l/2$;
\item"(4)" $-h''(0)=h''(l/2)=d_*$.
\endroster
\endproclaim

Let $\Cal C_1$ be the set of $(d_*,l,h)$
such that $d_*$ and $l$ are positive
constants and $h$ is a $C^{\infty}$ function on
$\R/l\Z$ satisfying the
conditions $(1),\dots,(4)$ in Lemma 8.1. We say that
two elements $(d_*,l,h)$
and $(\wt d_*,\wt l,\wt h)$ are equivalent if
$d_*=\wt d_*$, $l=\wt l$, and either
$h(t)=\wt h(t)$ or $h(t)=1-\wt h(l/2-t)$. For consistency with
the definition of
$\Cal C_n$, we shall also write $(\{1,0\},d_*,l,h)$ instead
of $(d_*,l,h)$.

\proclaim{Lemma 8.2} The assignment above gives the one-to-one
correspondence
between the set of the isomorphism classes of
1-dimensional compact toric
K\"ahler manifolds and the set of the equivalence
classes of elements of
$\Cal C_1$.
\endproclaim

The proof is easy.
Now, we state the main theorem in this
section, which will imply the existence
of compact K\"ahler-Liouville manifold of type (A) whose
structure of
toric variety is prescribed. Let $M$ be a
toric variety of KL-A type, and let
$\A$ be the associated partially ordered set. Let
$m_{\alpha,\nu}$
($\alpha\in\A$, not minimal, $0\le\nu\le |\p(\alpha)|$) be numbers satisfying
Proposition 4.21 with which $M$ is defined. Let
$c_{\alpha,\nu}$ $(0\le\nu\le
|\alpha|,\ \alpha\in\A)$, $e_{\beta\alpha}$
$(\beta\pcc\alpha)$, and $d_{\alpha}$
$(\alpha\in\A)$ be constants that satisfy the conditions (4.1),
(4.2), (4.8), and (4.9).
In this case we say that the constants
$\{c_{\alpha,\nu},e_{\beta\alpha},
d_{\alpha}\}$ are {\it compatible} with the toric variety
$M$ (of KL-A type).
Note that the compatibility has no meaning in
case $\#\A=1$.

\remark{Remark} 1. $M$ determines only the differences
$m_{\alpha,\nu}-m_{\alpha,0}$ for $\alpha$ such that $\p(\alpha)$ is minimal.
Hence for such $\alpha$ one can choose $m_{\alpha,0}$
arbitrary so that they
satisfy Proposition 4.21 (5).\par
2. $\{m_{\alpha,\nu}\}$ just determine every $e_{\beta\alpha}$, every ratio
$d_{\p(\alpha)}/d_{\alpha}$, and $\{c_{\alpha,\nu}\}$ for every non-maximal
$\alpha$. Hence one can choose $d_{\alpha}>0$ arbitrary for
minimal $\alpha$,
and also $c_{\alpha,\nu}$ arbitrary for maximal $\alpha$ so
that they satisfy
(4.1).
\endremark

\proclaim{Theorem 8.3} Let $M$ be a compact toric
variety of KL-A type, and let
$\A$ be the associated partially ordered set. Let
$\{c_{\alpha,\nu},e_{\beta\alpha},
d_{\alpha}\}$ be constants compatible with $M$. For each
$\alpha\in\A$, choose
$l_{\alpha}>0$ and $h_{\alpha}\in C^{\infty}(\R/l_{\alpha}\Z)$ so that
$(\{c_{\alpha,\nu}\}, |d_{\alpha}|,l_{\alpha},h_{\alpha})
\in\Cal C_{|\alpha|}$.
Then there is a unique structure of K\"ahler-Liouville
manifold
$(M,\Cal F)$ of type (A) over the toric
variety $M$ possessing the following properties:
\roster
\item"(1)" The associated structure of toric variety is
identical with the
given one;
\item"(2)" the fundamental constants, the conjunction constants, and
the scaling
constants are equal to $\{c_{\alpha,\nu}\}$, $\{e_{\beta\alpha}\}$, and
$\{d_{\alpha}\}$ respectively;
\item"(3)" for each $\alpha\in\A$, the induced K\"ahler-Liouville manifold
(the toric K\"ahler manifold if $|\alpha|=1$) $(M_{\alpha},\Cal F_{\alpha})$
corresponds to the equivalence class represented by the
given element
$$(\{c_{\alpha,\nu}\},|d_{\alpha}|, l_{\alpha}, h_{\alpha})
\in\Cal C_{|\alpha|}.$$
\endroster
\endproclaim
\demo{Proof} We prove this theorem by induction on
$\#\A$. The case $\#\A=1$
follows from Theorem 7.2. Let $k\ge 2$, and
assume that the theorem is true
for the case where the number of elements
of the associated partially ordered
set is less than $k$.

Now, let $M$ and $\A$ be as above,
and suppose $\#\A=k$. We may assume that
$\A$ is connected. Let $\alpha_0\in\A$ be the minimal
element, and put
$\A'=\{\alpha_0\}$, $\A''=\A-\A'$. As before,
let $\A''=\cup_{s=1}^r \A''_s$ be
the
decomposition into connected components, and $\alpha_s$ the minimal
elements
of $\A''_s$. Let $M'$ and $M''$ be the
associated toric varieties. As is easily seen,
the constants $c_{\alpha,\nu}$ $(\alpha\in\A'')$,
$e_{\beta\alpha}$ $(\beta\in\A'',\beta\pcc\alpha)$,
$\epsilon(\alpha_s)d_{\alpha}$
$(\alpha\in\A''_s, 1\le s\le r)$ are compatible with the
toric variety $M''$ of KL-A
type, where $\epsilon(\alpha_s)$ is the sign of $d_{\alpha_s}$.
So, by induction
assumption we obtain a unique structure of K\"ahler-Liouville
manifold $(M'',\Cal F'')$
over the toric variety $M''$ possessing the properties
stated in the theorem.

Also, by Theorem 7.2 and Lemma 8.2 there
is a unique structure of K\"ahler-Liouville
manifold (or toric K\"ahler manifold if $|\alpha_0|=1$) $(M',\Cal
F')$ over the toric
variety $M'$ corresponding to the element
$$(\{c_{\alpha_0,\nu}\},d_{\alpha_0},l_{\alpha_0},h_{\alpha_0})\in
\Cal C_{|\alpha_0|}.$$
Then, by Theorem 6.11 we obtain a structure
of K\"ahler-Liouville manifold $(M,\Cal F)$
over the toric variety $M$ such that $(M',\Cal
F')$ and $(M'',\Cal F'')$ are isomorphic
to the ones induced from $(M,\Cal F)$. It
is clear that $(M,\Cal F)$ possesses the
properties (1) and (2). (3) follows from the
fact that the K\"ahler-Liouville manifold
(or the toric K\"ahler manifold)
$(M_{\alpha},\Cal F_{\alpha})$ $(\alpha\in\A'')$
induced from $(M,\Cal F)$ is isomorphic to the
one induced from $(M'',\Cal F'')$.

This fact also proves the uniqueness of $(M,\Cal
F)$. In fact, let $(\wt M,\wt{\Cal F})$
be another K\"ahler-Liouville manifold possessing the properties stated
in the theorem,
and let $(\wt M',\wt{\Cal F'})$ and $(\wt M'',\wt{\Cal
F''})$ be the induced ones. Then
the fact mentioned above and the induction assumption
indicate that
$(\wt M'',\wt{\Cal F''})$ is isomorphic to $(M'',\Cal F'')$,
and Theorem 7.2 and
Lemma 8.2 indicate that $(\wt M',\wt{\Cal F'})$
is isomorphic to $(M',\Cal F')$. Hence by Theorem
6.11, $(\wt M,\wt{\Cal F})$ is
isomorphic to $(M,\Cal F)$. This completes the proof.
\qed
\enddemo

\enddocument